\documentclass[11pt]{article}
\usepackage{amsmath, amsfonts}

\numberwithin{equation}{section}

\def\B1{B_{1/2}}
\def\Box{\hfill\rule{2.5mm}{2.5mm}}

\def\F{{\cal F}}
\def\H{{\cal H}}

\def\L{{\cal L}}
\def\N{{\mathbb {N}}}

\def\R{{\mathbb {R}}}
\def\RR{{\cal R}}
\def\SS{{\cal S}}

\def\ap{{\alpha_1^i}}

\def\b{{\frac 2{p-1}}}
\def\bc{{\hat{\cal T}}}

\def\bt{{\bar{\cal T}}}

\def\build#1_#2^#3{\mathrel{
\mathop{\kern 0pt#1}\limits_{#2}^{#3}}}

\def\d{\displaystyle}

\def\h1{\mathop{\rm H^1_{\rm loc,\rm u}}}
\def\iint{\displaystyle\int_{-1}^1}

\def\l2{\mathop{\rm L^2_{\rm loc,\rm u}}}

\def\p{{\bar p}}
\def\pd{\partial_d}
\def\ps{\partial_s}
\def\py{\partial_y}
\def\q{{\bf q}}

\def\r{{\bf r}}
\def\sgn{\mathop{\rm sgn}}

\def\ts{{\tilde s}}

\def\tt{{\tilde t}}

\newcommand{\aref}[1]{(\ref{#1})}
\newcommand{\un}[2]{1_{\{#1<y<#2\}}}
\newcommand{\vc}[2]{
\left(
\begin{array}{l}
#1\\
#2
\end{array}
\right)
}

\title{\bf Existence and classification of characteristic points at blow-up for a semilinear wave equation in one space dimension
\footnote{Both authors are supported by a grant from the french Agence Nationale de la Recherche, project ONDENONLIN, reference ANR-06-BLAN-0185.} }
\author{Frank Merle\\
{\it \small Universit\'e de Cergy Pontoise, IHES and CNRS}\\
Hatem Zaag\\
{\it \small CNRS LAGA Universit\'e Paris 13}}
\date{September 23, 2009}

\newtheorem{cor}{Corollary}[section]

\newtheorem{cl}[cor]{Claim}
\newtheorem{lem}[cor]{Lemma}
\newtheorem{prop}[cor]{Proposition}

\newtheorem{propo}{Proposition}
\newtheorem{theor}[propo]{Theorem}
\newtheorem{coro}[propo]{Corollary}

\begin{document}

\maketitle

{\small {\bf Abstract}: We consider the semilinear wave equation with power nonlinearity in one space dimension. We first show the existence of a blow-up solution with a characteristic point. Then, we consider an arbitrary blow-up solution $u(x,t)$, the graph $x\mapsto T(x)$ of its blow-up points and $\SS\subset \R$ the set of all characteristic points and show that $\SS$ has an empty interior.
 Finally, given $x_0\in \SS$, we show that in selfsimilar variables, the solution decomposes into a decoupled sum of (at least two) solitons, with alternate signs and that $T(x)$ forms a corner of angle $\frac \pi 2$ at $x_0$.}

\medskip


\medskip


\section{Introduction}
\subsection{Known results and the case of non characteristic points}

We consider the one dimensional semilinear wave equation
\begin{equation}\label{equ}
\left\{
\begin{array}{l}
\partial^2_{t t} u =\partial^2_{xx} u+|u|^{p-1}u,\\
u(0)=u_0\mbox{ and }u_t(0)=u_1,
\end{array}
\right.
\end{equation}
where $u(t):x\in\R \rightarrow u(x,t)\in\R$, $p>1$, $u_0\in \rm H^1_{\rm loc,u}$
and $u_1\in \rm L^2_{\rm loc,u}$ with
 $\|v\|_{\rm L^2_{\rm loc,u}}^2=\d\sup_{a\in \R}\int_{|x-a|<1}|v(x)|^2dx$ and $\|v\|_{\rm H^1_{\rm loc,u}}^2 = \|v\|_{\rm L^2_{\rm loc,u}}^2+\|\nabla v\|_{\rm L^2_{\rm loc,u}}^2$.\\
We will also consider the following equation for $p>1$,
\begin{equation}\label{equ'}
\left\{
\begin{array}{l}
\partial^2_{t t} u =\partial^2_{xx} u+|u|^p,\\
u(0)=u_0\mbox{ and }u_t(0)=u_1.
\end{array}
\right.
\end{equation}
The Cauchy problem for equations \aref{equ} and \aref{equ'} in the space ${\rm H}^1_{\rm loc,u}\times {\rm L}^2_{\rm loc, u}$ follows from the finite speed of propagation and the wellposedness in
${\rm H}^1\times {\rm L}^2$ (see Ginibre, Soffer and Velo \cite{GSVjfa92}). The existence of blow-up solutions for equation \aref{equ} follows from Levine \cite{Ltams74}.
More blow-up results can be found in
Caffarelli and Friedman \cite{CFtams86}, \cite{CFarma85},
Alinhac \cite{Apndeta95}, \cite{Afle02} and Kichenassamy and Litman \cite{KL1cpde93}, \cite{KL2cpde93}.

\bigskip

If $u$ is a blow-up solution of \aref{equ}, we define (see for example Alinhac \cite{Apndeta95}) a 1-Lipschitz curve $\Gamma=\{(x,T(x))\}$ 
such that $u$ cannot be extended beyond the set called the maximal influence domain of $u$:
\begin{equation}\label{defdu}
D=\{(x,t)\;|\; t< T(x)\}.
\end{equation}
$\bar T=\inf_{x\in \R}T(x)$ and $\Gamma$ are called the blow-up time and the blow-up graph of $u$. An important notion for the blow-up graph is the notion of characteristic point (even though the existence of characteristic points remained unknown before this paper). 
A point $x_0$ is a non characteristic point (or a regular point)
if
\begin{equation}\label{nonchar}
\mbox{there are }\delta_0\in(0,1)\mbox{ and }t_0<T(x_0)\mbox{ such that }
u\;\;\mbox{is defined on }{\cal C}_{x_0, T(x_0), \delta_0}\cap \{t\ge t_0\}
\end{equation}
where 
${\cal C}_{\bar x, \bar t, \bar \delta}=\{(x,t)\;|\; t< \bar t-\bar \delta|x-\bar x|\}$.
We denote by $\RR$ (resp. $\SS$) the set of non characteristic (resp. characteristic) points.

\bigskip

Following our earlier work (\cite{MZajm03}-\cite{MZjfa07}), we aim at describing the blow-up behavior for {\it any} blow-up solution, especially $\Gamma$ and the solution near $\Gamma$. 

\medskip


\bigskip

Given some $(x_0,T_0)$ such that $0< T_0\le T(x_0)$, a natural tool is to introduce the following self-similar 
change of variables:
\begin{equation}\label{defw}
w_{x_0,T_0}(y,s)=(T_0-t)^{\frac 2{p-1}}u(x,t),\;\;y=\frac{x-x_0}{T_0-t},\;\;
s=-\log(T_0-t).
\end{equation}
If $T_0=T(x_0)$, then we simply write $w_{x_0}$ instead of $w_{x_0, T(x_0)}$. 
The function $w=w_{x_0,T_0}$ satisfies the 
following equation for all $y\in B=B(0,1)$ and $s\ge -\log T_0$:
\begin{equation}\label{eqw}
\partial^2_{ss}w= \L w-\frac{2(p+1)}{(p-1)^2}w+|w|^{p-1}w
-\frac{p+3}{p-1}\partial_sw-2y\partial^2_{y,s} w
\end{equation} 
\begin{equation}\label{defro}
\mbox{where }\L w = \frac 1\rho \py \left(\rho(1-y^2) \py w\right)\mbox{ and }
\rho(y)=(1-y^2)^{\frac 2{p-1}}.
\end{equation}
The Lyapunov functional for equation \aref{eqw}
\begin{equation}\label{defenergy}
E(w(s))= \iint \left(\frac 12 \left(\partial_s w\right)^2 + \frac 12  \left(\partial_y w\right)^2 (1-y^2)+\frac{(p+1)}{(p-1)^2}w^2 - \frac 1{p+1} |w|^{p+1}\right)\rho dy
\end{equation}
is defined for $(w,\partial_s w) \in \H$ where  
\begin{equation}\label{defnh0}
\H = \left\{q 
\;\;|\;\;\|q\|_{\H}^2\equiv \int_{-1}^1 \left(q_1^2+\left(q_1'\right)^2  (1-y^2)+q_2^2\right)\rho dy<+\infty\right\}.
\end{equation}
We will note
\begin{equation}\label{h0}
\H_0 = \{ r \in H^1_{loc}(-1,1)\;|\; \|r\|_{\H_0}^2\equiv\int_{-1}^1 \left( r'^2 (1-y^2) + r^2\right) \rho dy < + \infty\}.
\end{equation}
We also introduce for all $|d|<1$ the following stationary solutions of \aref{eqw} defined by 
\begin{equation}\label{defkd}
\kappa(d,y)=\kappa_0 \frac{(1-d^2)^{\frac 1{p-1}}}{(1+dy)^{\frac 2{p-1}}}\mbox{ where }\kappa_0 = \left(\frac{2(p+1)}{(p-1)^2}\right)^{\frac 1{p-1}} \mbox{ and }|y|<1.
\end{equation}
When $x_0$ is non characteristic, we have a good understanding of the solution's behavior near $x_0$. More precisely, we established the following results in \cite{MZjfa07} and \cite{MZcmp08}:

\medskip

{\bf (Blow-up behavior for $x_0\in \RR$}, {\it see Corollary 4 in \cite{MZjfa07},  Theorem 1 (and the following remark) and Lemma 2.2 in \cite{MZcmp08})}.

{\it (i)}
\label{threg}
{\it 
 The set of non characteristic points $\RR$ is non empty and open.}

{\it (ii)} {\bf (Selfsimilar blow-up profile for $x_0\in \RR$)}
{\it 
There exist positive $\mu_0$ and $C_0$ such that if $x_0 \in \RR$, then there exist $\delta_0(x_0)>0$, $d(x_0)\in (-1, 1)$, $|\theta(x_0)|=1$, $s_0(x_0)\ge - \log T(x_0)$ such that for all $s\ge s_0$:
\begin{eqnarray}
&&\left\|\vc{w_{x_0}(s)}{\partial_s w_{x_0}(s)}-\theta(x_0)\vc{\kappa(d(x_0), .)}{0}\right\|_{\H}\le C_0 e^{-\mu_0(s-s_0)},\label{profile}\\
&&\left\|\vc{w_{x_0}(s)}{\partial_s w_{x_0}(s)}-\theta(x_0)\vc{\kappa(d(x_0), .)}{0}\right\|_{H^1\times L^2(|y|<1+\delta_0)}\to 0 \mbox{ as }s\to \infty.\label{profilelebesgue}
\end{eqnarray}
Moreover, $E(w_{x_0}(s)) \to E(\kappa_0)$ as $s\to \infty$.

(iii)  The function $T(x)$ is $C^1$ on $\RR$ and for all $x_0\in \RR$, $T'(x_0)=d(x_0)\in (-1,1)$. Moreover, $\theta(x_0)$ is constant on connected components of $\RR$.
}

\subsection{Existence of characteristic points}
For characteristic points, the only available result about existence or non existence is due to Caffarelli and Friedman \cite{CFtams86} and \cite{CFarma85} who proved (using the maximum principle) the non existence of characteristic points for equation \aref{equ}:\\
- under conditions on initial data that ensure that for all $x\in\R$ and $t\ge 0$, $u\ge 0$ and $\partial_t u \ge (1+\delta_0)|\partial_x u|$ for some $\delta_0>0$,\\
- for $p\ge 3$ with $u_0\ge 0$, $u_1\ge 0$ and $(u_0,u_1)\in C^4\times C^3(\R)$.

\medskip

From this example, it was generally conjectured by most people that there were no blow-up solutions for equation \aref{equ} with characteristic points: {\it for all $(u_0,u_1)$ which lead to blow-up, $\RR=\R$.}\\
Our first result is to disprove this fact. 
Existence of characteristic points is seen as a consequence of two facts:\\
- on the one hand, the study of the blow-up profile at a non characteristic point,\\
- on the other hand, connectedness arguments related to the sign of the blow-up profile.

\medskip

To state our results, let us consider $u(x,t)$ a blow-up solution of equation \aref{equ} (take for example initial data $(u_0, u_1) \in H^1\times L^2(\R)$ satisfying
 $$\int_{\R}\left(\frac 12 |\partial_x u_0|^2+ \frac 12 u_1^2 -\frac 1{p+1}|u_0|^{p+1}\right)dx <0,$$ which gives blow-up by Levine \cite{Ltams74}).
The first result follows from the study near a regular point, that ensures the existence of an explicit signed profile.
\begin{propo}
\label{propelem}
If the initial data $(u_0,u_1)$ is odd and $u(x,t)$ blows up in finite time, then $0\in\SS$.
\end{propo}
The second one follows from the continuity of the profile on the connected components of $\RR$ (see Theorem 1 in \cite{MZcmp08}).
\begin{theor}[Existence and generic stability of characteristic points]\label{thexis}$ $\\
(i) {\bf (Existence)} Let $a_1<a_2$ be two non characteristic points such that
\[
w_{a_i}(s) \to \theta(a_i)\kappa(d_{a_i},\cdot)\mbox{ as }s\to \infty\mbox{ with }\theta(a_1)\theta(a_2)=-1
\]
for some $d_{a_i}$ in $(-1,1)$, in the sense \aref{profile}. Then, there exists a characteristic point $c\in (a_1,a_2)$.\\
(ii) {\bf (Stability)} There exists $\epsilon_0>0$ such that if $\|(\tilde u_0,\tilde u_1)- (u_0, u_1)\|_{\h1\times \l2}\le \epsilon_0$, then, $\tilde u(x,t)$ the solution of equation \aref{equ} with initial data $(\tilde u_0,\tilde u_1)$ blows up and has a characteristic point $\tilde c\in [a_1,a_2]$. 
\end{theor}
{\bf Remark}: It is enough to take $(u_0,u_1)$  with large plateaus of opposite signs to guarantee that $u(x,t)$ blows up satisfying the hypotheses of this theorem.

\medskip

Since a solution in one space dimension is also a solution in higher dimensions, we get from the finite speed of propagation the following existence result in $N$ dimensions:
\begin{coro}\label{corN}{\bf (Existence of characteristic points in higher dimensions)} Consider $\tilde u(x_1,t)$ a blow-up solution of \aref{equ} in one space dimension with a characteristic point. Then, for $R$ large enough, initial data $(u_0, u_1)$ such that $u_i(x)=\tilde u_i(x_1)$ for $|x|<R$, the solution $u(x,t)$ of equation \aref{equ} with initial data $(u_0, u_1)$ blows up and has a characteristic point.
\end{coro}


\subsection{Non existence results for characteristic points}
In this section, we give sufficient conditions under which no characteristic point can occur. Our analysis in fact relates the fact that $x_0$ is a characteristic point to sign changes of the solution in a neighborhood of $(x_0, T(x_0))$. We claim the following: 
\begin{theor}\label{thnonexis} Consider $u(x,t)$ a blow-up solution of \aref{equ} such that $u(x,t)\ge 0$ for all $x\in (a_0,b_0)$ and $t_0\le t< T(x)$ for some real $a_0$, $b_0$ and $t_0\ge 0$. Then, $(a_0, b_0)\subset \RR$.
\end{theor}
{\bf Remark}: This result can be seen as a generalization of the result of Caffarelli and Friedman, with no restriction on initial data. Indeed, from our result, taking nonnegative initial data suffices to exclude the occurrence of characteristic points.

\bigskip

Considering the equation \aref{equ'}, we get the following twin result of Theorem \ref{thnonexis}:

\medskip

\noindent {\bf Theorem \ref{thnonexis}'} {\it The set of characteristic points is empty for any blow-up solution of equation \aref{equ'}}.

\medskip

\subsection{Shape of the blow-up set near characteristic points and properties of $\SS$}
We have the following proposition:
\begin{propo}\label{thfini} {\bf ($\SS$ has an empty interior)} Consider $u(x,t)$ a blow-up solution of \aref{equ}. The set of characteristic points $\SS$ has an empty interior.
\end{propo}
{\bf Remark}: This implies in particular that $\SS$ has zero measure.
Direct arguments give no more than the fact that $\SS\neq \R$ (a point $x_0$ such that $T(x_0)$ is the blow-up time is non characteristic; see in page \pageref{threg} point (i) of the result of \cite{MZjfa07} and \cite{MZcmp08}).
The proof of Proposition \ref{thfini} uses 
the 
description 
of the solution in the $w$ variable as 
 a non trivial decoupled sum of (at least 2) $\pm \kappa(d_i(s))$ (see Theorem \ref{new} below). 

\medskip

Now, we have the following theorem, which is the main result of our analysis, where for a given $x_0\in \SS$, we are able to give the precise behavior of the solution near $(x_0, T(x_0))$: 
\begin{theor}\label{new}{\bf (Description of the behavior of $w_{x_0}$ where $x_0$ is characteristic)} Consider $u(x,t)$ a blow-up solution of \aref{equ} and $x_0\in \SS$. Then,  it holds that
\begin{equation}\label{cprofile00}
\left\|\vc{w_{x_0}(s)}{\ps w_{x_0}(s)} - \vc{\d\sum_{i=1}^{k(x_0)} e^*_i\kappa(d_i(s),\cdot)}0\right\|_{\H} \to 0\mbox{ and }E(w_{x_0}(s))\to k(x_0)E(\kappa_0)
\end{equation}
as $s\to \infty$, for some 
\begin{equation*}
k(x_0)\ge 2,\;\;e^*_i=e^*_1(-1)^{i+1}
\end{equation*}
 and continuous $d_i(s)=-\tanh \zeta_i(s)\in (-1,1)$ for $i=1,...,k(x_0)$. Moreover, for some $C_0>0$, for all $i=1,...,k(x_0)$ and $s$ large enough,
\begin{equation}\label{equid}
\left|\zeta_i(s)-\left(i-\frac{(k(x_0)+1)}2\right)\frac{(p-1)}2\log s\right|\le C_0. 
\end{equation}
\end{theor}
{\bf Remark}: In \cite{MZjfa07}, we proved a much weaker version of this result, with \aref{cprofile00} valid just with $k(x_0)\ge 0$ and no information on the signs of $e_i^*$, $\zeta_1(s)$ and $\zeta_k(s)$. Note that eliminating the case 
\[
k(x_0)=0
\]
 is the most difficult part in our analysis. In some sense, we put in relation the notion of characteristic point at $x_0$ and the notion of decomposition of $w_{x_0}$ in a decoupled sum of (at least 2) $\pm \kappa(d_i(s))$. This result can be seen as a result of decomposition up to dispersion into sum of decoupling solitons in dispersive problems. According to the value of $k(x_0)$, this sum appears to have a {\it multipole} nature (dipole if $k(x_0)=2$, tripole if $k(x_0)=3$,....). In some sense, it says that the space $\H$ is a critical space to measure dispersion and blow-up in the cone (i.e. if $E(w_{x_0, T_0})<E(\kappa_0)$, then:\\
- $T(x_0)>T_0$ and the solution can be extended in a strictly greater cone;\\
- $(w_{x_0, T_0}, \ps w_{x_0, T_0})$ converges as $s\to \infty$ to zero in $H^1\times L^2$ of the cone slice).\\
{\bf Remark}: Estimate \aref{equid} comes from the fact that the centers' positions $\zeta_i(s)$ satisfy the finite-dimensional Toda lattice system given in \aref{eqz} (which is an integrable system).
From \aref{equid}, we see that the distance between two solitons $\zeta_{i+1}(s)-\zeta_i(s)\sim \frac{(p-1)}2 \log s = \frac{(p-1)}2 \log|\log(T-t)|$ as $t\to T(x_0)$.\\
If $r=E\left(\frac {k(x_0)}2\right)$
then we see from \aref{equid} that $r$ solitons go to $-\infty$ as $s\to \infty$ (for $i=1,..,r$) and $r$ solitons go to $+\infty$ (for $i=k(x_0)+1-r,...,k(x_0)$). If $k(x_0)=2r+1$, then the central soliton (for $i=r+1=\frac{k(x_0)+1}2$) stays bounded for $s$ large.

\bigskip

Extending the definition of $k(x_0)$ defined for $x_0\in \SS$ in Theorem \ref{new} by setting 
\[
k(x_0)=1\mbox{ for all }x_0\in \RR
\]
and using the monotonicity of the Lyapunov functional $E(w)$, we have the following consequence from the blow-up behavior in the characteristic case (Theorem \ref{new}) and in the non characteristic case (the result of \cite{MZjfa07} cited here in page \pageref{threg}):
\begin{coro}[A criterion for non characteristic points]\label{corcriterion}$ $\\
(i) For all $x_0\in \R$ and $s_0\ge -\log T(x_0)$, we have
\[
E(w_{x_0}(s_0))\ge k(x_0)E(\kappa_0).
\]
(ii) If for some $x_0$ and $s_0\ge -\log T(x_0)$, we have 
\[
E(w_{x_0}(s_0))<2 E(\kappa_0),
\]
then $x_0\in \RR$.
\end{coro}
We also have the following consequences in the original variables: 
\begin{propo}\label{pdes}{\bf (Description of $T(x)$ for $x$ near $x_0$)}\\
(i) If $x_0\in \SS$ and $0<|x-x_0|\le \delta_0$, then 
\begin{equation}\label{chapeau0}
0<T(x)- T(x_0)+|x-x_0| \le \frac{C_0|x-x_0|}{|\log(x-x_0)|^{\frac{(k(x_0)-1)(p-1)}2}}
\end{equation} 
for some $\delta_0>0$ and $C_0$, where $k(x_0)$ is defined in Theorem \ref{new}.\\
(ii) If $x_0\in \SS$, then $T(x)$ is right and left differentiable at $x_0$, with
\[
T'_l(x_0)=1\mbox{ and }T'_r(x_0)=-1.
\]
(iii) For all $t\in [T(x_0)-\tau_0, T(x_0))$ for some $\tau_0>0$, there exist $z_1(t)<...<z_k(t)$ continuous in $t$ such that
\[
e_1^*(-1)^{i+1}u(z_i(t),t)>0
\]
and $z_i(t) \to x_0$ as $t\to T(x_0)$.
\end{propo}
{\bf Remark}: From (iii), we have the existence of zero lines $x_1(t)< ... < x_{k-1}(t)$ (not necessarily continuous in $t$) such that $u(x_i(t),t)=0$ and $x_i(t)\to x_0$ as $t\to T(x_0)$.\\
{\bf Remark}: In a forthcoming paper \cite{MZisol09}, we improve \aref{chapeau0} by finding a lower bound of the same type as the upper bound.

\bigskip

The paper is organized as follows. Section \ref{section2} is devoted to the proofs of Proposition \ref{propelem} and Theorem \ref{thexis} (note that Corollary \ref{corN} follows straightforwardly from Theorem \ref{thexis} and the finite speed of propagation). In Section \ref{section3}, we consider a characteristic point and study the equation in selfsimilar variables. As for Section \ref{sec2}, it is devoted to the proof of Theorems \ref{new}, \ref{thnonexis} and \ref{thnonexis}', as well as Propositions \ref{thfini} and \ref{pdes} (note that Corollary \ref{corcriterion} is a direct consequence of Theorem \ref{new} and the result of \cite{MZjfa07} cited here in page \pageref{threg}). 

\section{Existence and stability of characteristic points}\label{section2}


Here in this section, we consider $u(x,t)$ a blow-up solution of equation \aref{equ}.
As mentioned in the introduction, we prove in this section the existence of characteristic points (Proposition \ref{propelem} and Theorem \ref{thexis}).

\bigskip

{\it Proof of Proposition \ref{propelem}}: Assuming that $(u_0, u_1)$ is odd, we would like to prove that $0\in \SS$. Arguing by contradiction, we assume that $0\in\RR$.\\
On the one hand, using the result of \cite{MZjfa07} stated in \aref{profilelebesgue}, we see that for some $d(0)\in (-1,1)$,
\[
\|w_0(s)-\kappa(d(0),\cdot)\|_{L^\infty(-1,1)}\le C \|w_0(s)-\kappa(d(0),\cdot)\|_{H^1(-1,1)}\to 0\mbox{ as }s\to \infty.
\]
In particular, 
\begin{equation}\label{wali}
|w_0(0,s)|\to \kappa(d(0),0)>0\mbox{  as }s\to \infty.
\end{equation}
On the other hand, since the initial data is odd, the same holds for the solution, in particular, $u(0,t)=0$ for all $t\in [0,T(0))$, hence $w_0(0,s)=0$ for all $s\ge -\log T(0)$, which contradicts \aref{wali}. This concludes the proof of Proposition \ref{propelem}. \Box\\
{\bf Remark}: We don't need to know that for $x_0\in \RR$, $w_{x_0}$ converges to a particular profile to derive this result. It is enough to know that $w_{x_0}$ approaches the set $\{\theta(x_0) \kappa(d,\cdot)\;|\; |d|<1-\eta\}$ for some $\eta>0$, which is a much weaker result.

\bigskip


We now turn to the proof of Theorem \ref{thexis}. It is a consequence of three results from our earlier work:

\medskip

- the continuity with respect to initial data of the blow-up time at $x_0\in \RR$. 
\begin{prop}[Continuity with respect to initial data at $x_0\in \RR$]\label{propcontid} There exists $A_0>0$ such that $\tilde T(x_0)\to T(x_0)$ as $(\tilde u_0, \tilde u_1)\to (u_0, u_1)$ in $H^1\times L^2(|x|<A_0)$, where $\tilde T(x_0)$ is the blow-up time of $\tilde u(x,t)$ at $x=x_0$, the solution of equation \aref{equ} with initial data $(\tilde u_0, \tilde u_1)$.
\end{prop}
{\it Proof}: This is a direct consequence of the Liouville Theorem and its applications given in \cite{MZcmp08}. See Appendix \ref{appcont} for a sketch of the proof.\Box

\medskip

- the continuity of the blow-up profile on $\RR$ proved in Theorem 1 in \cite{MZcmp08} (in particular, the fact that $\theta(x_0)$ given in \aref{profile} is constant on the connected components of $\RR$).

\medskip

- the following trapping result from \cite{MZjfa07}: 
\begin{prop}{\bf (See Theorem 3 in \cite{MZjfa07} and its proof)}\label{thtrap}
There exists $\epsilon_0>0$ such that if $w\in C([s^*, \infty), \H)$ for some $s^*\in \R$ is a solution of equation \aref{eqw} such that 
\begin{equation}\label{trap}
\forall s\ge s^*,\;\;E(w(s))\ge E(\kappa_0)\mbox{ and }\left\|\vc{w(s^*)}{\partial_s w(s^*)}-\omega^*\vc{\kappa(d^*,\cdot)}{0}\right\|_{\H}\le \epsilon^*
\end{equation}
for some $d^*=-\tanh \xi^*$, $\omega^*=\pm 1$ and $\epsilon^*\in (0, \epsilon_0]$, then there exists $d_\infty=-\tanh \xi_\infty$ such that
\begin{equation}\label{profile0}
\left|\xi_\infty - \xi^*\right|\le C_0 \epsilon^*\mbox{ and }
\left\|\vc{w(s)}{\partial_s w(s)}-\omega^*\vc{\kappa(d_\infty, \cdot)}{0}\right\|_{\H}\to 0.
\end{equation}
\end{prop}
 
Let us use these results to prove Theorem \ref{thexis}.

\bigskip

{\it Proof of Theorem \ref{thexis}}: We consider $a_1<a_2$ two non characteristic points such that $w_{a_i}(s) \to \theta(a_i)\kappa(d_{a_i},\cdot)$ with $\theta(a_1)\theta(a_2)=-1$ for some $d_{a_i}$ in $(-1,1)$, in the sense \aref{profile}. Up to changing $u$ in $-u$, we can assume that $\theta(a_1)=1$ and $\theta(a_2)=-1$. We aim at proving that $(a_1,a_2)\cap \SS\neq \emptyset$ and the stability of such a property with respect to initial data.

\medskip

(i) If we assume by contradiction that $[a_1,a_2]\subset \RR$, then the continuity of $\theta(x_0)$ where $x_0\in [a_1,a_2]$ implies that $\theta(x_0)$ is constant on $[a_1,a_2]$. This is a contradiction, since $\theta(a_1)=1$ and $\theta(a_2)=-1$. 

\medskip

(ii) By hypothesis and estimate \aref{profilelebesgue}, there is $\delta_0>0$ and $s_0\in \R$ such that
\[
\left\|\vc{w_{a_i}(s_0)}{\partial_s w_{a_i}(s_0)}-\theta(a_i)\vc{\kappa(d_{a_i},\cdot)}{0}\right\|_{H^1\times L^2(|y|<1+\delta_0)}\le \frac{\epsilon_0}2
\]
where $\epsilon_0$ is defined in Proposition \ref{thtrap}. From the continuity with respect to initial data for equation \aref{equ} at the fixed time $T(a_i)- e^{-s_0}$, we see there exists $\eta(\epsilon_0)>0$ such that if 
\[
\|(\tilde u_0,\tilde u_1)- (u_0, u_1)\|_{\h1\times \l2}\le \eta,
\]
then $\tilde u(x,t)$ the solution of equation \aref{equ} with initial data $(\tilde u_0,\tilde u_1)$ is such that $\tilde w_{a_i}(y,s_0)$ is defined for all $|y|<1+\delta_0/2$ and 
\[
\left\|\vc{\tilde w_{a_i}(s_0)}{\partial_s \tilde w_{a_i}(s_0)}-\theta(a_i)\vc{\kappa(d_{a_i},\cdot)}{0}\right\|_{H^1\times L^2(|y|<1+\delta_0/2)}\le \frac 34 \epsilon_0,
\]
where $\tilde w_{a_i}$ is the selfsimilar version defined from $\tilde u(x,t)$ by \aref{defw}.\\
From Proposition \ref{propcontid}, we have $\tilde T(a_i)\to T(a_i)$ as $\eta\to 0$, where $\tilde T(a_i)$ is the blow-up time of $\tilde u(t)$ at $a_i$. We then have for $\eta$ small enough,  
\begin{equation}\label{18,5}
\left\|\vc{\tilde w_{a_i}(s_0)}{\partial_s \tilde w_{a_i}(s_0)}-\theta(a_i)\vc{\kappa(d_{a_i},\cdot)}{0}\right\|_{\H}\le \epsilon_0.
\end{equation}
Two cases then arise (by the way, we will prove later in Theorem \ref{new} that  the Lyapunov functional stays above $2 E(\kappa_0)$ at a characteristic point, which means by \aref{18,5} that $a_i$ and $a_2$ are non characteristic points for $\eta$ small enough, but we cannot use Theorem \ref{new} for the moment):\\
- If $a_1$ or $a_2$ is a characteristic point of $\tilde u(t)$, then the proof is finished.\\
- Otherwise, \aref{profile} holds for $\tilde w_{a_i}$ from the fact that the point is non characteristic. Thus, from the monotonicity of $E(\tilde w_{a_i}(s))$, \aref{trap} holds with $\omega^*= \theta(a_i)$. 
Applying Proposition \ref{thtrap}, we see that $\tilde w_{a_i}(s) \to \theta(a_i)\kappa(\tilde d_{a_i},\cdot)$ as $s\to \infty$, for some $\tilde d_{a_i}\in (-1,1)$. Noting that $\theta(a_1)=1$ and $\theta(a_2)=-1$, we apply (i) to get the result. This concludes the proof of Theorem \ref{thexis}.\Box

\section{Refined behavior for $w_{x_0}$ where $x_0$ is characteristic}\label{section3}
In this section, we consider $x_0\in \SS$. We know from \cite{MZcmp08} that
\begin{equation}\label{cprofile0}
\left\|\vc{w_{x_0}(s)}{\ps w_{x_0}(s)} - \d\sum_{i=1}^{k(x_0)} e_i\vc{\kappa(d_i(s),\cdot)}0\right\|_{\H} \to 0\mbox{ as }s\to \infty
\end{equation}
for some $k(x_0)\ge 0$, $e_i=\pm 1$ and continuous $d_i(s)=-\tanh \zeta_i(s)\in (-1,1)$ for $i=1,...,k(x_0)$ with 
\begin{equation}\label{div}
\zeta_1(s)< ...< \zeta_{k(x_0)}(s)\mbox{ and }\zeta_{i+1}(s)-\zeta_i(s)\to \infty\mbox{ for all }i=1,...,k-1.
\end{equation}
Since $w_{x_0}(s)$ is convergent when $k(x_0)\le 1$ (to $0$ when $k(x_0)=0$ and to some $\kappa(d_\infty)$ by Proposition \ref{thtrap}), we focus throughout this section on the case 
\[
k(x_0)\ge 2.
\]
 For simplicity in the notations, we forget the dependence of $w_{x_0}$ and $k(x_0)$ on $x_0$.\\
This section is organized as follows. In Subsection \ref{sub1}, assuming an ODE on the solitons' center, we find their behavior. Then, in Subsection \ref{subseclemode}, we study equation \aref{eqw} around the solitons' sum and derive in Subsection \ref{applemode0} the ODE satisfied by the solitons' center. Finally, we prove in Subsection \ref{sub4} the corner property near characteristic points.

\subsection{Time behavior of the solitons' centers}\label{sub1}
We will prove the following:
\begin{prop}\label{propref}{\bf (Refined behavior of $w_{x_0}$ where $x_0\in \SS$)}Assuming that $k\ge 2$, there exists another set of parameters (still denoted by $\zeta_1(s)$, ... $\zeta_k(s)$) such that \aref{cprofile0} and \aref{div} hold and:\\
(i)  For all $i=1,...,k$, $e_i=(-1)^{i+1}e_1$.\\
(ii) For some $C_0>0$, for all $i=1,...,k$ and $s$ large enough, we have:
\[
\left(i-\frac{(k+1)}2\right)\frac{(p-1)}2\log s -C_0 \le \zeta_i(s) \le \left(i-\frac{(k+1)}2\right)\frac{(p-1)}2\log s +C_0.
\]
\end{prop}
The following finite-dimensional Toda lattice system  
satisfied by $(\zeta_i(s))$ is crucial in our proof: 
\begin{prop}\label{lemode}{\bf (Equations satisfied by the solitons' centers)}
Assuming that $k\ge 2$, there exists another set of parameters (still denoted by $\zeta_1(s)$, ... $\zeta_k(s)$) such that \aref{cprofile0} and \aref{div} hold and for all $i=1,...,k$ and $s$ large enough:
\begin{eqnarray}
\frac 1{c_1} \zeta_i' &=& - e_{i-1}e_i e^{-\frac 2{p-1}(\zeta_i-\zeta_{i-1})}
+ e_i e_{i+1} e^{-\frac 2{p-1}(\zeta_{i+1}-\zeta_i)}+R_i\label{eqz}\\
\mbox{ where }
|R_i| &\le & CJ^{1+\delta_0},\;\; J(s)=\sum_{j=1}^{k-1} e^{-\frac 2{p-1}(\zeta_{j+1}(s)-\zeta_j(s))},\label{defJ}
\end{eqnarray}
 $e_0=e_{k+1}=0$, for some $c_1>0$ and $\delta_0>0$.
\end{prop}
{\it Proof}: See subsection \ref{applemode0}.\Box

\medskip

Let us now give the proof of Proposition \ref{propref}.

\medskip

{\it Proof of Proposition \ref{propref}}:  We proceed in two parts, proving (i) and then (ii).

\bigskip

{\bf Part 1: Proof of (i)}

Given some $s_0\in\R$, we first define for all $s\ge s_0$, $J_0(s)$ and $j_0(s)\in \{1,...,k-1\}$, 
\begin{equation}\label{defj0}
J_0(s)\equiv \max_{i=1,..k-1}\int_{s_0}^s e^{-\beta(\zeta_{i+1}(s')-\zeta_i(s'))}ds'=\int_{s_0}^s e^{-\beta(\zeta_{j_0(s)+1}(s')-\zeta_{j_0(s)}(s'))}ds'
\end{equation}
where $\beta= \frac 2{p-1}$. Then, we claim that
\begin{equation}\label{lim}
J_0(s)\to \infty\mbox{ as }s\to \infty.
\end{equation}
Indeed, we write from \aref{eqz} and \aref{defj0}
$|\zeta_i(s)-\zeta_i(s_0)|\le C\sum_{i=1}^{k-1} \int_{s_0}^s e^{-\beta(\zeta_{i+1}(s')-\zeta_i(s'))}ds'\le C J_0(s)$
 and \aref{div} implies \aref{lim}.\\
Integrating equation \aref{eqz}, this yields as $s\to \infty$ for all $i=1,...,k$:
\begin{eqnarray}
\frac{\zeta_1(s)+...+\zeta_{i}(s)}{i}&=&e_{i}e_{i+1}\frac{c_1}{i}\int_{s_0}^{s} e^{-\beta(\zeta_{i+1}(s')-\zeta_{i}(s'))}ds'+ o\left(J_0(s)\right),\label{30}\\
\frac{\zeta_i(s)+...+\zeta_k(s)}{k-i+1}&=&-e_{i}e_{i-1}\frac{c_1}{k-i+1}\int_{s_0}^{s} e^{-\beta(\zeta_{i}(s')-\zeta_{i-1}(s'))}ds'+ o\left(J_0(s)\right)\label{31}.
\end{eqnarray}
%
%
%
Using \aref{div}, we write for $s$ large,
\begin{eqnarray*}
\frac{\zeta_1(s)+...+\zeta_{j_0(s)}(s)}{j_0(s)}&<& \frac{\zeta_{j_0(s)+1}(s)+...+\zeta_k(s)}{k-j_0(s)},\\
\mbox{if }i< j_0(s),\;\;\frac{\zeta_1(s)+...+\zeta_{i}(s)}{i}&<&\frac{\zeta_1(s)+...+\zeta_{j_0(s)}(s)}{j_0(s)},\\
\mbox{if }i> j_0(s),\;\;\frac{\zeta_{j_0(s)+1}(s)+...+\zeta_{k}(s)}{k-j_0(s)}&<&\frac{\zeta_{i+1}(s)+...+\zeta_{k}(s)}{k-i}.
\end{eqnarray*}
Then, using \aref{30}, \aref{31} and \aref{defj0}, we write for $s$ large,
\begin{eqnarray}
&& e_{j_0(s)}e_{j_0(s)+1}\frac{c_1}{j_0(s)}J_0(s)\le -e_{j_0(s)}e_{j_0(s)+1}\frac{c_1}{k-j_0(s)}J_0(s)+o(J_0(s)),\nonumber\\
&\mbox{if}&i< j_0(s),\nonumber\\
&& e_{i}e_{i+1}\frac{c_1}{i}\int_{s_0}^{s} e^{-\beta(\zeta_{i+1}(s')-\zeta_{i}(s'))}ds'\le e_{j_0(s)}e_{j_0(s)+1}\frac{c_1}{k-j_0(s)}J_0(s)+o(J_0(s)),\label{26,1}\\
&\mbox{if}&i> j_0(s),\nonumber\\
&& -\frac{c_1e_{j_0(s)}e_{j_0(s)+1}}{k-j_0(s)}J_0(s)+o(J_0(s))\le 
-e_{i}e_{i+1}\frac{c_1}{k-i}\int_{s_0}^{s} e^{-\beta(\zeta_{i+1}(s')-\zeta_{i}(s'))}ds'.\label{26,2}
\end{eqnarray}
Therefore, for $s$ large, 
$J_0(s)\left(e_{j_0(s)}e_{j_0(s)+1}(\frac{c_1}{j_0(s)}+\frac{c_1}{k-j_0(s)})+o(1)\right)\le 0$,
hence, 
\[
e_{j_0(s)}e_{j_0(s)+1}=-1.
\]
Then, \aref{26,1} and \aref{26,2} write together with \aref{defj0}
\[
\forall i,\;\;\frac 1{2k} J_0(s) \le -e_{i}e_{i+1}\int_{s_0}^{s} e^{-\beta(\zeta_{i+1}(s')-\zeta_{i}(s'))}ds'\le J_0(s),
\]
which gives for all $i$ and $s$ large,
\begin{equation}\label{equiv}
e_ie_{i+1}=-1\mbox{ and }\frac {J_0(s)}{C_0}\le \int_{s_0}^s e^{-\beta(\zeta_{j+1}(s')-\zeta_j(s'))}ds' \le J_0(s)\to \infty.
\end{equation}
Using a finite induction, we get $e_i=(-1)^{i+1}e_1$. This closes the proof of (i) of Proposition \ref{propref}.

\bigskip

{\bf Part 2: Proof of (ii)}

%
%
Using (i), we rewrite system \aref{eqz} in the following:
\begin{cor}\label{corlemode}{\bf (Equations satisfied by the solitons' centers)}
Assuming that $k\ge 2$, it holds 
that
\begin{equation}\label{eqz1}
\begin{array}{rcllll}
\d\frac 1{c_1} \zeta_1' &=&&
-e^{-\frac 2{p-1}(\zeta_2-\zeta_1)}&+R_1\\
\\
\d\mbox{if }i=2,..,k-1,\;\;\frac 1{c_1} \zeta_i' &=&e^{-\frac 2{p-1}(\zeta_i-\zeta_{i-1})}
&-e^{-\frac 2{p-1}(\zeta_{i+1}-\zeta_i)}&+R_i\\
\\
\d\frac 1{c_1} \zeta_k' &=&e^{-\frac 2{p-1}(\zeta_k-\zeta_{k-1})}
&&+R_k
\end{array}
\end{equation}
with $|R_i|\le C J^{1+\delta_0}$ for $s$ large enough.
\end{cor}
Introducing
\begin{equation}\label{defli}
L_i(s)=\zeta_{i+1}(s)-\zeta_i(s),\mbox{ where }i=1,...,k-1
\end{equation}
and 
\begin{equation}
\sigma(s)= \sum_{j=1}^r \epsilon_0^{j-1 }S_j(s)\mbox{ where }S_j(s) = \sum_{i=j}^{k-j} L_i(s),\;\;\epsilon _0 = \frac 1{1000}\label{defsigma}
\end{equation}
and 
\begin{equation}\label{defr}
r=E\left(\frac k2\right)\in \N^*\mbox{ (note that either }k=2r+1\mbox{ or }k=2r\mbox{)},
\end{equation}
we claim that (ii) follows from the following: 
\begin{lem}\label{lemli}
Assuming that $k\ge 2$, the following holds for $s$ large enough:\\
(i) For some $C_0>0$ and for all $i=1,..,k-1$,
\[
|L_i(s)-L_1(s)|\le C_0.
\]
(ii)
\[
\frac{\epsilon_0^{r-1}}2 \sum_{j=1}^{k-1}e^{-\b L_j(s)}\le \frac{\sigma'(s)}{c_1}\le 3 \sum_{j=1}^{k-1}e^{-\b L_j(s)}.
\]
\end{lem}
Let us first derive (ii) from Lemma \ref{lemli} and then prove Lemma \ref{lemli}.\\
%
%
Using (i) of Lemma \ref{lemli}, we see that \aref{defsigma} and 
(ii) of Lemma \ref{lemli}
 give for some $C_0>0$ and $s$ large enough,
\begin{eqnarray}
L_1(s)-C_0 &\le &\frac {\sigma(s)}{\gamma} \le L_1(s)+C_0\mbox{ with }\gamma =\sum_{j=1}^r\epsilon_0^{j-1}(k+1-2j)>0,\label{nejm}\\
\frac 1{C_0} e^{-\b L_1(s)} &\le& \sigma'(s) \le C_0 e^{-\b L_1(s)}.\nonumber
\end{eqnarray}
Therefore, for $s$ large enough, we have
\[
\frac 1{C_0} e^{-\b\frac{\sigma(s)}\gamma} \le \sigma'(s) \le C_0 e^{-\b\frac{\sigma(s)}\gamma}\mbox{ and }
\frac 1{C_0} \le (e^{\b \frac{\sigma(s)}\gamma})'\le C_0.
\]
Integrating this and using \aref{nejm}, we see that for $s$ large, we have
\begin{eqnarray*}
\frac s{C_0} - C_0&\le& e^{\b\frac{\sigma(s)}\gamma} \le C_0 s +C_0,\\
\log s - C_0 &\le& \b\frac{\sigma(s)}\gamma \le \log s +C_0,\\
\log s - C_0 &\le& \b L_1(s) \le \log s +C_0.
\end{eqnarray*}
Using (i) of Lemma \ref{lemli}, we see that for all $i=1,...,k-1$ and $s$ large enough, we have
\begin{equation}\label{boundli}
\frac{(p-1)}2\log s - C_0 \le L_i(s) \le \frac{(p-1)}2\log s +C_0.
\end{equation}
Therefore, 
we write from Corollary \ref{corlemode}
\[
\frac{(\zeta_1(s)+...+ \zeta_k(s))'}k= \frac{c_1}k\sum_{i=1}^k R_i(s) \le CJ(s)^{1+\delta_0} \le \frac {C_0}{s^{1+\delta_0}}
\] 
for $s$ large enough. 
Hence, 
\begin{equation}\label{dim}
\bar \zeta(s) \equiv \frac{\zeta_1(s)+...+ \zeta_k(s)}k\mbox{ converges to some }\bar \zeta_0\mbox{ as }s\to \infty.
\end{equation}
Now, according to \aref{defr}, we consider two cases: $k=2r$ and $k=2r+1$.

\medskip

{\bf Case $k=2r$}.
Using the definition \aref{defli} of $L_i$ and (i) of Lemma \ref{lemli}, we see that for all $j=1,...,r$,
\begin{equation}\label{133}
-C_0 \le \left(\zeta_{r+j}(s) - \frac{\zeta_r(s)+\zeta_{r+1}(s)}2\right) - \left(\frac{\zeta_r(s)+\zeta_{r+1}(s)}2 - \zeta_{r+1-j}(s)\right)\le C_0.
\end{equation}
Since we have from \aref{dim}, 
\begin{eqnarray*}
&&\sum_{j=1}^r \left(\zeta_{r+j}(s) -  \frac{\zeta_r(s)+\zeta_{r+1}(s)}2\right)- \left(\frac{\zeta_r(s)+\zeta_{r+1}(s)}2 - \zeta_{r+1-j}(s)\right)\\
&=& \sum_{j=1}^r (\zeta_{r+j}(s)+\zeta_{r+1-j}(s))-2r \frac{(\zeta_r(s)+\zeta_{r+1}(s))}2\\
&=& k\left(\bar \zeta(s) - \frac{\zeta_r(s)+\zeta_{r+1}(s)}2\right),
\end{eqnarray*}
we see from \aref{133} that for $s$ large enough, 
\[
-C_0 \le \bar \zeta(s) - \frac{\zeta_r(s)+\zeta_{r+1}(s)}2 \le C_0,
\]
and from \aref{dim}, we see that
\begin{equation}\label{134}
-C_0 \le \frac{\zeta_r(s)+\zeta_{r+1}(s)}2\le C_0.
\end{equation}
Now, since we have from the definition \aref{defli} of $L_i$ and \aref{boundli} for all $i=1,...,k$ and $s$ large enough, 
\[
(i-(r+\frac 12))\frac{(p-1)}2\log s - C_0 \le \zeta_i(s) - \frac{\zeta_r(s)+\zeta_{r+1}(s)}2 \le (i-(r+\frac 12))\frac{(p-1)}2\log s + C_0
\]
and $r+\frac 12= \frac{k+1}2$, the conclusion of (ii)
follows from \aref{134}, when $k=2r$.

\medskip

{\bf Case $k=2r+1$}. We omit the proof since it is quite similar to the case $k=2r$ (one has just to handle $\zeta_{r+1}$ instead of $\frac{\zeta_r+\zeta_{r+1}}2$).

\medskip

It remains to prove Lemma \ref{lemli} in order to finish the proof of Proposition \ref{propref}. 

\bigskip

{\it Proof of Lemma \ref{lemli}}: We first derive from Corollary \ref{corlemode} the following ODE system satisfied by the $L_i(s)$:
\begin{cl}\label{clode} Assuming that $k\ge 2$, it holds that
for all $i=1,...,k-1$ and  $s$ large enough, 
\[
\d\frac 1{c_1}L_i'(s) = - e^{-\b L_{i-1}(s)}+  2 e^{-\b L_i(s)}- e^{-\b L_{i+1}(s)}+\tilde R_i,
\]
where $|\tilde R_i| \le C J^{1+\delta_0}$ with the convention $L_0(s)\equiv +\infty$ and $L_k(s)\equiv +\infty$.
\end{cl}
%
The proof of (i) in Lemma \ref{lemli} is much longer than the proof of (ii). Therefore, we first derive (ii) assuming that (i) is true, then we prove (i).

\medskip

(ii) Using Claim \ref{clode}, the fact that $\tilde R_i=o(J)$ and \aref{defsigma}, we see that for $j=2,..,r$, we have 
\begin{eqnarray*}
\frac 1{c_1}S_j'(s)&=& e^{-\b L_j(s)}+e^{-\b L_{k-j}(s)}- e^{-\b L_{j-1}(s)}-e^{-\b L_{k-j+1}(s)}+o(J(s)),\\
\frac 1{c_1}S_1'(s)&=& e^{-\b L_1(s)}+e^{-\b L_{k-1}(s)}+o(J(s))
\end{eqnarray*}
as $s\to +\infty$. Therefore, using \aref{defsigma}, we write
\begin{eqnarray*}
\sigma'(s)= \sum_{j=1}^{r-1} \epsilon_0^{j-1}S_j'(s)
&=&  \sum_{j=1}^{r-1} \epsilon_0^{j-1}(1-\epsilon_0)(e^{-\b L_j(s)}+e^{-\b L_{k-j}(s)})\\
&+&\epsilon_0^{r-1}(e^{-\b L_r(s)}+e^{-\b L_{k-r}(s)})+o(J(s))
\end{eqnarray*}
as $s\to \infty$. Since for all $j=1,...,r-1$, $\frac{\epsilon_0^{r-1}}2 < \epsilon_0^{j-1}(1-\epsilon_0)<1$ and $\epsilon_0^{r-1}\le 1$, the conclusion of 
(ii) of Lemma \ref{lemli} follows by taking $s$ large enough.

\medskip

(i) Proceeding by contradiction, we consider a sequence $\varphi_n \to \infty$ and find two sequences $s_n <s_n'$ both going to $+\infty$ as $n\to +\infty$ such that
\begin{equation}\label{colere}
\bar L(s_n) = \varphi_n,\;\;\bar L(s_n') =1+ \varphi_n\mbox{ and }\forall s\in (s_n, s_n'],\;\;\bar L(s) >\varphi_n,
\end{equation}
where 
\begin{equation}\label{deflbar}
\bar L(s) = \max_{i=1,...,k-1}L_i(s) - \min_{i=1,...,k-1}L_i(s).
\end{equation}
Introducing
\begin{equation}\label{defmn}
\bar m_n = \min_{j=1,...,k-1}L_i(s_n)\mbox{ and }\bar M_n = \max_{j=1,...,k-1}L_i(s_n),
\end{equation}
we see from \aref{colere} that  
\begin{equation}\label{mm}
\bar M_n - \bar m_n = \bar L(s_n)=\varphi_n \to \infty\mbox{ as }n\to \infty.
\end{equation}
Up to extracting a subsequence, we assume that for all $i=1,...,k-1$, we have
\begin{equation}\label{convli}
L_i(s_n)-\bar m_n \to l_i\in [0, +\infty]\mbox{ as }n\to \infty.
\end{equation}
Introducing
\begin{eqnarray}
&&I_0=\{ i=1,...,k-1\;|\; l_i=0\},\;\;
I_+=\{ i=1,...,k-1\;|\; l_i\in (0, +\infty)\}\nonumber\\
&\mbox{and}&I_\infty=\{ i=1,...,k-1\;|\; l_i=+\infty\},\label{defi0}
\end{eqnarray}
we see that $\{1,...,k-1\}=I_0\cup I_+\cup I_\infty$ and that $I_0\neq \emptyset$.\\
Let us now introduce the following change of variables for all $i=1,...,k-1$ and $s\in [s_n, s_n']$,
\begin{equation}\label{defain}
a_{i,n}(\theta)=L_i(s)-\bar m_n\mbox{ where }s=s_n+\theta \frac{e^{\frac 2{p-1}\bar m_n}}{c_1}
\end{equation}
where $c_1$ appears in Claim \ref{clode}. Using Claim \ref{clode} and the fact that $\tilde R_i=o(J)$ as $s\to \infty$, we see that for all $i=1,...,k-1$ and $\theta \ge 0$,
\begin{equation}\label{eqain}
a_{i,n}'(\theta) = -e^{-\b a_{i-1,n}(\theta)}+2 e^{-\b a_{i,n}(\theta)}-e^{-\b a_{i+1,n}(\theta)}
+o\left(\sum_{j=1}^{k-1}e^{-\b a_{j,n}(\theta)}\right)
\end{equation}
as $n\to \infty$, where we take by convention $a_{0,n}(\theta) \equiv +\infty$ and $a_{k,n}(\theta) \equiv +\infty$. 
Now, fixing
\[
\theta_0=e^{-\frac 8{p-1}},
\]
we claim the following:
\begin{cl}\label{clest}
For all $\eta_0>0$, there exists $n_0\in \N$ such that for all $n\ge n_0$, $i=1,...,k-1$ and $\theta \in [0, \theta_0]$, we have
\begin{equation}\label{hadaf0}
-4\le a_{i,n}(\theta) \le \varphi_n +\eta_0.
\end{equation}
\end{cl}
{\it Proof}: See below.\Box\\
%
%
%
%
Using equation \aref{eqain} and Claim \ref{clest} (with $\eta_0=1$), we see that for some $C_0$, for all $n\ge n_0$, $i=1,...,k-1$ and $\theta \in [0, \theta_0]$, we have
\begin{equation}\label{dha}
|a_{i,n}'(\theta)|\le C_0.
\end{equation}
Using Ascoli's theorem, \aref{convli} and the definition \aref{defain} of $a_{i,n}$, we see that for all $i=1.,,,.k-1$, 
\begin{equation}\label{convain}
a_{i,n}(\theta)\to a_i(\theta)\mbox{ as }n\to \infty,\mbox{ uniformly for }\theta\in [0, \theta_0],
\end{equation}
where $a_i(\theta)\equiv +\infty$ if $i\in I_\infty$, and $a_i\in C([0, \theta_0])$ if not.
Moreover, getting to the limit in equation \aref{eqain} and introducing $a_0(\theta) \equiv +\infty$ and $a_{k-1}(\theta) \equiv +\infty$, we see that the vector $(a_i(\theta))_i$ satisfies the following system (first in an integral sense, then in a strong sense):
\begin{equation}\label{eqai}
a_i'(\theta) = -e^{-\b a_{i-1}(\theta)}+2 e^{-\b a_i(\theta)}-e^{-\b a_{i+1}(\theta)}
\end{equation}
for all $\theta\in [0, \theta_0]$, with initial data
\begin{equation}\label{defai0}
a_i(0)=l_i\in [0, +\infty].
\end{equation}
Using \aref{dha}, we see that
\begin{equation}\label{bound'}
\forall i\in I_0\cup I_+,\;\;\forall \theta \in [0, \theta_0],\;\;|a_i'(\theta)|\le C_0.
\end{equation}
Now, we claim that for all $i=1,...,k-1$ and for some $\theta_i\in (0, \theta_0]$, we have
\begin{equation}\label{infai}
a_i(\theta) >0\mbox{ for all }\theta\in (0, \theta_i].
\end{equation}
Indeed, if $i\in I_\infty$, this is clear from the convergence \aref{convain}.\\
If $i\in I_+$, then this is clear from \aref{bound'}.\\
Now, if $i\in I_0$, then we define $r_1$ and $r_2$ in $\{0,k\}\cup I_\infty\cup I_+$ such that $r_1<i<r_2$ and for all $j=r_1+1,...,r_2-1$, $j\in I_0$. Therefore, we have from \aref{defi0} and \aref{defai0}, 
\begin{equation}\label{hyppp}
a_{r_1}(0)>0,\;\;a_{r_2}(0)>0\mbox{ and if }r_1<j<r_2,\mbox{ then }a_j(0)=0.
\end{equation}
The following Claim allows us to conclude:
\begin{cl} \label{clo}
For some $\delta_0>0$,
we have 
\begin{equation}\label{con}
\forall j=r_1,...,r_2,\;\;\forall \theta\in (0, \delta_0],\;\;a_j(\theta)> 0.
\end{equation}
\end{cl}
{\it Proof}: See below.\Box\\
Applying Claim \ref{clo}, we see that \aref{infai} holds.

\medskip

Using \aref{infai} and introducing $\hat \theta= \inf_{i=1,...,k-1}\theta_i\in (0, \theta_0]$, we see that for some $\hat a>0$, we have
\[
\forall i=1,...,k-1,\;\;a_i(\hat \theta)\ge 4\hat a>0\mbox{ and }a_i(\theta) \ge 0\mbox{ for all }\theta\in [0, \hat \theta].
\]
From the convergence \aref{convain}, we see that for $n$ large enough, we have 
\[
\forall i=1,...,k-1,\;\;a_{i,n}(\hat \theta)\ge 3\hat a>0\mbox{ and }a_i(\theta) \ge -\frac 14\mbox{ for all }\theta\in [0, \hat \theta].
\]
Using Claim \ref{clest} (with $\eta_0=\min(\hat a, \frac 14)$), we see from the definitions \aref{defain} and \aref{deflbar} of $a_{i,n}$ and $\bar L(s)$ that for 
\[
\hat s_n =s_n+\hat \theta \frac{e^{\frac 2{p-1}\bar m_n}}{c_1},\;
\]
we have
\begin{eqnarray}
\bar L(\hat s_n)&\le& \bar m_n+\varphi_n +\eta_0-(\bar m_n+3\hat a)
= \varphi_n +\eta_0 - 3\hat a\le \varphi_n-2\hat a,\label{bem1}\\
\forall s\in [s_n, \hat s_n],\;\;\bar L(s)&\le& \bar m_n +\varphi_n+\eta_0-(\bar m_n-\frac 14)= \varphi_n +\eta_0+\frac 14\le \varphi_n +\frac 12.\label{bem2}
\end{eqnarray}
Since we have from \aref{colere} and \aref{bem2} that $\hat s_n \le s_n'$, we see that \aref{bem1} is in contradiction with \aref{colere}. It remains to prove Claims \ref{clest} and \ref{clo} in order to finish the proof of Lemma \ref{lemli}.

\bigskip

{\it Proof of Claim \ref{clest}}: Note first from the definitions \aref{defain} and \aref{defmn} of $a_{i,n}$, $\bar m_n$ and $\bar M_n$, and from \aref{mm}, that we have for all $n\in \N$ and $i=1,...,k$,
\begin{equation}\label{init}
0\le a_{i,n}(0) \le \varphi_n.
\end{equation}
Consider $\eta_0>0$. Proceeding by contradiction to prove Claim \ref{clest}, we define a subsequence (still denoted by $(a_{i,n})_n$) such that for all $n\in \N$, there is $\theta^*=\theta^*(n)\in (0, \theta_0)$ such that \aref{hadaf0} holds for all $i=1,...,k-1$ and $\theta \in [0, \theta^*]$, and for some $i^*=i^*(n)$, we have
\begin{equation}\label{b01}
\mbox{either }a_{i^*, n}(\theta^*)=-4\mbox{ or }a_{i^*, n}(\theta^*)=\varphi_n +\eta_0.
\end{equation}
From \aref{init}, we can define $\bar \theta = \bar \theta(n) \in [0,\theta^*]$ such that
\begin{equation}\label{inf}
a_{i^*,n}(\bar \theta)=\varphi_n\mbox{ and }\forall \theta \in [\bar \theta, \theta^*], a_{i^*, n}(\theta) \ge \varphi_n.
\end{equation} 
Since  \aref{hadaf0} holds for all $i=1,...,k-1$ and $\theta \in [0, \theta^*]$, we derive from equation \aref{eqain} that for $n$ large enough and $\theta \in [0, \theta^*]$,
\begin{eqnarray*}
a_{i^*,n}'(\theta) &\ge& -e^{-\b a_{i^*-1,n}(\theta)}+2 e^{-\b a_{i^*,n}(\theta)}-e^{-\b a_{i^*+1,n}(\theta)}
-\frac 1{k-1}\sum_{j=1}^{k-1}e^{-\b a_{j,n}(\theta)}\\
&\ge & -e^{\frac 8{p-1}}+0-e^{\frac 8{p-1}}-\frac 1{k-1}.(k-1)e^{\frac 8{p-1}}=-3e^{\frac 8{p-1}}.
\end{eqnarray*}
Integrating this for $\theta \in [0, \theta^*]$ and using \aref{init}, we get
\begin{equation}\label{b11}
a_{i^*,n}(\theta^*) \ge a_{i^*,n}(0)-3\theta^*e^{\frac 8{p-1}} \ge -3\theta_0e^{\frac 8{p-1}} = -3,
\end{equation}
on the one hand.\\
 On the other hand, from equation \aref{eqain}, \aref{hadaf0} (valid for $\theta \in [0, \theta^*]$) and \aref{inf}, we have for $n$ large enough and for all  $\theta \in [\bar \theta, \theta^*]$,
\begin{eqnarray*}
a_{i^*,n}'(\theta) &\le & -e^{-\b a_{i^*-1,n}(\theta)}+2 e^{-\b a_{i^*,n}(\theta)}-e^{-\b a_{i^*+1,n}(\theta)}
+\frac {\eta_0}{4(k-1)}\sum_{j=1}^{k-1}e^{-\b a_{j,n}(\theta)}\\
&\le & 0+2e^{-\frac {2\varphi_n}{p-1}}+0+\frac {\eta_0}{4(k-1)}(k-1)e^{\frac 8{p-1}}\le \frac{\eta_0}2e^{\frac 8{p-1}}.
\end{eqnarray*}
Integrating this for $\theta \in [\bar \theta, \theta^*]$, using \aref{inf} and the fact that $\theta^*- \bar \theta\le \theta_0= e^{-\frac 8{p-1}}$, we get
\begin{equation}\label{b21}
a_{i^*,n}(\theta^*) \le a_{i^*,n}(\bar \theta)+\frac{\eta_0}2e^{\frac 8{p-1}}(\theta^*- \bar \theta)
\le \varphi_n +\frac{\eta_0}2.
\end{equation}
Since \aref{b11} and \aref{b21} are in contradiction with \aref{b01}, this concludes the proof of Claim \ref{clest}.\Box

\bigskip

{\it Proof of Claim \ref{clo}}: 
%
%
%
%
%
%
%
%
In order to conclude, it is enough to prove that for all $l=0,..,l^*\equiv E(\frac{r_2-r_1}2)$, we have
\begin{equation}\label{obj}
a_{r_1+l}^{(l)}(0)>0,\;\;a_{r_2-l}^{(l)}(0)>0\mbox{ and if }r_1+l+1\le j\le r_2-l-1,\mbox{ then }a_j^{(l)}(0)=0,
\end{equation}
where we use the notation
\[
f^{(0)}=f\mbox{ and }f^{(l)}\mbox{ is the }l-\mbox{th derivative of }f.
\]
Indeed, if \aref{obj} holds, then it is easy to see that for any $j=r_1,...,r_2$, we have
\[
a_j^{(l)}(0)>0\mbox{ and if }0\le j'\le l-1,\mbox{ then }a_i^{(j')}(0)=0
\] 
where $l=j-r_1$ if $j\le l^*$ and $l=r_2-j$ if $j\ge l^*+1$, hence
\[
a_j(\theta)\sim \frac{a_j^{(l)}(0)}{l!}\theta^l\mbox{ as }\theta \to 0\mbox{ with }a_j^{(l)}(0)>0,
\]
and \aref{con} holds. It remains then to prove \aref{obj} in order to conclude the proof of Claim \ref{clo}.\\
If $l=0$, then \aref{obj} is true by \aref{hyppp}.\\
In order to prove \aref{obj} for all $l=1,..,l^*$, we proceed by induction.\\
- If $l=1$, then using \aref{eqai} and \aref{obj} (with $l=0$), we write
\begin{eqnarray*}
a_{r_1+1}'(0)&=&  -e^{-\b a_{r_1}(0)}+2-e^{-\b a_{r_1+2}(0)}
\ge  1-e^{-\b a_{r_1}(0)}>0,\\
a_{r_2-1}'(0)&=&  -e^{-\b a_{r_2-2}(0)}+2-e^{-\b a_{r_2}(0)}
\ge  1-e^{-\b a_{r_2}(0)}>0,\\
\mbox{if }r_1+2\le j \le r_2-2,\;
a_j'(0)&=&-1+2-1=0.
\end{eqnarray*}
- Take $l=2,...,l^*$ and assume that \aref{obj} holds for any $l'=0,...,l-1$. 
Therefore, it is easy to see that
\begin{equation}\label{anne}
\begin{array}{rcl}
&&a_{r_1+l-1}^{(l-1)}(0)>0\mbox{ and if }0\le j'\le l-2,\;\;a_{r_1+l-1}^{(j')}(0)=0,\\
\\
\mbox{if }r_1+l\le j\le r_2-l,&& a_j^{(j')}(0)=0\mbox{ for all }j'=0,...,l-1,\\
\\
&&a_{r_2-l+1}^{(l-1)}(0)>0\mbox{ and if }0\le j'\le l-2,\;\;a_{r_2-l+1}^{(j')}(0)=0.
\end{array}
\end{equation}
In the following, we prove that \aref{obj} holds for $l$. Starting from equation \aref{eqai}, we prove thanks to a straightforward induction that 
\begin{equation}\label{ail}
\begin{array}{rcl}
&&\frac{(p-1)}2a_j^{(l)}=[a_{j-1}^{(l-1)}+P_l(a_{j-1}',...,a_{j-1}^{(l-2)})]e^{-\b a_{j-1}}\\
\\
&-&2 [a_{j}^{(l-1)}+P_l(a_{j}',...,a_{j}^{(l-2)})]e^{-\b a_{j}}
+[a_{j+1}^{(l-1)}+P_l(a_{j+1}',...,a_{j+1}^{(l-2)})]e^{-\b a_{j+1}}
\end{array}
\end{equation}
where $P_2\equiv 0$ and for $l\ge 3$, $P_l$ is a polynomial of $l-2$ variables satisfying
\begin{equation}\label{condP}
P_l(0,...,0)=0.
\end{equation}
Using \aref{ail}, \aref{anne} and \aref{condP}, we write
\begin{eqnarray*}
a_{r_1+l}^{(l)}(0)&=&a_{r_1+l-1}^{(l-1)}(0)+a_{r_1+l+1}^{(l-1)}(0)
\ge  a_{r_1+l-1}^{(l-1)}(0) >0,\\
a_{r_2-l}^{(l)}(0)&=&a_{r_2-l-1}^{(l-1)}(0)+a_{r_2-l+1}^{(l-1)}(0)
\ge  a_{r_2-l+1}^{(l-1)}(0) >0,\\
\mbox{if }r_1+l+1\le j \le r_2-l-1,\;\;a_j^{(l-1)}(0)&=&0.
\end{eqnarray*}
which is the conclusion of \aref{obj} with the index $l$. Thus, \aref{obj} holds. This concludes the proof of Claim \ref{clo} and Lemma \ref{lemli} too. Thus, Proposition \ref{propref} holds. \Box






\subsection{Refinement of \aref{cprofile0} for $k\ge 2$}\label{subseclemode}
Note that the case $k=1$ has been already treated in \cite{MZjfa07} giving rise to estimate \aref{profile}. As announced in the beginning of the section, we assume that $k \ge 2$ and claim the following: 
\begin{prop}\label{lemq}{\bf (Size of $q$ in terms of the distance between solitons)}
There exists another set of parameters (still denoted by $\zeta_1(s)$, ... $\zeta_k(s)$) such that \aref{cprofile0} and \aref{div} hold and for some  $s^*\in \R$ and for all $s\ge s^*$, 
\begin{equation}\label{defq}
\|q(s)\|_{\H}\le C \sum_{i=1}^{k-1}h(\zeta_{i+1}(s)-\zeta_i(s)),
\end{equation}
 where
\[
q=\vc{q_1}{q_2}= \vc{w}{\ps w} - \d\sum_{i=1}^{k} e_i\vc{\kappa(d_i)}0,
\]
and 
\begin{equation}\label{defh}
h(\zeta)= e^{-\frac p{p-1} \zeta} \mbox{ if }p<2,\;\;h(\zeta)= e^{-2 \zeta} \sqrt \zeta\mbox{ if }p=2 \mbox{ and }h(\zeta)= e^{-\frac 2{p-1}\zeta}\mbox{ if }p>2.
\end{equation}
\end{prop}
Before proving the estimate, we need to use a modulation technique to slightly change the $\zeta_i(s)$ in order to guarantee some orthogonality conditions. In order to do so, we need to introduce for $\lambda=0$ or $1$, for any $d\in (-1,1)$ and $r\in \H$, 
\begin{equation}\label{defpdi}
\pi_\lambda^d(r) =\phi\left(W_\lambda(d), r\right)
\end{equation}
where: 
\begin{equation}\label{defPhi}
\phi(q,r)
= \int_{-1}^1 \left(q_1r_1+q_1' r_1' (1-y^2)+q_2r_2\right)\rho dy
=\int_{-1}^1 \left(q_1\left(-\L r_1+r_1\right) +q_2 r_2\right)\rho dy,
\end{equation}
$W_\lambda(d,y)= (W_{\lambda,1}(d,y), W_{\lambda,2}(d,y))$,
\begin{equation}\label{defWl2-0}
W_{1,2}(d,y)(y)= c_1(d) \frac{1-y^2}{(1+dy)^{\frac 2{p-1}+1}},\;\;
W_{0,2}(d,y) = c_0\frac {y+d}{1+dy}\kappa(d,y),
\end{equation}
with $0<c_1(d) \le C(1-d^2)^{\frac 1{p-1}}$, $c_0>0$,\\
 and $W_{\lambda,1}(d,y)\in \H_0$ is uniquely determined as the solution of 
\begin{equation}\label{eqWl1-0}
-\L r + r =\left(\lambda - \frac{p+3}{p-1}\right)W_{\lambda,2}(d) - 2 y\py W_{\lambda,2}(d)+ \frac 8{p-1} \frac{W_{\lambda,2}(d)}{1-y^2}
\end{equation}
(in \cite{MZjfa07}, we defined $W_{0,2}(d,y)$ by $\frac {c_0(d)(y+d)}{(1+dy)^{\frac 2{p-1}+1}}$ with  
\[
1= c_0(d)(1-d^2)^{\frac 1{p-1}}\frac 4{p-1} \int_{-1}^1 \frac{(y+d)^2}{(1+dy)^{\frac 4{p-1}+2}} \frac \rho{1-y^2} dy.
\]
 Setting $y=\tanh \xi$, we compute the integral and get $c_0(d)= c_0'(1-d^2)^{\frac 1{p-1}}$. Using \aref{defkd}, we get \aref{defWl2-0}).
Recall from Lemma 4.4 page 85 in \cite{MZjfa07} that
\begin{equation}\label{normw}
\forall d\in(-1,1),\;\;\|W_\lambda(d)\|_{\H} \le C.
\end{equation}
We now have the following:
\begin{lem}[Modulation technique]\label{clmod}Assume that $k\ge 2$.\\
(i) {\bf (Choice of the modulation parameters)} There exist other values of the parameters (still denoted by $d_i(s)$) of class $C^1$, such that $\zeta_{i+1}(s) - \zeta_i(s)\to \infty$ as $s\to \infty$ where $d_i(s) = -\tanh \zeta_i(s)$,
\begin{equation}\label{kill}
\|q(s)\|_{\H} \to 0\mbox{ and }\pi_0^{d_i(s)}(q)=0 \mbox{ for all }i=1,..,k,
\end{equation}
where $\pi_0^d$ and $q$ are defined in \aref{defpdi} and \aref{defq} respectively.\\
(ii) {\bf (Equation on $q$)} For $s$ large, we have
\begin{equation}\label{eqq}
\d\frac \partial {\partial s}
\left(
\begin{array}{l}
q_1\\
q_2
\end{array}
\right)
=L
\left(
\begin{array}{l}
q_1\\
q_2
\end{array}
\right)
-\sum_{j=1}^k e_jd_j'(s)\vc{\partial_d \kappa(d_j(s),y)}{0}
+\vc{0}{R}
+\left(
\begin{array}{l}
0\\
f(q_1)
\end{array}
\right)
\end{equation}
\begin{equation*}
\begin{array}{rcl}
\mbox{where }
L\vc{q_1}{q_2}&=&\vc{q_2}{\L q_1+\psi q_1-\frac{p+3}{p-1}q_2-2y\py q_2},\\\\
\psi(y,s)&=&p|K(y,s)|^{p-1} -\frac{2(p+1)}{(p-1)^2},\;\;
K(y,s) = \sum_{j=1}^k e_j \kappa(d_j(s),y),\\\\
f(q_1)&=&|K+q_1|^{p-1}(K+q_1)- |K|^{p-1}K- p|K|^{p-1} q_1,\\\\
R&=& |K|^{p-1}K- \sum_{j=1}^k e_j \kappa(d_j)^p.
\end{array}
\end{equation*}
\end{lem}
{\bf Remark}: From the modulation technique, it is clear that the distance between old and new parameter $\zeta_i(s)$ goes to zero as $s\to \infty$.\\
{\it Proof}: See the proof of Proposition 5.1 in \cite{MZjfa07} where the case $k=1$ is treated. There is no difficulty in adapting the proof to $k\ge 2$.\Box

\bigskip

In the following, we will show that Proposition \ref{lemq} holds with the set of parameters $\zeta_1(s)$,..,$\zeta_k(s)$ given by the modulation technique of Lemma \ref{clmod}. Before giving the proof, we start by reformulating the problem.

Let us first remark that equation \aref{eqq} can be localized near each soliton's center which allows us to view it locally as a perturbation of the case of one soliton already treated in \cite{MZjfa07}. For this, given $i=1,...,k$, we need to expand the linear operator of equation \aref{eqq} as 
\begin{eqnarray}
&&L(q) = L_{d_i(s)}(q) + (0, V_i(y,s) q_1)\mbox{ with }\nonumber\\
L_d\vc{q_1}{q_2}&=&\vc{q_2}{\L q_1+(p\kappa(d_i(s),y)^{p-1}-\frac{2(p+1)}{(p-1)^2}) q_1-\frac{p+3}{p-1}q_2-2yq_2'},\nonumber\\
V_i(y,s)&=& p|K(y,s)|^{p-1} - p \kappa(d_i(s),y)^{p-1}.\label{defvi-0}
\end{eqnarray}
Since the solitons' sum is decoupled (remember from (i) of Lemma \ref{clmod} that
\begin{equation}\label{decouple}
\xi_{i+1}-\xi_i \to \infty\mbox{ as }s\to \infty),
\end{equation}
the properties of $L_{d_i(s)}$ will be essential in our analysis.\\
From section 4 in \cite{MZjfa07}, we know that for any $d\in (-1,1)$, the operator $L_d$ has $1$ and $0$ as eigenvalues, the rest of the eigenvalues are negative. More precisely, introducing 
\begin{equation}\label{deffld}
F_1(d,y)=(1-d^2)^{\frac p{p-1}}\vc{(1+dy)^{-\frac 2{p-1}-1}}{(1+dy)^{-\frac 2{p-1}-1}},\;\; F_0(d,y)=(1-d^2)^{\frac 1{p-1}}\vc{\d\frac{y+d}{(1+dy)^{\frac 2{p-1}+1}}}{0},
\end{equation}
we have
\begin{equation}\label{79bis}
L_d (F_\lambda(d)) = \lambda F_\lambda(d)\mbox{ and }
\|F_1(d)\|_{\H}+\|F_0(d)\|_{\H}\le C.
\end{equation}
The projection on $F_\lambda(d)$ is defined in \aref{defpdi} by $\pi_\lambda^d(r) =\phi\left(W_\lambda(d), r\right)$. 
%
%
Of course,
\begin{equation}\label{79bis*}
L_d^* W_\lambda(d) = \lambda W_\lambda(d)
\end{equation}
where $L_d^*$ is the conjugate of $L_d$ with respect to the inner product $\phi$, and the choice of the constants $c_1(d)$ and $c_0$ guarantees the orthogonality condition
\begin{equation}\label{orth}
\pi_\lambda^d(F_\mu(d)) = \phi(W_\lambda(d), F_\mu(d)) = \delta_{\lambda,\mu}.
\end{equation}
In the following, we give a decomposition of the solution which is well adapted to the proof:
\begin{lem}[Decomposition of $q$]\label{lemdecomp}
If we introduce for all $r$ and $\r$ in $\H$ the operator
 $\pi_-(r)\equiv r_-(y,s)$ defined by
\begin{equation}\label{decomp}
r(y,s)= \sum_{i=1}^{k} \left(\pi^{d_i(s)}_1(r) F_1(d_i(s),y) + \pi^{d_i(s)}_0(r) F_0(d_i(s),y)\right) + \pi_-(r)
\end{equation}
and the bilinear form 
\begin{equation}\label{defphi}
\varphi(r,\r) = \iint \left(r_1'\r_1' (1-y^2)-\psi r_1\r_1+r_2\r_2\right)\rho dy
\end{equation}
 where $\psi(y,s)$ is defined in \aref{eqq}, then:\\
(i) for $s$ large enough and for all $r$ and $\r$ in $\H$, we have
\begin{equation}\label{contphi}
\left| \varphi(r,\r) \right| \le C\|r\|_{\H} \|\r\|_{\H},
\end{equation}
(ii) for some $C_0>0$ and for all $s$ large enough, we have:
\begin{eqnarray}
q(y,s)&=& \sum_{i=1}^{k} \alpha^i_1(s)F_1(d_i(s),y) + q_-(y,s),\label{qdecomp}\\
\frac 1{C_0}\|q_-(s)\|_{\H}^2 - C_0\bar J(s)^2\|q(s)\|_{\H}^2&\le& A_-(s)\le C_0\|q_-(s)\|_{\H}^2,\label{48bis}\\
\frac 1{C_0}\|q(s)\|_{\H}^2&\le& \sum_{i=1}^k \left(\alpha^i_1(s)\right)^2+A_-(s)\le C_0 \|q(s)\|_{\H}^2\label{equivq}
\end{eqnarray}
where
\begin{equation}\label{defa-}
\bar J(s)= \sum_{j=1}^{k-1} 
(\zeta_{j+1}-\zeta_j)e^{-\frac 2{p-1}(\zeta_{j+1}-\zeta_j)},\;\;
\alpha^i_\lambda(s)= \pi_\lambda^{d_i(s)}(q(s))\mbox{ and }A_-(s)=\varphi(q_-(s),q_-(s)).
\end{equation} 
\end{lem}
{\bf Remark}: Note that the choice of $d_i(s)$ made in \aref{kill} guarantees that for $s$ large enough,
\begin{equation}\label{peur}
\forall i =1,...,k,\;\;\alpha_0^i(s)\equiv {\alpha_0^i}'(s)\equiv 0.
\end{equation} 
Moreover, we see from \aref{48bis} that $A_-(s)$ is nearly positive and nearly equivalent to $\|q_-\|_{\H}^2$. Note that from \aref{imen2} (proved in the proof of Claim \ref{lemdecomp}), \aref{qdecomp}, \aref{79bis} and \aref{48bis} we have for $s$ large,
\begin{equation}\label{91,5}
|\alpha_1^i(s)|\le C\|q(s)\|_{\H}\mbox{ and }\|q_-(s)\|_{\H}\le C\min\left(\|q(s)\|_{\H}, \sqrt{|A_-(s)|}+\bar J(s)\|q(s)\|_{\H}\right).
\end{equation}
{\bf Remark}: The operator $\pi_-$ depends on the time variable $s$. In \cite{MZjfa07}, we had only one soliton, and we decomposed $q$ as follows:
\begin{equation}\label{defpi-d}
q(y,s)= \pi^{d}_1(q) F_1(d,y) + \pi^{d}_0(q) F_0(d,y) + \pi^d_-(q),
\end{equation}
where we had only one $d(s)$ (note that this decomposition is in fact a definition of the operator $\pi^d_-$). Here, due to \aref{decouple}, we have a decoupling effect, in the sense that $\pi_\lambda^{d_j(s)}(q)$ for $j\neq i$ cannot be ``seen'' when $y$ is close to $-d_i(s)$, the ``center'' of the soliton $\kappa(d_i(s), y)$. Hence, $\pi_-(q)$ is more or less $\pi^{d_i(s)}_-(q)$ and we are reduced to the situation of one soliton already treated in \cite{MZjfa07}. This idea will be essential in our proof since given some $i=1,...,k$, we have two types of terms in equation \aref{eqq}:\\
- terms involving the soliton $\kappa(d_i(s),y)$ for which we refer the reader to \cite{MZjfa07},\\
- interaction terms involving a different soliton $\kappa(d_j(s),y)$ which we treat in details.\\ 
{\it Proof of Lemma \ref{lemdecomp}}: See Appendix \ref{appdecomp}.\Box

\bigskip

In order to prove Proposition \ref{lemq}, we project equation \aref{eqq} according to the decomposition \aref{decomp}. More precisely, we have the following: 
\begin{lem}\label{lemproj}
For $s$ large enough, the following holds:\\ 
(i) {\bf (Control of the positive modes and the modulation parameters)} 
\begin{equation}\label{eqa-star}
\forall i=1,...,k,\;\;\left|\ap'(s)- \ap(s)\right|+|\zeta_i'(s)|\le C\|q(s)\|_{\H}^2 + C J(s)
\end{equation}
where $J(s)$ is defined in \aref{defJ}.\\
(ii) {\bf (Control of the negative part)}
\begin{eqnarray}
\left(R_-+\frac 12 A_-\right)'&\le& - \frac 3{p-1}\iint q_{-,2}^2 \frac \rho{1-y^2} dy
+o\left(\|q(s)\|_{\H}^2\right) + C \sum_{m=1}^{k-1} \left(h(\zeta_{m+1}-\zeta_m)\right)^2\nonumber\\
& +& C J(s)\sqrt{|A_-(s)|}\label{eqam}
\end{eqnarray}
for some $R_-(s)$ satisfying 
\begin{equation}\label{boundr-}
|R_-(s)|\le C \|q(s)\|_{\H}^{\bar p+1}
\end{equation}
where $\bar p = \min(p,2)$ and $h(s)$ is defined in \aref{defh}.\\ 
(iii) {\bf (An additional relation)}
\begin{equation} \label{153}
\frac d{ds}\iint q_1q_2 \rho \le 
-\frac 45 A_-+C \sum_{m=1}^{k-1} h(\zeta_{m+1}-\zeta_m)^2+ C \iint q_{-,2}^2 \frac \rho{1-y^2}+C\sum_{i=1}^k\left(\ap\right)^2.
\end{equation}
\end{lem}
{\it Proof}: See Appendix \ref{appproj}.\Box

\bigskip

With Lemma \ref{lemproj}, we are ready to prove Proposition \ref{lemq}.

\medskip

{\it Proof of Proposition \ref{lemq}}: 
 We proceed as in section 5.3 page 113 in \cite{MZjfa07}, though the situation is a bit different because of the presence of the forcing terms $J(s)$ and $\sum_{i=1}^k h(\zeta_{i+1} - \zeta_i)^2$ in the differential inequalities in Lemma \ref{lemproj}. 

If we introduce
\begin{equation}\label{defab}
a(s) = \sum_{i=1}^k \alpha_1^i(s)^2,\;\; b(s) = A_-(s) + 2 R_-(s)\mbox{ and } H(s) = \sum_{m=1}^{k-1} h(\zeta_{m+1}-\zeta_m)^2
\end{equation}
where $h$ is defined in \aref{defh}, then we see from (i) of Lemma \ref{clmod} and \aref{boundr-} that
\begin{equation}\label{limitab}
a(s)+b(s)+H(s) \to 0\mbox{ as }s\to \infty.
\end{equation}
Moreover, we see from \aref{boundr-} and \aref{equivq} that  $|b- A_-|\le \frac1{1000}\left(A_- + \sum_{i=1}^k (\alpha_1^i)^2\right)$ for $s$ large enough, hence
\begin{equation}\label{equivab}
\frac{99}{100} A_- - \frac 1{100} a \le b \le \frac{101}{100} A_- + \frac 1{100} a.
\end{equation}
Therefore, since $J(s) \le H(s)$ by \aref{defJ} and \aref{defh}, we have for $s$ large, 
\[
\forall \epsilon>0,\;\;CJ\sqrt{|A_-|} \le \epsilon (a+b) + \frac C \epsilon H(s).
\]
Using \aref{defab}, \aref{equivab} and \aref{limitab}, we rewrite estimates \aref{equivq} and Lemma \ref{lemproj} with the new variables, in the following:
\begin{cor}\label{corproj}{\bf (Equations in the new framework)}
There exists $K_0\ge 1$ such that for all $\epsilon>0$, there exists $s_0(\epsilon)\in \R$ such that for all $s\ge s_0(\epsilon)$, the following holds:\\
(i) {\bf (Size of the solution)}
\begin{eqnarray}
\frac 1{K_0}(a+b) \le \|q\|_{\H}^2 &\le & K_0 (a+b),\label{104bis}\\
 \left|\iint q_1 q_2 \rho dy \right| &\le & K_0 (a+b).\label{qq}
\end{eqnarray}
(ii) {\bf (Equations)} 
\begin{eqnarray}
\frac 32 a - \epsilon b-K_0 H \le a' &\le & \frac 52 a + \epsilon b+K_0 H,\nonumber\\
b' &\le& -\frac 6{p-1} \iint q_{-,2}^2 \frac \rho{1-y^2} dy + \epsilon(a+b) + \frac{K_0}\epsilon H,\nonumber\\
\frac d{ds} \iint q_1 q_2 \rho dy &\le & - \frac 3 5 b +K_0 \iint q_{-,2}^2 \frac \rho{1-y^2} dy +K_0 a + K_0H,\nonumber\\
|H'|&\le & \epsilon H.\label{eqg}
\end{eqnarray}
\end{cor}
We proceed in 2 steps:\\
- In Step 1, we show that $a$ is controlled by $b+H$.\\
- In Step 2, we show that $b$ is controlled by $H$ and conclude the proof using \aref{104bis}.

\bigskip

{\bf Step 1: $a$ is controlled by $b+H$}

We claim that for $\epsilon$ small enough, we have:
\begin{equation}\label{cl1}
\forall s\ge s_0(\epsilon),\;\;a(s)\le \epsilon b(s) + \frac{K_0}\epsilon H(s).
\end{equation}
Indeed, from Corollary \ref{corproj}, we see that for all $s\ge s_0(\epsilon)$, we have
\begin{eqnarray*}
a' &\ge & \frac 32 a - (\epsilon b+ \frac{K_0}\epsilon H),\\
(\epsilon b+ \frac{K_0}\epsilon H)' &\le & 2\epsilon(\epsilon b+ \frac{K_0}\epsilon H)+\epsilon^2 a.
\end{eqnarray*}
Introducing $\gamma_1(s) = a(s) - (\epsilon b(s)+ \frac{K_0}\epsilon H(s))$, we see that for all $s\ge s_0(\epsilon)$,
\begin{eqnarray*}
\gamma_1' = a' - (\epsilon b'+ \frac{K_0}\epsilon H')&\ge& \frac 32 a - (\epsilon b+ \frac{K_0}\epsilon H) 
-2\epsilon(\epsilon b+ \frac{K_0}\epsilon H)-\epsilon^2 a\\
& =& (\frac 32  -\epsilon^2 -1-2\epsilon)a +(1+2\epsilon)\gamma_1 \ge \gamma_1
\end{eqnarray*}
if $\epsilon$ is small enough. 
Since $\gamma_1(s) \to 0$ as $s\to \infty$ (see \aref{limitab}), this implies $\gamma_1(s) \le 0$, hence \aref{cl1} follows.

\bigskip

{\bf Step 2: $b$ is controlled by $H$}

We claim that in order to conclude, it is enough to prove for some $K_1>0$ that
\begin{equation}\label{cl2}
\forall s\ge s_0(\epsilon),\;\;
f(s) \le K_1 H(s)
\end{equation}
where 
\begin{equation}\label{deff}
f= b +\eta \iint q_1 q_2 \rho dy\mbox{ and }\eta = \frac 12 \min\left( \frac 1{2K_0}, \frac 6{(p-1) K_0}\right).
\end{equation}
Indeed, using \aref{qq} and \aref{cl1}, and taking $\epsilon$ small enough, we get for all $s\ge s_0(\epsilon)$,\\
 $\left|\iint q_1 q_2 \rho dy\right| \le  2K_0 b+\frac{K_0^2}\epsilon H$ and $|f-b|\le 2K_0\eta b+\eta \frac{K_0^2}\epsilon H\le \frac b2 +\eta \frac{K_0^2}\epsilon H$, hence
\begin{equation}\label{boundf}
\frac b2 - \eta \frac{K_0^2}\epsilon  H \le f \le 2b + \eta \frac{K_0^2}\epsilon H.
\end{equation}
Therefore, if \aref{cl2} holds, then using \aref{104bis}, 
\aref{boundf} and \aref{cl1}, we see that for some $K_2>0$ and for all $s\ge s_0(\epsilon)$, $\|q(s)\|_{\H}^2 \le K_0(a(s)+b(s)) \le K_2 H(s)$ which is the desired conclusion of Proposition \ref{lemq}. It remains to prove \aref{cl2}.\\
Using Corollary \ref{corproj}, \aref{cl1}, \aref{deff} and the fact that $K_0\ge 1$, and taking $\epsilon$ small enough, we get for all $s\ge s_0(\epsilon)$:
\begin{eqnarray}
b' &\le& -\frac 6{p-1} \iint q_{-,2}^2 \frac \rho{1-y^2} dy + 2\epsilon b + 2\frac{K_0}\epsilon H,\label{eqb}\\
\frac d{ds} \iint q_1 q_2 \rho dy &\le & - \frac 2 5 b +K_0 \iint q_{-,2}^2 \frac \rho{1-y^2} dy + 2\frac{K_0^2}\epsilon H,\label{eqqq}\\
f' &\le& -(\frac 25 \eta- 2\epsilon) b
- \left(\frac 6{p-1}- K_0\eta\right)\iint q_{-,2}^2 \frac \rho{1-y^2} dy
+(2\frac{K_0}\epsilon +2\frac{K_0^2}\epsilon \eta)H\nonumber\\ 
&\le& -\frac \eta 4 b +3\frac{K_0}\epsilon H
 \le -\frac \eta 8 f+ 4 \frac{K_0}\epsilon H.\label{eqf}
\end{eqnarray}
If $\gamma_2(s) = f(s) - \frac{64 K_0}{\eta\epsilon} H(s)$, then we write from \aref{eqg} and \aref{eqf}, for all $s\ge s_0(\epsilon)$,
\[
\gamma_2' = f' - \frac{64 K_0}{\eta\epsilon}  H' \le -\frac \eta 8 f +4\frac{K_0}\epsilon H+\frac{64 K_0}{\eta\epsilon}  \epsilon H
= -\frac \eta 8 \gamma_2 + \frac{K_3}\epsilon H \le -\frac \eta 8 \gamma_2 
\]
because $K_3 = -\frac{64 K_0}\eta \frac \eta 8 +\frac{64 K_0}\eta \epsilon +4K_0
=-4K_0+\frac{64 K_0}\eta \epsilon \le 0$ if $\epsilon$ is small enough. Therefore, for all $s\ge s_0(\epsilon)$, $\gamma_2(s) \le e^{-\frac \eta 8 (s-s_0)}\gamma_2(s_0)$, hence
\begin{equation}\label{bach1}
f(s) \le  \frac{64 K_0}{\eta\epsilon} H(s) + e^{-\frac \eta 8 (s-s_0)}|\gamma_2(s_0)|.
\end{equation}
Since we have from \aref{eqg} and \aref{defab}, 
\begin{equation}\label{bach2}
H(s) \ge e^{-\epsilon(s-s_0)}H(s_0)\mbox{ and }H(s_0)>0,
\end{equation}
taking $\epsilon \le \frac \eta 8$, we see that
\aref{cl2} follows from \aref{bach1} and \aref{bach2}. This concludes the proof of Proposition \ref{lemq}.\Box

\subsection{An ODE system satisfied by the solitons' centers}\label{applemode0}
With Proposition \ref{lemq}, we are ready to prove Proposition \ref{lemode} now. The proof consists in refining the projection of equation \aref{eqq} with the projector $\pi^d_0$ \aref{defpdi}, already performed in the proof of (i) of Lemma \ref{lemproj} (see Part 1 page \pageref{adel1-9}). 

\bigskip

{\it Proof of Proposition \ref{lemode}}: 
Using \aref{adel1-9}, \aref{nonl-9}, \aref{69bis-9}, \aref{70-9}, \aref{74-9}, the differential inequality \aref{eqa-star} on $\zeta_i$ and the fact that $\alpha_0^i(s)\equiv {\alpha_0^i}'(s)\equiv 0$ (see \aref{peur}), we write for some $\delta_2>0$ and for $s$ large enough, 
\begin{equation}\label{zaid}
\left|-e_i \frac{2\kappa_0}{p-1} \zeta_i'(s)+\pi_0^{d_i(s)}\vc{0}{R}\right|\le C\|q(s)\|_{\H}^2 +CJ(s)^{1+\delta_2}.
\end{equation}
Since we have from Proposition \ref{lemq},
\begin{equation}\label{amr}
\|q(s)\|_{\H}^2 \le C\sum_{i=1}^{k-1}\left(h(\zeta_{i+1}-\zeta_i)\right)^2\le C J(s)^{1+\delta_3},
\end{equation}
for some $\delta_3>0$ and for $s$ large enough, 
where $h(\zeta)$ is defined in \aref{defh}, it is clear that if one proves that for some $c_1'>0$, $\delta_4>0$ and for $s$ large enough,
\begin{equation}\label{coc}
\frac 1{c_1'}\pi_0^{d_i(s)}\vc{0}{R} =-e_{i-1} e^{-\frac 2{p-1}|\zeta_i(s)-\zeta_{i-1}(s)|}
+ e_{i+1} e^{-\frac 2{p-1}|\zeta_{i+1}(s)-\zeta_i(s)|}+CJ(s)^{1+\delta_4}
\end{equation}
(with the convention $\zeta_0(s)\equiv-\infty$ and $\zeta_{k+1}(s)\equiv+\infty$), 
then, Proposition \ref{lemode} immediately follows from \aref{zaid} and \aref{amr} (with $\delta_0 =\min(\delta_2, \delta_3, \delta_4)$). It remains to prove \aref{coc} in order to conclude the proof of Proposition \ref{lemode}.

\bigskip
{\it Proof of \aref{coc}}: 
We claim first that 
\begin{eqnarray}
&&\left|R- \sum_{j=1}^kp\kappa(d_j(s))^{p-1}\un{y_{j-1}(s)}{y_{j}(s)} \sum_{l\neq j}e_l \kappa(d_l(s)) \right|\nonumber\\
&\le& C\sum_{j=1}^k \kappa(d_j(s))^{p-\p}\un{y_{j-1}(s)}{y_{j}(s)} \sum_{l\neq j} \kappa(d_l(s))^\p \label{krim}
\end{eqnarray}
where $y_i$ are the solitons' separators defined in \aref{0defy}.\\
 Indeed, let us take $y\in (y_{j-1}(s), y_j(s))$ and set $X = (\sum_{l\neq j}e_l \kappa(d_l(s)))/e_j \kappa(d_j(s))$. From the fact that $\zeta_{j+1}(s)-\zeta_j(s)\to \infty$ and \aref{0defy}, we have $|X|\le 3$ hence
\[
||1+X|^{p-1}(1+X) -1 -pX|\le CX^2
\]
and for $y\in (y_{j-1}(s), y_j(s))$ and $s$ large,
\begin{eqnarray*}
&&||K|^{p-1}K - e_j \kappa(d_j(s))^p -p\kappa(d_j(s))^{p-1} \sum_{l\neq j}e_l \kappa(d_l(s)) |\\
&\le& C \kappa(d_j(s))^{p-2} \sum_{l\neq j} \kappa(d_l(s))^2 
\end{eqnarray*}
Since for $y\in (y_{j-1}(s), y_j(s))$,  $|\sum_{l\neq j} e_l \kappa(d_l(s))^p|\le \sum_{l\neq j}\kappa(d_l(s))^p$ and $\kappa(d_j(s))\ge \kappa(d_l(s))$ if $l\neq j$, this concludes the proof of \aref{krim}.\\
 Now, we prove \aref{coc}. Using \aref{krim}, \aref{defpdi}, \aref{defWl2-0} and the notations of Lemma \ref{lemtech}, we write
\[
\left|\pi_0^{d_i}\vc{0}{R}-pc_0\sum_{j=1}^k\sum_{l\neq j}A_{i,j,l}\right|\le C \sum_{j=1}^k\sum_{l\neq j} B_{i,j,l}.
\]
Since we have from (iii) and (iv) of Lemma \ref{lemtech}, $|A_{i,j,l}|+|B_{i,j,l}|\le CJ^{1+\delta_7}$, except for $A_{i,i,l}$ with $l=i\pm 1$ where we have $|A_{i,i,l}- c_1''' \sgn(l-i)e^{-\frac 2{p-1}|\zeta_i-\zeta_{i\pm 1}|}|\le CJ^{1+\delta_4}$ where $\delta_4=\min(\delta_5, \delta_6)>0$, we get \aref{coc}.
Since Proposition \ref{lemode} follows form \aref{coc} and \aref{zaid}, this concludes the proof of Proposition \ref{lemode} too. \Box

\subsection{The blow-up set has the corner property near $x_0\in \SS$ when $k(x_0)\ge 2$}\label{sub4}
We derive here the following consequence of Proposition \ref{propref}:
\begin{prop}\label{phat}{ \bf (Existence of signed lines and the corner property near $x_0\in \SS$)}  If $x_0\in \SS$ with $k(x_0)\ge 2$, then:\\
(i) For all $j=1,..,k$, 
\[
u(z_j(t),t)\sim e_j^* \kappa_0 \cosh^{\frac 2{p-1}}\zeta_j(s)(T(x_0)-t)^{-\frac 2{p-1}}\mbox{ as }t\to T(x_0),
\]
 where $t\mapsto z_j(t)$ is continuous and defined by
\begin{equation}\label{defzj}
z_j(t)= x_0+(T(x_0)-t)\tanh \zeta_j(s)\mbox{ with }s=-\log(T(x_0)-t).
\end{equation}
(ii) We have for some $\delta_0>0$ and $C_0>0$, 
\begin{equation}\label{chapeau}
\mbox{if }|x-x_0|\le \delta_0,\mbox{ then }T(x_0)-|x-x_0| \le T(x)\le T(x_0)-|x-x_0| + \frac{C_0|x-x_0|}{|\log(x-x_0)|^{\frac{(k-1)(p-1)}2}}.
\end{equation}
\end{prop}
{\bf Remark}: The point $z_j(t)$ corresponds in the original variables to the center of the $j$-th soliton in the description \aref{cprofile0}.\\
{\bf Remark}: In the next section, we prove that for all $x_0\in \SS$, we have $k(x_0)\ge 2$. Thus, the result will hold for all $x_0\in \SS$.\\
{\it Proof}: It follows from Proposition \ref{propref} and the following:
\begin{lem}\label{pupper}{\bf (Upper blow-up bound for equation \aref{equ})}
For all real $x_1$, $x_2$ and $t$ satisfying
\begin{equation}\label{elef}
(1-e^{-3})T(x_1)\le t< T(x_1)\mbox{ and }|y|<1\mbox{ where }y=\frac{x_2-x_1}{T(x_1)-t},
\end{equation}
it holds that  
\[
|u(x_2,t)| \le K(x_1)(T(x_1)-t)^{-\frac 2{p-1}}(1-y^2)^{-\frac 1{p-1}},
\]
where $K(x_1)$ depends only on $p$ and an upper bound on $T(x_1)$ and $1/T(x_1)$.
\end{lem}
{\it Proof}: Take $x_1\in \R$. Using Proposition 3.5 page 66 in \cite{MZjfa07}, we see that for all  $s\ge -\log T(x_1)+3$, we have $\|w_{x_1}(s)\|_{\H_0}\le K(x_1)$, where $K(x_1)$ depends only on $p$ and an upper bound on $T(x_1)$ and $1/T(x_1)$.
 Using (i) of Lemma \ref{lemsobolev-9}, we see that 
\begin{equation}\label{carla}
\mbox{if }|y|<1,\mbox{ then }|w_{x_1}(y,s)|\le K(x_1)(1-y^2)^{-\frac 1{p-1}}.
\end{equation}
Now, taking $x_1$, $x_2$ and $t$ satisfying \aref{elef} and introducing $s=-\log(T(x_1)-t)$, we see that $s\ge -\log T(x_1)+3$ and the conclusion follows from \aref{carla} thanks to the selfsimilar transformation \aref{defw}.
\Box 

\bigskip

{\it Proof of Proposition \ref{phat}}:\\
(i) Using Proposition \ref{propref} and (i) of Lemma \ref{lemsobolev-9}, we see that we have
\begin{equation}\label{damien}
\sup_{|y|<1}\left|(1-y^2)^{\frac 1{p-1}}w_{x_0}(y,s)- \sum_{i=1}^{k(x_0)}e_i^*(1-y^2)^{\frac 1{p-1}}\kappa(d_i,y)\right|\to 0 \mbox{ as }s\to \infty.
\end{equation}
Since we have 
\begin{equation}\label{tkappa}
\kappa(d_i(s),y)(1-y^2)^{\frac 1{p-1}}= \kappa_0\cosh^{-\frac 2{p-1}}(\xi-\zeta_i(s))\mbox{ if }y=\tanh \xi
\end{equation}
and $\zeta_{i+1}(s)-\zeta_i(s) \to \infty$ as $s\to \infty$, we apply \aref{damien} with $y=-d_j(s)=\tanh \zeta_j(s)$ to get 
\begin{equation}\label{limj}
(1-d_j(s)^2)^{\frac 1{p-1}}w_{x_0}(-d_j(s),s) \to e_j^*\kappa_0\mbox{ as }s\to \infty.
\end{equation}
Since $(1-d_j(s)^2)^{-\frac 1{p-1}}= \cosh^{\frac 2{p-1}}\zeta_j(s)$, $e_j^* = e_1^*(-1)^{j+1}$ and $s\mapsto d_j(s)\in (-1,1)$ is continuous, using the selfsimilar transformation \aref{defw}, we see that (i) follows. 

\medskip

\noindent (ii) Let us introduce $B(x)$ by
\begin{equation}\label{defBB}
\frac{T(x_0)-T(x)}{|x_0-x|}=1-B(x).
\end{equation}
From the fact that $x\mapsto T(x)$ is $1$-Lipschitz, we see that
\begin{equation}\label{onelip}
\frac{|T(x_0)-T(x)|}{|x_0-x|}\le 1\mbox{ hence } 0\le B(x)\le 2,
\end{equation}
and the left-hand inequality of \aref{chapeau} follows. Using \aref{defBB}, we see that in order to prove the right-hand inequality, it is enough to prove that for $x_0-x$ small,
\begin{equation}\label{butB}
B(x) \le Cl^{-\frac{(k-1)(p-1)}2}\mbox{ where }l=|\log(x_0-x)|.
\end{equation}
To prove \aref{butB}, we consider only the case $x<x_0$ since the case $x>x_0$ follows by considering $u(-x,t)$ instead of $u(x,t)$. 
The key idea is to derive lower and upper bounds for $|u(z_1(\tt),\tt)|$ where $z_1(\tt)$ defined in \aref{defzj} is the ``center'' of the first soliton and 
\begin{equation}\label{defttt}
\tt=\tt(x)=T(x_0)-2(x_0-x).
\end{equation}
The following claim allows us to conclude:
\begin{cl}\label{cllubound}It holds that:\\
(i) $|u(z_k(\tt),\tt)| \ge \frac 1C(x_0-x)^{-\frac 2{p-1}}l^{\frac{k-1}2}$.\\
(ii) $|u(z_k(\tt),\tt)| \le C(x_0-x)^{-\frac 2{p-1}}B(x)^{-\frac 1{p-1}}$.
\end{cl}
Indeed, if Claim \ref{cllubound} holds, then we see that
\[
\frac 1C(x_0-x)^{-\frac 2{p-1}}l^{\frac{k-1}2} 
\le |u(z_k(\tt),\tt)| 
\le C(x_0-x)^{-\frac 2{p-1}}B(x)^{-\frac 1{p-1}}
\mbox{ hence } 
B(x)\le Cl^{-\frac{(k-1)(p-1)}2}
\]
and \aref{butB} follows. It remains to prove Claim \ref{cllubound} to conclude the proof of Proposition \ref{phat}.

\medskip

{\it Proof of Claim \ref{cllubound}}:\\
(i) Using (i) of Proposition \ref{phat} with $j=1$, we get for $x_0-x$ small enough,
\begin{equation}\label{ulower}
|u(z_k(\tt),\tt)| \ge  \frac{\kappa_0}2\cosh^{\frac 2{p-1}}\zeta_1(\ts)(T(x_0)-\tt)^{-\frac 2{p-1}}
\end{equation}
where 
\begin{equation}\label{defts}
\ts=-\log(T(x_0)-\tt)=-\log2 - \log(x_0-x).
\end{equation}
Recalling from (ii) of Proposition \ref{propref} that
\begin{equation}\label{estzeta1}
|\zeta_1(\ts)+\frac{(k-1)(p-1)}4\log(\ts)|\le C_0,
\end{equation}
hence, from \aref{defts}
\[
\cosh \zeta_1(\ts)\ge \frac {e^{-\zeta_1(\ts)}}2\ge C \ts^{\frac{(k-1)(p-1)}4}\ge Cl^{\frac{(k-1)(p-1)}4},
\]
the conclusion follows from \aref{ulower} and \aref{defttt}.\\
(ii) The idea is to apply Lemma \ref{pupper} for some well-chosen $x_1$, $x_2$ and $t$. We claim that condition \aref{elef} holds with $x_1=x$, $x_2=z_1(\tt)$ and $t=\tt$. Indeed:\\
- using \aref{onelip} and \aref{defttt}, we write for $x_0-x$ small enough,
\begin{equation}\label{o1}
T(x)\ge T(x_0)-(x_0-x) > T(x_0)-2(x_0-x) =\tt(x) \ge T(x_0)(1-e^{-3});
\end{equation}
- using \aref{defzj}, \aref{defttt}, \aref{estzeta1} and \aref{defts}, we write
\begin{eqnarray}
z_1(\tt)-x &=& x_0 +2(x_0-x)(-1 +2e^{2\zeta_1(\ts)}+O(e^{4\zeta_1(\ts)}))-x\nonumber\\
&= & (x_0-x)(-1+e^{O(1)}\ts^{-\frac{(k-1)(p-1)}2})\nonumber\\
&=&  (x_0-x)(-1+e^{O(1)}l^{-\frac{(k-1)(p-1)}2}).\label{souchon}
\end{eqnarray}
Since we have from \aref{onelip} and the definitions \aref{defBB} and \aref{defttt} of $B(x)$ and $\tt$,
\begin{equation}\label{o3}
T(x)-\tt =T(x)-T(x_0)+2(x_0-x)=(x_0-x)(1+B(x))\in [x_0-x,3(x_0-x)],
\end{equation}
we deduce that $y=\frac{z_1(\tt)-x}{T(x)-\tt}$ satisfies
\begin{equation}\label{o2}
y=\frac{-1+e^{O(1)}l^{-\frac{(k-1)(p-1)}2}}{1+B(x)}\in \left[-1+\frac{l^{-\frac{(k-1)(p-1)}2}}C, -\frac 13 +Cl^{-\frac{(k-1)(p-1)}2}\right].
\end{equation}
Thus, from \aref{o1} and \aref{o2}, condition \aref{elef} holds and Lemma \ref{pupper} applies and gives
\begin{equation}\label{upperb}
|u(z_1(\tt), \tt)|\le K(x)(T(x)-\tt)^{-\frac 2{p-1}}(1-y^2)^{-\frac 1{p-1}}
\end{equation}
where $K(x) \le K_0$ uniformly for $x$ close to $x_0$. 
Since we have from \aref{o2} and \aref{onelip},
\begin{equation}\label{o4}
(1-y^2) \ge \frac 1C(1+y)=\frac 1C\frac{B(x)+e^{O(1)}l^{-\frac{(k-1)(p-1)}2}}{1+B(x)} \ge \frac {B(x)}C,
\end{equation}
the conclusion follows from \aref{o4}, \aref{o3} and \aref{upperb}. This concludes the proof of Claim \ref{cllubound} and Proposition \ref{phat} too. \Box

\section{Properties of $\SS$}\label{sec2}
We proceed in 3 subsections. We first prove that the interior of $\SS$ is empty which is the desired conclusion of Proposition \ref{thfini}. Then, we give the proofs of Theorems \ref{thnonexis}, \ref{thnonexis}' and \ref{new} as well as Proposition \ref{pdes}. 

\subsection{Soliton characterization on $\SS$}
This subsection is devoted to the proof of the following result:
\begin{prop}\label{pbut}$ $\\
(i) The interior of $\SS$ is empty.\\
(ii) For all $x_0\in \SS$, $k(x_0)\ge 2$.
\end{prop}
Before proving this proposition, let us first state the following Lemmas:
\begin{lem}\label{cns}{\bf (Characterization of the interior of $\SS$)}
For any $x_1<x_2$, the following statements are equivalent:

(a) $(x_1, x_2) \in \SS$.

(b) There exists $x^*\in [x_1, x_2]$ such that for all $x\in [x_1, x_2]$, $T(x) = T(x^*) - |x-x^*|$.
\end{lem}
\begin{lem}\label{maxslope}
Consider $x_1<x_2$ such that $e\equiv \frac{T(x_2)-T(x_1)}{x_2-x_1}=\pm 1$. Then,\\
(i) for all $x\in [x_1, x_2]$, $T(x)=T(x_1)+e(x-x_1)$,\\
(ii) $(x_1, x_2) \in \SS$.
\end{lem}
\begin{lem}[Boundary properties of $\SS$]\label{lemfrank}$ $\\
(i) For all $x_0\in \partial \SS$, $k(x_0) \neq 0$.\\
(ii) Consider $x_0\in \partial \SS$ with $k(x_0)=1$.
If there exists a sequence $x_n\in \RR$ converging from the left (resp. the right) to $x_0$, then $x_0$ is left-non-characteristic (resp. right-non-characteristic).
\end{lem}
{\bf Remark}: We mean by $x_0$ is left-non-characteristic (resp. right-non-characteristic) that it satisfies condition \aref{nonchar} only for $x<x_0$ (resp. for $x>x_0$).\\
{\bf Remark}: It is not possible to prove by a direct argument that $k(x_0) \ge 1$ when $x_0$ is arbitrary in $\SS$. We need to prove it first for $x_0\in\partial \SS$ and then prove that the interior is empty. See the derivation of Proposition \ref{pbut} from Lemma \ref{lemfrank}.

\bigskip

We now give the proofs of Lemmas \ref{cns}, \ref{maxslope} and \ref{lemfrank}.

\medskip

{\it Proof of Lemma \ref{cns}}:

(a) $\Longrightarrow$ (b): Let us introduce $x^*\in [x_1, x_2]$ such that
\[
T(x^*)=\max_{x_1\le x\le x_2} T(x).
\]
We claim that $T(x)$ is nondecreasing on $[x_1,x^*]$, and nonincreasing on $[x^*, x_2]$. Indeed, let us prove the first fact, the second being similar. If for some $x'\le x''$ in $[x_1, x^*]$, we have $T(x')>T(x'')$, then $\min_{x'\le x\le x^*}T(x)\le T(x'') <T(x') \le T(x^*)$. Therefore, this minimum is achieved at a point $\tilde x$ different from $x'$ and $x^*$, hence
\[
\tilde x\in (x', x^*) \subset (x_1, x_2).
\]
In other words, $\tilde x$ is a local minimum, hence non characteristic, which is in contradiction with (a).\\
The result clearly follows if we prove that
\begin{eqnarray}
\forall x\in (x_1, x^*),&&T(x) = T(x_1)+(x-x_1),\label{triest}\\
\forall x\in (x^*, x_2),&&T(x)= T(x_2)-(x-x_2).\label{triouest}
\end{eqnarray}
We only prove \aref{triest} since \aref{triouest} follows similarly.\\
Assume by contradiction that for some $x'\in (x_1, x^*)$, we have  
\begin{equation}\label{contr}
\d\frac{T(x')-T(x_1)}{x'-x_1}=m_0\not\in \{-1,1\}.
\end{equation}
 Then, since $x\mapsto T(x)$ is 1-Lipschitz and nondecreasing, it follows that $0\le m_0<1$. Considering a family of lines of slope $\frac{1+m_0}2$ growing from below, we find $\lambda_0 \in \R$ and $x_0\in [x_1,x']$ such that
\begin{equation}\label{nazih}
\forall x\in [x_1,x'],\;\;T(x) \ge \frac{(1+m_0)}2(x-x_1)+\lambda_0\mbox{ and }T(x_0)= \frac{(1+m_0)}2(x_0-x_1)+\lambda_0.
\end{equation}
If $x_0\in (x_1,x')$, then for all $x\in [x_1, x']$, $T(x)\ge \frac{(1+m_0)}2(x-x_0)+T(x_0)$, hence $x_0$ is non characteristic (the cone of slope $\frac{1+m_0}2$ is convenient).\\
If $x_0=x'$, then since $T(x)$ is non decreasing on $(x_1, x^*)$, it follows that $x_0$ is again non characteristic. In these two cases, we have a contradiction with the fact that $(x_1, x_2) \in \SS$.\\
If $x_0=x_1$, then we have from \aref{nazih}, $T(x_1)= \lambda_0$ and $T(x')\ge \frac{(1+m_0)}2(x'-x_1)+T(x_1)$, in contradiction with \aref{contr}.\\
Thus, \aref{triest} holds. Since \aref{triouest} follows similarly, (b) follows too. 

\medskip

(b) $\Longrightarrow$ (a): For any $x\in (x_1, x_2)$, the left-slope of $x\mapsto T(x)$ is $1$ or $-1$, hence, by definition, $x\in \SS$ and (a) follows. 
This concludes the proof of Lemma \ref{cns}.\Box

\bigskip

We now give the proof of Lemma \ref{maxslope}:

{\it Proof of Lemma \ref{maxslope}}: Up to replacing $u(x,t)$ by $u(-x,t)$, we can assume that $x_1<x_2$ and
\begin{equation}\label{nej}
e\equiv \frac{T(x_2)-T(x_1)}{x_2-x_1}=1. 
\end{equation}
(i) If $x\in (x_1, x_2)$, we use the fact that $x\mapsto T(x)$ together with \aref{nej} to write:
\begin{eqnarray*}
T(x) &\le & T(x_1)+(x-x_1),\\
T(x) &\ge & T(x_2)-(x_2-x)=T(x_1)+(x-x_1)
\end{eqnarray*}
and (i) follows.\\
(ii) It follows from (i), just by applying the fact that (b) implies (a) in Lemma \ref{cns} (take $x^*=x_2$). This concludes the proof of Lemma \ref{maxslope}.\Box 

\bigskip

We now give the proof of Lemma \ref{lemfrank}.

\medskip

{\it Proof of Lemma \ref{lemfrank}}: 
Consider $x_0\in \partial \SS$. Up to replacing $u(x,t)$ by $u(-x,t)$, we can assume that for some sequence
\begin{equation}\label{hypon}
x_n\in \RR,\mbox{ we have }x_n<x_0\mbox{ and }x_n \to x_0\mbox{ as }n\to \infty.
\end{equation}
Therefore, we have
\begin{equation}\label{espace0}
\forall x< x_0,\;\;T(x)>T(x_0)-(x_0-x).
\end{equation}
Indeed, note first form the fact that $x\mapsto T(x)$ is $1$-Lipschitz that for all $x<x_0$, $T(x)\ge T(x_0)-(x_0-x)$. By contradiction, if for some $\hat x<x_0$, we have $T(\hat x)=T(x_0)-(x_0-\hat x)$, then we see from Lemma \ref{maxslope} that $(\hat x, x_0)\subset \SS$, in contradiction with \aref{hypon}. Thus, \aref{espace0} holds.

\medskip

To prove (i) and (ii), we proceed  by contradiction. Assume then that

\medskip

{\it \noindent 
either $k(x_0)=0$ ({\bf case 1}),\\
or $k(x_0)=1$ and $x_0$ is not left-non-characteristic ({\bf case 2}).}

\medskip

\noindent 
Using \aref{cprofile0} when $k(x_0)=0$ and Proposition \ref{thtrap} when $k(x_0)=1$, we see that 
\begin{eqnarray}
&&\left\|\vc{w_{x_0}(s)}{\partial_s w_{x_0}(s)} - \vc{w_\infty}0\right\|_\H\to 0\mbox{ as }s\to \infty,\label{carprofile01}\\
\mbox{where}&& w_\infty(y)=0\mbox{ if }k(x_0)=0
\mbox{ and }w_\infty(y)=e^* \kappa(d(x_0),y)\mbox{ if }k(x_0)=1,\label{defg}
\end{eqnarray} 
for some $e^*=\pm 1$ and $d(x_0)\in (-1,1)$. Now, we claim the following continuity result:
\begin{cl} \label{clrakia}
For all $\epsilon_0>0$, there exists $\tilde t<T(x_0)$ and $\tilde x<x_0$ such that for all $x'\in (\tilde x, x_0)$,
\begin{equation}\label{hadaf}
\left\|\vc{w_{x'}(\tilde s_0(x'))}{\partial_s w_{x'}(\tilde s_0(x'))} - \vc{w_\infty}0\right\|_\H \le \epsilon_0
\end{equation}
where $\tilde s_0(x') = - \log(T(x')-\tilde t)$.
\end{cl}
{\it Proof}: See Appendix \ref{appclrakia}.\Box 

Let us first use this lemma to find a contradiction.

\medskip

{\bf Case 1: $k(x_0)=0$}. Consider some $\epsilon_0>0$ (to be fixed small enough later). Using this claim, \aref{defg} and \aref{hypon}, we see that for some $\tilde t<T(x_0)$ and for $n$ large enough, we have 
\[
x_n \in \RR\mbox{ and }\left\|\vc{w_{x_n}(\tilde s_0(x'))}{\partial_s w_{x_n}(\tilde s_0(x'))}\right\|_\H \le \epsilon_0.
\]
Using  the continuity of $E(w)$ in $\H$ (which is a consequence of Lemma \ref{lemsobolev-9}), we see that
\[
E(w_{x_n}(\tilde s_0(x')))\le C\epsilon_0\le \frac 12 E(\kappa_0)
\]
 (if $\epsilon_0$ is small enough), on the one hand. On the other hand, since $x_n\in \RR$, we know from the limit and the monotonicity of $E(w_{x_n}(s))$ stated in page \pageref{threg} that $E(w_{x_n}(s))\ge E(\kappa_0)>0$, which is a contradiction.

\medskip

{\bf Case 2: $k(x_0)=1$ and $x_0$ is not left-non-characteristic}. Since $x_0$ is not left-non-characteristic, we see from \aref{espace0} that there exists a sequence  $\hat x_n$ such that 
\begin{equation}\label{defmnt}
\hat x_n< x_0,\;\;\hat x_n\to x_0\mbox{ and }\hat m_n\equiv \frac{T(\hat x_n)-T(x_0)}{\hat x_n -x_0}\in [1-\frac 1n, 1).
\end{equation}
Considering a family of lines of slope $\frac{1+\hat m_n}2$, we can select one such that
\begin{equation}\label{leilat}
\forall x\in [\hat x_n, x_0],\;\;T(x)\ge \frac{(1+\hat m_n)}2(x-\hat x_n)+\lambda_n\mbox{ and }T(\tilde x_n)= \frac{(1+\hat m_n)}2(\tilde x_n-\hat x_n)+\lambda_n
\end{equation}
for some $\lambda_n\in \R$ and $\tilde x_n \in [\hat x_n, x_0]$.\\
If $\tilde x_n=x_0$, then for all $x\in [\hat x_n, x_0]$, $T(x) \ge \frac{(1+\hat m_n)}2(x-x_0)+T(x_0)$, which is in contradiction with the fact that $x_0$ is left-non-characteristic.\\
If $\tilde x_n = \hat x_n$, then we have from \aref{leilat}, $T(\hat x_n) = \lambda_n$ and $T(x_0) \ge \frac{1+\hat m_n}2(x_0-\hat x_n)+T(\hat x_n)$, in contradiction with \aref{defmnt}.\\
If $\tilde x_n \in (\hat x_n, x_0)$, then $\tilde x_n \in \RR$ (the cone of slope $\frac{1+\hat m_n}2$ is convenient). Since $\tilde x_n\to x_0$ and $\tilde x_n \in\RR$, we see from Claim \ref{clrakia} that for some $\tilde t<T(x_0)$ and $n$ large enough, we have
\begin{equation}\label{ebde}
\left\|\vc{w_{\tilde x_n}(\tilde s_0(\tilde x_n))}{\partial_s w_{\tilde x_n}(\tilde s_0(\tilde x_n))} - e^*\vc{\kappa(d(x_0),\cdot)}0\right\|_\H \le \epsilon^*
\end{equation}
where $\epsilon^*$ is introduced in Proposition \ref{thtrap}. Since the energy barrier follows from the fact that $\tilde x_n\in \RR$, Proposition \ref{thtrap} applies and we have for $n$ large enough, 
\[
(w_{\tilde x_n}(s), \partial_s w_{\tilde x_n}(s)) \to e^* (\kappa(d_n), 0)\mbox{ in }\H\mbox{ as }s\to \infty,
\]
where $|d_n-d(x_0)|<\eta_0$ for some $\eta_0>0$ small enough so that $|d(x_0)\pm \eta_0|<1$. The use of the geometrical interpretation of $d_n$ is crucial for the conclusion. Indeed, from \aref{ebde} and the regularity result of \cite{MZjfa07} cited in page \pageref{threg}, we see that $x\mapsto T(x)$ is differentiable at $x=\tilde x_n$ and that
\[
T'(\tilde x_n)= d_n \le d(x_0)+\eta_0<1
\]
on the one hand. On the other hand, using \aref{leilat} and \aref{defmnt}, we see that
\[
T'(\tilde x_n)=\frac{1+\hat m_n}2 \to 1\mbox{ as } n\to \infty
\]
which is a contradiction. This concludes the proof of Lemma \ref{lemfrank}.\Box

\medskip

Now, we are ready to prove Proposition \ref{pbut}.

\medskip

{\it Proof of Proposition \ref{pbut}}:\\
(i) Let us assume by contradiction that $\SS$ contains some non empty interval $(a', b')$. Since $\SS \neq \R$ by the result of \cite{MZjfa07} cited in page \pageref{threg}, by maximizing this interval and up to replacing $u(x,t)$ by $u(-x,t)$, we can assume that:
\[
(a,b) \subset \SS \mbox{ with }a \in \partial \SS,\;\;b>a
\]
and, either $b\in \partial \SS$ or $b=+\infty$. If $b$ is finite, then up to replacing $u(x,t)$ by $u(-x,t)$, we can assume that $T(b) \ge T(a)$.
Using Lemma \ref{cns} and the fact that $T(x) \ge 0$, we see that for some $\tilde b<b$, we have
\begin{equation}\label{line}
\forall x\in (a,\tilde b),\;\;T(x)=T(a)+(x-a).
\end{equation}
We consider three cases and find a contradiction in each case.\\
- If $k(a)=0$, then a contradiction occurs from (i) of Lemma \ref{lemfrank}.\\
- If $k(a) =1$, then from the fact that $a\in \partial \SS$, there exists a sequence $x_n \in \RR$ such that $x_n \to a$ as $n\to \infty$. Since $(a,b)\subset \SS$, it follows that $x_n <a$ for $n$ large enough. Therefore, applying (ii) of Lemma \ref{lemfrank}, we see that $a$ is left-non-characteristic. Since it is clearly right-non-characteristic by \aref{line},  $a$ is in fact non-characteristic, which contradicts the fact that $a\in \partial \SS \subset \SS$ (note that $\SS$ is closed since its complementary set $\RR$ is open by the result of \cite{MZcmp08} cite in page \pageref{threg}).\\
- If $k(a) \ge 2$, then the corner property stated in Proposition \ref{propref} is in contradiction with \aref{line}.\\
Thus, (i) of Proposition \ref{pbut} follows.

\medskip

\noindent (ii) Consider $x_0\in \SS$. From (i), we have $x_0 \in \partial\SS$.
Using (i) of Lemma \ref{lemfrank}, we see that $k(x_0)\ge 1$. The result follows if we rule out the case $k(x_0)=1$.\\
 Assume by contradiction that $k(x_0)=1$. Since the interior of $\SS$ is empty, we can construct 2 sequences $x_n$ and $y_n$ in $\RR$, such that $x_n\to x_0$ from the left, and $y_n\to x_0$ from the right. Applying (ii) of Lemma \ref{lemfrank}, we see that $x_0$ is in fact left-non-characteristic and right-non-characteristic, hence non characteristic. This contradicts the fact that $x_0\in \SS$. Thus, (ii) follows. This concludes the proof of Proposition \ref{pbut}\Box

\subsection{On characteristic points for equation \aref{equ}}
We prove Proposition \ref{pdes}, as well as Theorems \ref{new} and \ref{thnonexis} here.

\medskip

{\it Proof of Theorem \ref{new}}: Consider $u(x,t)$ a solution of equation \aref{equ} and $x_0\in \SS$. Using  (ii) of Proposition \ref{pbut}, we see that  $k(x_0) \ge 2$. Therefore,  Proposition \ref{propref} applies and directly gives the conclusion of Theorem \ref{new}. 

\medskip

{\it Proof of Proposition \ref{pdes}}:  Consider $u(x,t)$ a solution of equation \aref{equ} and $x_0\in \SS$. Using  (ii) of Proposition \ref{pbut}, we see that  $k(x_0) \ge 2$. Therefore, Proposition \ref{propref} and Proposition \ref{phat} apply and give the conclusion of Proposition \ref{pdes}, except for the strict inequality in \aref{chapeau0}, which we prove now.\\
Assume by contradiction that for some $x_1< x_0$, we have equality in the left-hand inequality of \aref{chapeau0}. Then, we see from Lemma \ref{maxslope} that $(x_1, x_0) \subset \SS$, which contradicts the fact that the interior of $\SS$ is empty (see (i) of Proposition \ref{pbut}). Thus, \aref{chapeau0} follows and Proposition \ref{pdes} follows too.\Box

\medskip

{\it Proof of Theorem \ref{thnonexis}}: Consider $u(x,t)$ a solution of equation \aref{equ} that blows up on a graph $x\mapsto T(x)$ such that for some $a_0<b_0$ and some $t_0\ge 0$, we have 
\begin{equation}\label{hypp}
\forall x\in (a_0, b_0)\mbox{ and }t\in [t_0, T(x)),\;\; u(x,t)\ge 0.
\end{equation}
 We would like to prove that $(a_0, b_0)\subset \RR$. Proceeding by contradiction, we assume that there exists $x_0\in (a_0, b_0)\cap \SS$. Using Proposition \ref{pdes}, we see that for some $e_1=\pm 1$ and $t_1\in [t_0, T(x_0))$, there are continuous $t\mapsto z_i(t)$ where $i=1$ and $2$ such that $z_i(t) \to x_0$ as $t\to T(x_0)$ and
\[
\forall t\in [t_1, T(x_0)),\;\;e_1 u(z_1(t), t)>0\mbox{ and } e_1u(z_2(t), t)<0.
\]
Therefore, $u$ changes sign in $(a_0, b_0)\times [t_1, T(x_0))$, which is in contradiction with \aref{hypp}. Thus, $(a_0, b_0)\subset \RR$ and Theorem \ref{thnonexis} follows.\Box

\subsection{Non existence of characteristic points for equation \aref{equ'}}
This subsection is dedicated to the study of equation \aref{equ'} which we recall here:
\begin{equation}\label{equ2}
\left\{
\begin{array}{l}
\partial^2_{t t} u =\partial^2_{xx} u+|u|^p,\\
u(0)=u_0\mbox{ and }u_t(0)=u_1.
\end{array}
\right.
\end{equation}
We take $p>1$ and $(u_0, u_1) \in \h1\times \l2$.
Our aim is to prove Theorem \ref{thnonexis}', which asserts that the set of characteristic points is empty, for any blow-up solution of \aref{equ2}. To do so, we need to perform for equation \aref{equ2}, an almost identical analysis to what we did for equation \aref{equ}, in our previous papers, including this last one. Therefore, we only give the main steps and stress only the novelties.

\medskip

Consider $u(x,t)$ a solution of equation \aref{equ'} that blows up on some graph $x\mapsto T(x)$.  As for equation \aref{equ}, we denote the set of non characteristic points by $\RR$ and the set of characteristic points by $\SS$. Our aim is to show that $\SS =\emptyset$.\\
Given $x_0\in \R$ and $T_0\in (0, T(x_0)]$, we define $w_{x_0, T_0}$ as in \aref{defw} by 
\begin{equation}\label{defw2}
w_{x_0,T_0}(y,s)=(T_0-t)^{\frac 2{p-1}}u(x,t),\;\;y=\frac{x-x_0}{T_0-t},\;\;
s=-\log(T_0-t).
\end{equation}
If $T_0=T(x_0)$, then we simply write $w_{x_0}$ instead of $w_{x_0, T(x_0)}$. 
The function $w=w_{x_0,T_0}$ satisfies the 
following equation for all $y\in B=B(0,1)$ and $s\ge -\log T_0$:
\begin{equation}\label{eqw'}
\partial^2_{ss}w= \L w-\frac{2(p+1)}{(p-1)^2}w+|w|^p
-\frac{p+3}{p-1}\partial_sw-2y\partial^2_{y,s} w.
\end{equation} 
The Lyapunov functional for this equation is defined in $\H$ \aref{defnh0} and is given by
\[
\tilde E(w(s))= \iint \left(\frac 12 \left(\partial_s w\right)^2 + \frac 12  \left(\partial_y w\right)^2 (1-y^2)+\frac{(p+1)}{(p-1)^2}w^2 - \frac 1{p+1} |w|^pw\right)\rho dy
\]
and satisfies 
\[
\frac d{ds}\tilde E(w(s)) = -\frac 4{p-1}\iint  \left(\partial_s w(y,s)\right)^2 \frac{\rho}{1-y^2}dy.
\]
We first give the following lower bound:
\begin{lem}{\bf (A lower bound on solutions of \aref{equ'})} \label{lembi} For all $R>0$, there exists $M(R)>0$ such that for all $x\in (-R, R)$ and $t\in [0, T(x))$, we have $u(x,t) \ge -M(R)$.
\end{lem}
{\it Proof}: 
 Using Duhamel's formula, we write for all $x\in \R$ and $t\in [0, T(x))$, 
\begin{equation}\label{duhamel}
u(x,t)=S(t)u_0(x)+S_1(t) u_1(x)+\int_0^t(S_1(t-\tau)|u(\tau)^p|)(x) d\tau
\end{equation}
where
\[
S(t)h(x) = \frac 12(h(x+t)+h(x-t))\mbox{ and }S_1(t)h(x)=\frac 12 \int_{x-t}^{x+t}h(x')dx'.
\]
Take $R>0$ and introduce $R_0= R+\max_{|x|\le R}T(x)$. Since we have by hypothesis,
\[
u_0\in H^1(-R_0, R_0) \subset L^\infty(-R_0, R_0)\mbox{ and }u_1\in L^2(-R_0, R_0),
\]
we use the continuity of $x\mapsto T(x)$ to write from \aref{duhamel}:
\[
u(x,t)\ge S(t)u_0(x)+S_1(t) u_1(x) \ge -\|u_0\|_{L^\infty(-R_0, R_0)} - \sqrt{R_0}\|u_1\|_{L^2(-R_0, R_0)}, 
\]
which concludes the proof of Lemma \ref{lembi}.\Box

\bigskip

Using Lemma \ref{lembi}, we get the following consequences:
\begin{cl}\label{lempert} {\bf (A lower bound on solutions of \aref{eqw'})}If $x_0\in \R$, and $T_0\in (0, T(x_0)]$, then:\\
(i) For all $y\in (-1,1)$ and $s\ge -\log T_0$, we have $w(y,s) \ge -M_0 e^{-\frac {2s}{p-1}}$, where $M_0=M(|x_0|+T(x_0))$.\\
(ii) For all $s\ge -\log T_0$, 
\begin{eqnarray*}
\frac 12 \iint  |w|^{p+1} \rho dy - CM_0^{p+1}e^{-\frac {2(p+1)s}{p-1}}&\le& \iint  |w|^pw \rho dy \le \iint  |w|^{p+1} \rho dy,\\
(1-2M_0e^{-\frac {2s}{p-1}})\iint  |w|^{p+1} \rho dy - CM_0e^{-\frac {2s}{p-1}}&\le& \iint  |w|^pw \rho dy.
\end{eqnarray*}
\end{cl}
{\it Proof}:\\
 (i) It follows straightforwardly from Lemma \ref{lembi}.\\
(ii) Consider $s\ge -\log T_0$. The right-hand side of the first line is obvious. For the left-hand side inequality of both lines, we use (i) to write 
\[
|z|^pz-|z|^{p+1}\ge -\epsilon_0 |z|^p \ge \max(-\frac 12z^{p+1} -2^p\epsilon_0^{p+1}, - \epsilon_0(1+z^{p+1}))
\]
where $z=w(y,s)$ and $\epsilon_0 = 2M_0e^{-\frac {2s}{p-1}}$. By integration, (ii) follows.\Box

\bigskip

In the following, we give the following blow-up criterion for equation \aref{equ'}:
\begin{cl}[A blow-up criterion for equation \aref{equ'}]Consider $W(y,s)$ a solution to equation \aref{eqw'} locally bounded in $L^{p+1}(-1,1)$ such that 
$\tilde E(W(s_0))<0$ for some $s_0\in \R$. Then, $W(y,s)$ cannot exist for all $(y,s)\in (-1,1)\times[s_0, \infty)$.
\end{cl}
{\it Proof}: The proof is the same as the proof of Theorem 2 in Antonini and Merle \cite{AMimrn01} (of course, one need to use the Lyapunov functional $\tilde E(w)$). 
\Box  

\bigskip

Using the Lyapunov functional $\tilde E(w)$ together with the estimate in (ii) of Claim \ref{lempert}, one can adapt with no difficulty the analysis of our previous papers (\cite{MZajm03}, \cite{MZimrn05}, \cite{MZjfa07} and \cite{MZcmp08}, without forgetting the present paper) to equation \aref{equ'}, and get the same results, with the following new feature:\\
Due to the lower bound of (i) in Claim \ref{lempert}, only nonnegative objects appear in the limit at infinity of $w_{x_0}$ when $x_0\in \RR$ (take $\theta(x_0)=1$ in \aref{profile} and \aref{profilelebesgue}) and in the asymptotic decomposition of $w_{x_0}$ when $x_0\in \SS$ (take $e_i=1$ for all $i=1,...,k$ in \aref{cprofile0}). This is the main difference with the case of equation \aref{equ}, where different signs may appear. More precisely, we have the following:
\begin{cl}{\bf (Classification of the nonnegative stationary solutions of equation \aref{eqw'})} Consider $w\in \H_0$ a nonnegative stationary solution of \aref{eqw'}. Then
\begin{equation}\label{tzev}
\mbox{ either }w\equiv 0\mbox{ or }w(y,s) = \kappa(d,y)
\end{equation}
for some $d\in(-1,1)$, 
 where $\kappa(d,y)$ is defined in \aref{defkd}.
\end{cl}
{\bf Remark}: It is easy to prove that all stationary solutions of \aref{eqw'} in $\H_0$ are in fact nonnegative, hence characterized by \aref{tzev}.\\
{\it Proof}: If $w\in \H_0$ a nonnegative stationary solution of \aref{eqw'}, then it is also a stationary solution of \aref{eqw}. Using Proposition 1 of \cite{MZcmp08}, we see that either $w\equiv 0$ or $w(y,s) = e\kappa(d,y)$ for some $d\in(-1,1)$ and $e=\pm 1$. Since $w$ is nonnegative, we get $e=1$. \Box


\bigskip

Arguing as in \cite{MZcmp08} (note that the Liouville Theorem 2 and 2' of \cite{MZcmp08} hold for equation \aref{equ'} and that only nonnegative solutions are possible), we get that the set $\RR$ of non characteristic points is non empty, open and $x\mapsto T(x)$ is $C^1$ on $\RR$ (see page \pageref{threg} in this paper, see Theorem 1 and the following remark in \cite{MZcmp08}). In other words, the set $\SS$ of characteristic points is closed and 
\[
\partial \SS \subset \SS.
\]

\medskip

As a consequence of the fact that only nonnegative solitons appear in the asymptotic decomposition \aref{cprofile0} of $w_{x_0}$ when $x_0\in \SS$, we have the following result which is the main difference with equation \aref{equ}:
\begin{cl}\label{clk1}$ $\\
(i) For all $x_0\in \SS$, $k(x_0)=0$ or $1$.\\
(ii) For all $x_0\in \partial \SS$, $k(x_0)=1$.\\
(iii) Consider $x_0\in \partial \SS$.
If there exists a sequence $x_n\in \RR$ converging from the left (resp. the right) to $x_0$, then $x_0$ is left-non-characteristic (resp. right-non-characteristic).
\end{cl}
{\it Proof}:\\ 
(i) Proceeding by contradiction, we assume that for some $x_0\in \SS$, we have $k(x_0)\ge 2$. As for equation \aref{equ}, $w_{x_0}$ can be decomposed as $s\to \infty$ as a sum of decoupled solitons (take 
\begin{equation}\label{farid}
e_i=1\mbox{ for all }i=1,...,k,
\end{equation} 
in \aref{cprofile0}), and we can show that (up to slightly changing the solitons' location), the solitons' centers satisfy the same ODE system \aref{eqz} as in the case of equation \aref{equ}.\\  
Therefore, (i) of Proposition \ref{propref} holds and we have
\[
\forall i=1,...,k,\;\;e_i=(-1)^{i+1} e_1.
\]
Since $k(x_0)\ge 2$, this is in contradiction with \aref{farid}. Thus, (i) holds.\\
(ii) Using Lemma \ref{lemfrank}, we see that for all $x_0\in \partial \SS$, $k(x_0)\neq 0$. Using (i), we get the conclusion.\\
(iii) Consider $x_0\in \partial \SS$. From (ii), we have $k(x_0)=1$. Applying (ii) of Lemma \ref{lemfrank}, we get the result. This concludes the proof of Claim \ref{clk1}.\Box

\bigskip

Now, we are ready to give the proof of Theorem \ref{thnonexis}'.

\medskip

{\it Proof of Theorem \ref{thnonexis}'}: Assume by contradiction that $\SS\neq \emptyset$. Since $\RR\neq \emptyset$ (see the remark following Theorem 1 in \cite{MZcmp08}), it follows that $\partial \SS \neq \emptyset$. If $x_0\in \partial \SS$, then up to replacing $u(x,t)$ by $u(-x,t)$, we assume that there exists a sequence $x_n \in \RR \to x_0$ from the left as $n\to \infty$. Applying Claim \ref{clk1}, we see that 
\begin{equation}\label{ahmed}
k(x_0)= 1\mbox{ and }x_0\mbox{ is left-non-characteristic}.
\end{equation} 
Now, we consider 2 cases.\\
- If $[x_0, \infty)\subset \SS$, then we have from Lemma \ref{cns} and the positivity of $T(x)$ that
\[
\forall x\ge x_0, T(x)= T(x_0)+(x-x_0).
\]
Therefore, $x_0$ is right-non-characteristic, hence non characteristic. Contradiction with the fact that $x_0\in \partial \SS \subset \SS$.\\
- Now, if $[x_0, \infty)\not\subset \SS$, then we can define $x_1\ge x_0$ maximal such that 
\[
[x_0, x_1] \subset \SS.
\]
Since $x_1$ is maximal, it follows that $x_1\in \partial \SS$ and that there exists a sequence $y_n \in \RR \to x_1$ from the right as $n\to \infty$. Applying Claim \ref{clk1}, we see that 
\begin{equation}\label{ahmed2}
k(x_1)= 1\mbox{ and }x_1\mbox{ is right-non-characteristic}.
\end{equation} 
If $x_1=x_0$, then $x_0$ is non characteristic by \aref{ahmed}, which is a contradiction.\\
If $x_1>x_0$, then applying Lemma \ref{cns}, we see that for some $x^* \in [x_0, x_1]$, we have
\begin{equation*}
\forall x\in [x_0, x_1],\;\;T(x) = T(x^*) - |x-x^*|.
\end{equation*}
If $x^*>x_0$, then $x_0$ is right-non-characteristic, hence non characteristic by \aref{ahmed}. Contradiction again.\\
If $x^*=x_0$, then $x_1$ is left-non-characteristic, hence non characteristic by \aref{ahmed2}. This is in contradiction with the fact that $x_1\in \partial \SS\subset \SS$.\\
This concludes the proof of Theorem \ref{thnonexis}'.\Box

\appendix

\section{Continuity with respect to initial data of the blow-up time at a non characteristic point}
\label{appcont}
This section is devoted to the proof of Proposition \ref{propcontid}. The proof is more or less included in the arguments of the proof of Lemma 2.2 of \cite{MZcmp08}. We only give here a sketch of the proof (see \cite{MZcmp08} for more details).

\medskip

{\it Sketch of the proof of Proposition \ref{propcontid}}: We will prove the continuity in the norm $H^1\times L^2(\R)$ since the result with the norm $H^1\times L^2(|x|<A_0)$ follows from the finite speed of propagation.\\
Using the continuity of $u(t_0)$ for $t_0<T(x_0)$ with respect to initial data, it follows that $T(x_0)$ is lower semi-continuous as a function of initial data.\\
For the upper continuity, we consider $T_0>T(x_0)$ to be taken close enough to $T(x_0)$ and aim at proving that $\tilde u(x,t)$ blows up in finite time $\tilde T(x_0)<T_0$, where $\tilde u(x,t)$ is the solution of equation \aref{equ} with initial data $(\tilde u_0, \tilde u_1)$ close enough to $(u_0, u_1)$.

\medskip

\noindent Up to changing $u$ in $-u$, we know from \aref{profilelebesgue} that for some $d_0\in (-1,1)$ and $\delta_0>0$, 
\begin{equation}\label{profilelarge}
\left\|\vc{w_{x_0}(s)}{\partial_s w_{x_0}(s)}-\vc{\kappa(d_0, .)}{0}\right\|_{H^1\times L^2(|y|<1+\delta_0)}\to 0 \mbox{ as }s\to \infty.
\end{equation}
Consider $s_0<0$ to be fixed later and introduce $t_0<T(x_0)$ such that
$\frac{T(x_0)-t_0}{T_0-t_0}=1-e^{s_0}$.
Using the selfsimilar transformation \aref{defw} and \aref{profilelarge}, we see that 
\begin{eqnarray}
&&w_{x_0, T_0}(y,-\log(T_0-t_0)) = (1-e^{s_0})^{-\frac 2{p-1}}w_{x_0}\left(\frac y{1-e^{s_0}}, -\log(T(x_0)-t_0)\right),\nonumber\\
&&\left\|\vc{w_{x_0, T_0}(-\log(T_0-t_0))}{\partial_s w_{x_0, T_0}(-\log(T_0-t_0))}-\vc{w_-(s_0)}{\partial_s w_-(s_0)}\right\|_{H^1\times L^2(-1,1)}\to 0\label{conv}
\end{eqnarray}
as $T_0\to T(x_0)$, where $w_-(y,s) = \kappa_0 \frac{(1-d_0^2)^{\frac 1{p-1}}}{(1-e^{s_0}+ d_0 y)^{\frac 2{p-1}}}$ is a particular solution of equation \aref{eqw}.\\
Since $E(w_-(s_0))<0$ for some $s_0<0$ from Appendix B in \cite{MZcmp08}, we see from \aref{conv} that for $T_0$ close enough to $T(x_0)$ and $t_0$ defined 
above,
we have 
\begin{equation}\label{Eneg}
E(w_{x_0, T_0}(-\log(T_0-t_0)))<0.
\end{equation}
Using the blow-up criterion of Antonini and Merle (see Theorem 2 in \cite{AMimrn01}), we see that $w_{x_0, T_0}$ cannot be defined for all $(y,s)\in (-1,1)\times [-\log T_0, \infty)$, which means that $u$ blows up in finite time and that $T(x_0)<T_0$. This yields the upper semi-continuity and concludes the proof of Proposition \ref{propcontid}.\Box
\section{Estimates on the quadratic form $\varphi$}\label{appdecomp}
This section is devoted to the proof of Lemma \ref{lemdecomp}. We proceed in two subsections:\\
- in the first subsection, we give some preliminary results, in particular, we change the problem to the $\xi$ variable, where $y=\tanh \xi$.\\
- in the second subsection, we give the proof of Lemma \ref{lemdecomp}.

\subsection{Preliminaries and formulation in the $\xi$ variable with $y=\tanh \xi$}
We first recall the following result from \cite{MZjfa07}.
\begin{cl}\label{lemsobolev-9} 
(i) {\bf (A Hardy-Sobolev type identity)} For all $h\in \H_0$, it holds that
\[
\|h\|_{L^2_{\frac \rho{1-y^2}}}+\|h\|_{L^{p+1}_\rho}+ \|h(1-y^2)^{\frac 1{p-1}}\|_{L^\infty(-1,1)}
\le C\|h\|_{\H_0}.
\]
(ii) {\bf (Boundedness of $\kappa(d,y)$ in several norms)}
For all $d\in (-1,1)$, it holds that
\[
\|\kappa(d,y)\|_{L^{p+1}_\rho}+ \|\kappa(d,y)(1-y^2)^{\frac 1{p-1}}\|_{L^\infty(-1,1)}
\le C\|\kappa(d,y)\|_{\H_0}\le CE(\kappa_0).
\]
(iii) {\bf (Same energy level for $\kappa(d,y)$)} For all $d\in (-1,1)$, it holds that
\[
E(\kappa(d,y))=E(-\kappa(d,y))=E(\kappa_0).
\]
(iv) {\bf (Continuity of the Lyapunov functional)} If $(w_i(y,s), \ps w_i(y,s))\in \H$ for $i=1$ and $2$ and for some $s\in \R$, then 
\begin{eqnarray*}
&&|E(w_1(s))-E(w_2(s))|\\
&\le& C\left\|\vc{w_1(s)}{\ps w_1(s)}-\vc{w_2(s)}{\ps w_2(s)}\right\|_{\H}\left(1+\left\|\vc{w_1(s)}{\ps w_1(s)}\right\|_{\H}^p+\left\|\vc{w_2(s)}{\ps w_2(s)}\right\|_{\H}^p\right).
\end{eqnarray*}
\end{cl}
{\it Proof}: For (i), see Lemma 2.2 page 51 in \cite{MZjfa07}. For (ii), use (i) and identity (49) page 59 in \cite{MZjfa07}. For (iii), see (ii) of Proposition 1 page 47 in \cite{MZjfa07}. The proof of (iv) is straightforward from the definition \aref{defenergy} of $E(w)$.\Box

\bigskip

To prove estimates about $\varphi$, we take advantage of the decoupling in the solitons' sum (see \aref{decouple}) and use information we proved in \cite{MZjfa07} for the 1-soliton version of $\varphi$ defined for all $d\in (-1,1)$, $r$ and $r'$ in $\H$ by
\begin{equation}\label{defphid}
\varphi_d(r,\r) = \iint \left(r_1' \r_1' (1-y^2)-\left(p\kappa(d)^{p-1} - \frac {2(p+1)}{(p-1)^2}\right) r_1\r_1+r_2\r_2\right)\rho dy
\end{equation}
and satisfying (see estimate (138) page 91 in \cite{MZjfa07}): 
\begin{equation}\label{contphid}
|\varphi_d(r,\r)|\le C\|r\|_{\H}\|\r\|_{\H}.
\end{equation}
It happens that the proof is clearer in the $\xi$ variable where 
\begin{equation*}
y=\tanh \xi.
\end{equation*}
 More precisely, let us introduce the transformations
\begin{equation}\label{deftr}
r(y) \mapsto \bt r(\xi) = \bar r(\xi)= r(y)(1-y^2)^{\frac 1{p-1}}\mbox{ and }
r(y) \mapsto \bc r(\xi) = \hat r(\xi)=r(y)(1-y^2)^{\frac 1{p-1}+\frac 12},
\end{equation}
and for $r=(r_1,r_2)$, the notation
\[
\tilde {\cal T}(r) = \tilde r = \vc{\bar r_1}{\hat r_2}=\vc{\bt(r_1)}{\bc(r_2)}.
\]
In the following claim, we transform $\varphi$ and $\varphi_d$ in the new set of variables. Let us first introduce the quadratic forms (where $d\in (-1,1)$):
\begin{eqnarray}
\bar \varphi_d(q,\q)&=&\int_{\R}\left(q_1'\q_1'+ \beta_d(\xi)q_1\q_1+q_2\q_2\right) d\xi\label{defbphid}\\
\bar \varphi(q,\q)&=&\int_{\R}\left(q_1'\q_1'+ \beta(\xi,s)q_1\q_1+q_2\q_2\right) d\xi\label{defbphi}
\end{eqnarray}
where (using \aref{defkd})
\begin{eqnarray}
\beta_d(\xi)&=& \frac 4{(p-1)^2}-p(\bar \kappa(d,y))^{p-1}
= \frac 4{(p-1)^2}- p \bar \kappa_0(\xi-\zeta)^{p-1}
\mbox{ with }
d=-\tanh \zeta, \nonumber\\
\bar \kappa_0(\xi)&=&\kappa_0 \cosh^{-\frac 2{p-1}}(\xi),\label{defbk0}\\
\beta(\xi,s) &=& \frac 4{(p-1)^2}- p|\bar K(\xi,s)|^{p-1}
= \frac 4{(p-1)^2}-p\left|\sum_{i=1}^k e_i\bar \kappa_0(\xi-\zeta_i(s))\right|^{p-1}.\label{defbk}
\end{eqnarray}
In the following claim, we give the effect of the new transformation: 
\begin{cl}\label{cltrans} $ $\\
(i) There exists $C_0>0$ such that for all $r\in \H$, we have
\[
\frac 1{C_0} \|r\|_{\H} \le \|\tilde r\|_{H^1\times L^2(\R)}\le C_0 \|r\|_{\H}.
\]
(ii) If $r_1\in \H_0$, then $(1-y^2)\bt\left(\L r_1 - \frac{2(p+1)}{(p-1)^2}r_1\right)= 
\left(\partial^2_{\xi}\bar r_1 - \frac 4{(p-1)^2}\bar r_1\right)$.\\
(iii) For all $r$, $\r$ in $\H$ and $d\in (-1,1)$, we have
\[
\varphi(r,\r)= \bar \varphi\left(\tilde r,\tilde \r\right)
\mbox{ and } 
\varphi_d(r,\r)= \bar \varphi_d\left(\tilde r,\tilde \r\right)
\]
where $\bar \varphi$ and $\bar \varphi_d$ are 
introduced in 
\aref{defbphi} and \aref{defbphid}.
\end{cl}
{\it Proof}:

(i) Consider $r=(r_1, r_2) \in \H$. Using \aref{deftr}, we first write
\begin{equation}\label{da}
\int_{\R} \bar r_1(\xi)^2 d\xi = \iint r_1(y)^2 \frac{\rho(y)}{1-y^2} dy
\mbox{ and }\int_{\R} \hat r_2(\xi)^2 d\xi = \iint r_2(y)^2 \rho(y) dy.\\
\end{equation}
Using Lemma \ref{lemsobolev-9}, we obtain
\begin{equation}\label{li}
\|r_1\|_{L^2_\rho}^2 \le \int_{\R} \bar r_1(\xi)^2 d\xi \le \|r\|_{\H}^2.
\end{equation}
Now, using again \aref{deftr}, we write
\begin{equation*}
\partial_\xi \bar r_1(\xi) = \py r_1(y)(1-y^2)^{\frac 1{p-1}+1} -\frac{2y}{p-1}(1-y^2)^{\frac 1{p-1}}r_1(y),
\end{equation*}
therefore, 
\begin{eqnarray*}
|\partial_\xi \bar r_1|^2 &\le & 2 |\py r_1|^2\rho (1-y^2)^2+ C|r_1|^2\rho,\\
|\py r_1|^2\rho (1-y^2)^2 &\le & 2 |\partial_\xi \bar r_1|^2+C|\bar r_1|^2.
\end{eqnarray*}
Integrating this and using Lemma \ref{lemsobolev-9}, we write
\begin{eqnarray}
&&\int_{\R} |\partial_\xi \bar r_1|^2  d\xi 
\le  2 \iint |\py r_1|^2\rho (1-y^2) dy + C \iint r_1^2 \frac \rho{1-y^2} dy \le C \|r\|_{\H}^2,\nonumber\\
&&\iint |\py r_1|^2\rho (1-y^2) dy \le  2 \int_{\R}|\partial_\xi \bar r_1|^2d\xi + C \int_{\R} |\bar r_1|^2d\xi.\label{daa}
\end{eqnarray} 
Gathering \aref{da}, \aref{li} and \aref{daa}, we conclude the proof of (i).

(ii) See page 60 in \cite{MZjfa07}.

(iii) We only prove the estimate for $\varphi$ since it is even easier for 
$\varphi_d$. Using the definitions \aref{defphi}, \aref{defro} and \aref{eqq} of $\varphi$, $\L$ and $\psi$,  integration by parts and the change of variables \aref{deftr}, we write
\begin{eqnarray*}
&&\varphi(r,\r)= \iint \left[-\L r_1.\r_1-\psi r_1\r_1+r_2\r_2\right]\rho dy\\
&=& \iint (-\L r_1+\frac{2(p+1)}{(p-1)^2}r_1)\r_1\rho dy
-p\iint r_1\r_1|K|^{p-1}\rho dy + \iint r_2\r_2\rho dy\\
&=& \int_{\R} (1-y^2)\bt(-\L r_1+\frac{2(p+1)}{(p-1)^2}r_1)\bar \r_1 d\xi
- p\int_{\R}\bar r_1 \bar \r_1|\bar K|^{p-1} d\xi
+\int_{\R}\hat r_2 \hat \r_2 d\xi 
 \end{eqnarray*}
Using (ii)
and integration by parts, we see that
\begin{eqnarray*}
\varphi(r,\r)&=& -\int_{\R}\left(\partial^2_{\xi}\bar r_1 - \frac 4{(p-1)^2}\bar r_1\right)\bar \r_1 d\xi
-p \int_{\R}\bar r_1 \bar \r_1 |\bar K|^{p-1} d\xi
+\int_{\R}\hat r_2 \hat \r_2 d\xi\\
&=& \bar \varphi(\tilde r, \tilde \r)
\end{eqnarray*}
where $\bar \varphi$ is introduced in \aref{defbphi}. This concludes the proof of Claim \ref{cltrans}.\Box

\bigskip

As we said earlier, we take advantage of the decoupling in the solitons' sum. In the following claim, we give a localized estimate coming from an identity we proved in \cite{MZjfa07} for $\bar\varphi_d$, the 1-soliton version of $\varphi$ defined in \aref{defphid}, then we derive a global estimate for $\bar \varphi$. 
\begin{cl}[Identity for $\varphi_d$]\label{clident}
$ $\\
(i) There exist $\epsilon_0>0$ and $A_0>0$ such that for all $A>A_0$, $d\in (-1,1)$ and $q\in H^1\times L^2(\R)$, we have
\[
\bar \varphi_d(q\sqrt {\chi_{A,d}},q\sqrt {\chi_{A,d}}) \ge \epsilon_0 \|q \sqrt {\chi_{A,d}}\|_{H^1\times L^2}^2-\frac{\epsilon_0}{8k}\|q\|_{H^1\times L^2}^2-
\sum_{\lambda=0}^1 \left|\pi_\lambda^d(\tilde T^{-1}(q))\right|^2
\]
where $\chi_{A,d}(\xi) = \chi_{1,0}(\frac{\xi-\zeta}A)$, $\tanh \zeta = -d$ and $\chi_{1,0}\in C^\infty(\R, [0,1])$ is even, decreasing for $\xi>0$ with $\chi_{1,0}(\xi)=1$ if $|\xi|<1$ and $\chi_{1,0}(\xi)=0$ if $|\xi|>2$.\\
(ii) There exists $\epsilon_2>0$ such that for $s$ large enough and for all $q\in H^1\times L^2$, we have
\begin{equation*}
\bar \varphi(q,q) \ge \epsilon_2 \|q\|_{H^1\times L^2}^2 - \frac 1{\epsilon_2}\sum_{i=1}^k  \sum_{\lambda=1}^2 
\left|\pi_\lambda^{d_i}(\tilde T^{-1}(q))\right|^2
\end{equation*}
\end{cl}
{\it Proof}:\\
(i) Consider some $d\in (-1,1)$ and $r\in \H$. 
On the one hand, we write from Proposition 4.7 page 90 in \cite{MZjfa07},
\begin{equation}\label{salama}
\varphi_d(r_{-,d}, r_{-,d}) \ge 2\epsilon_1 \|r\|_{\H}^2 - \frac 1{\epsilon_1}\sum_{\lambda=0}^1 |\pi_1^d(r)|^2
\end{equation}
for some $\epsilon_1>0$ where $r_-^d = \pi_-^d(r)$ defined in \aref{defpi-d}. On the other hand, using the continuity of $\varphi_d$ stated in \aref{contphid} and \aref{79bis}, we write
\begin{eqnarray*}
\varphi_d(r_{-,d}, r_{-,d})&\le & \varphi_d(r,r) +C\sum_{\lambda=0}^1|\pi_\lambda^d(r)|^2 + C \|r\|_{\H}\sum_{\lambda=0}^1|\pi_\lambda^d(r)|\\
&\le & \varphi_d(r,r) +\frac C{\epsilon_1}\sum_{\lambda=0}^1|\pi_\lambda^d(r)|^2 + \epsilon_1 \|r\|_{\H}.
\end{eqnarray*}
Using \aref{salama}, we see that
\[
\varphi_d(r,r) \ge \epsilon_0 \|r\|_{\H}^2 - \sum_{\lambda=0}^1 \left|\pi_\lambda^d(r)\right|^2.
\]
Using the $\xi$ framework and Claim \ref{cltrans}, we get for all $q \in H^1\times L^2(\R)$,
\begin{equation}\label{(i)}
\bar \varphi_d(q,q) \ge \epsilon_0 \|q\|_{H^1\times L^2(\R)}^2-
\sum_{\lambda=0}^1 \left|\pi_\lambda^d(\tilde T^{-1}(q))\right|^2.
\end{equation}
Now, 
%
we claim that (i) 
follows from the fact that for all $d\in (-1,1)$ and $\lambda=0$ or $1$, we have
\begin{equation}\label{monastere}
\forall u\in H^1\times L^2(\R),\;\;|\pi_\lambda^d(\tilde T^{-1}(u))|\le C \int \bar \kappa_0(\xi-\zeta)|(|u_1(\xi)|+|u_2(\xi)|) d\xi\mbox{ where } 
d= -\tanh \zeta.
\end{equation}
Indeed, consider $q\in H^1\times L^2$, $d\in (-1,1)$, $A>0$ and $\lambda =0$ or $1$. Taking 
\[
u= \tilde T^{-1}(q(1-\sqrt{\chi_{A,d}})),
\]
 using the Cauchy-Schwartz inequality and performing the change of variables $z= \xi-\zeta$, we see that
\begin{eqnarray*}
|\pi_\lambda^d(\tilde T^{-1}(q(1-\sqrt{\chi_{A,d}})))|&\le & C \int \bar \kappa_0(\xi-\zeta)(1-\sqrt{\chi_{A,d}})(|q_1(\xi)|+|q_2(\xi)|) d\xi\\
&\le & C \left(\int \bar \kappa_0(z)^2(1-\sqrt{\chi_{A,0}})^2 dz\right)^{1/2}
\|q\|_{H^1\times L^2}.
\end{eqnarray*}
Using Lebesgue's theorem, we find $A_0>0$ such that if $A\ge A_0$, then 
\[
|\pi_\lambda^d(\tilde T^{-1}(q(1-\sqrt{\chi_{A,d}})))|
\le \sqrt{\frac{\epsilon_0}{16 k}} \|q\|_{H^1\times L^2}
\]
(uniformly in $d\in (-1,1)$ of course). Since $\pi_\lambda^d$ is linear, this gives 
\[
|\pi_\lambda^d(\tilde T^{-1}(q\sqrt{\chi_{A,d}}))|^2 \le 
2 |\pi_\lambda^d(\tilde T^{-1}(q))|^2+\frac{\epsilon_0}{8k}\|q\|_{H^1\times L^2}^2.
\]
Using \aref{(i)} with $q\sqrt{\chi_{A,d}}$, (i) follows. It remains to prove \aref{monastere} to finish the proof of (i) of Claim \ref{clident}.


\bigskip

{\it Proof of \aref{monastere}}:  Consider $d\in (-1,1)$, $\lambda=0$ or $1$ and $u\in H^1\times L^2$. If we introduce 
$r= \tilde T^{-1}(u)$ 
which is in $\H$ by (i) of Claim \ref{cltrans}, then we have from \aref{defpdi} and integration by parts
\begin{equation}\label{lina}
\pi_\lambda^d(r) = \iint \left[(-\L W_{\lambda,1}(d) + W_{\lambda,1}(d))r_1 + W_{\lambda,2}(d)r_2\right] \rho(y) dy.
\end{equation}
Since we have from \aref{defWl2-0}, \aref{eqWl1-0} and \aref{defkd}
\[
W_{\lambda,2}(d,y) \le C\kappa(d,y)\mbox{ and }
|- \L W_{\lambda,1}(d,y)+W_{\lambda,1}(d,y)|\le C\frac{\kappa(d,y)}{1-y^2},
\]
we get from \aref{lina} and the transformation \aref{deftr}
\begin{eqnarray*}
|\pi_\lambda^d(r)|&\le & C \iint \kappa(d,y) |r_1(y)|\frac {\rho(y)}{1-y^2} dy 
+ C \iint \kappa(d,y) |r_2(y)| \rho(y) dy\\
&\le & C \int \bar \kappa(d,\xi)|u_1(\xi)| d\xi
+ \int \hat \kappa(d,\xi)|u_2(\xi)| d\xi.
\end{eqnarray*}
Since we have from \aref{deftr} and \aref{defbk0}, 
\[
\hat \kappa(d,\xi) \le \bar \kappa(d,\xi) = \bar \kappa_0 \cosh^{-\frac 2{p-1}}(\xi-\zeta)= \kappa_0(\xi-\zeta)\mbox{ with } d= -\tanh,
\]
\aref{monastere} follows. This concludes the proof of (i) in Claim \ref{clident}.

\bigskip

\noindent (ii) Introducing the notation
\[
\chi_i = \chi_{A,d_i}(\xi)=\chi_{1,0}\left(\frac{\xi-\zeta_i}A\right)
\]
%
%
and using \aref{defbphi}, we write
\begin{eqnarray}
\bar\varphi(q,q)&=&
\int(\partial_\xi q_1)^2+ \frac{2(p+1)}{(p-1)^2}\int q_1^2+\int q_2^2
- p\int |\bar K|^{p-1} q_1^2\nonumber\\
&=& \sum_{j=1}^k \left[\int(\partial_\xi q_1)^2 \chi_j
+\frac{2(p+1)}{(p-1)^2} \int q_1^2 \chi_j + \int q_2^2 \chi_j
- p\int |\bar K|^{p-1} q_1^2\chi_j\right]\nonumber\\
&+& \int(\partial_\xi q_1)^2(1-\sum_{j=1}^k \chi_j)
+\frac{2(p+1)}{(p-1)^2} \int q_1^2(1-\sum_{j=1}^k \chi_j)
+ \int q_2^2(1-\sum_{j=1}^k \chi_j)\nonumber\\
&-& p \int |\bar K|^{p-1}q_1^2 (1-\sum_{j=1}^k\chi_j)\nonumber\\
&=& \sum_{j=1}^k \bar\varphi(q \sqrt{\chi_j}, q \sqrt{\chi_j})
+\bar\varphi\left(q\sqrt{1- \sum_{j=1}^k \chi_j}, q\sqrt{1- \sum_{j=1}^k \chi_j}\right) + I_1(s)\label{hokm}
\end{eqnarray}
where 
\begin{eqnarray}
I_1(s)&=&-\sum_{j=1}^k\left\{ \int q_1^2 \left(\partial_\xi \sqrt{\chi_j}\right)^2 
- 2 \int q_1 \partial_\xi q_1 \sqrt{\chi_j} \partial_\xi \sqrt{\chi_j}\right\}\nonumber\\
&-& \int q_1^2 \left(\partial_\xi \sqrt{1-\sum_{j=1}^k}\chi_j\right)^2 
-2 \int q_1 \partial_\xi q_1 \sqrt{1-\sum_{j=1}^k \chi_j}   \partial_\xi \sqrt{1-\sum_{j=1}^k \chi_j}.\label{122bis}
\end{eqnarray}
Using the definitions \aref{defbphi} and \aref{defbphid} 
of $\bar \varphi$ and $\bar \varphi_d$,
we write
\begin{eqnarray*}
&&\bar\varphi\left(q \sqrt{\chi_j}, q \sqrt{\chi_j}\right)
=\bar\varphi_{d_i(s)}\left(q\sqrt{\chi_j}, q\sqrt{\chi_j}\right)- I_2(s),\label{dal}\\
&&\bar\varphi\left(q\sqrt{1-\d\sum_{j=1}^k\chi_j}, q\sqrt{1-\d\sum_{j=1}^k \chi_j}\right)
\ge c_0(p)\left\|q\sqrt{1-\d\sum_{j=1}^k\chi_j}\right\|_{H^1\times L^2}^2-I_3(s)\label{ida}
\end{eqnarray*}
where $c_0(p) =\min \left(1, \frac{2(p+1)}{(p-1)^2}\right)$,
\[
I_2(s)= p\d\int \left(|\bar K(\xi,s)|^{p-1} - \bar \kappa_0(\xi-\zeta_i(s))^{p-1}\right) q_1^2\chi_j
\mbox{ and }
I_3(s) = p\d\int |\bar K|^{p-1}q_1^2 (1-\sum_{j=1}^k \chi_j).
\]
Since we have from Claim \ref{clident} and \aref{defbk}
\begin{eqnarray}
|\partial_\xi \chi_j|\le C/A,&&
\|\left(|\bar K(\xi,s)|^{p-1} - \bar \kappa_0(\xi-\zeta_i(s))^{p-1}\right)\chi_j\|_{L^\infty} \le C(A) J(s),\label{besame}\\
\mbox{and}&&
\||\bar K|^{p-1}(1-\sum_{j=1}^k \chi_j)\|_{L^\infty}
\le C e^{-2A}\nonumber
\end{eqnarray}
where $J(s)\to 0$ is defined in \aref{defJ}, it follows that for $A$ and $s$ large enough, 
\begin{equation}\label{boundIi}
|I_1(s)|+|I_2(s)|+|I_3(s)| \le \frac CA \|q_1\|_{H^1}^2.
\end{equation}
Therefore, using \aref{hokm}, \aref{122bis}, \aref{dal}, \aref{ida}, \aref{boundIi} and Claim \ref{clident}, we write for $A$ and $s$ large enough,
\begin{eqnarray}
\bar\varphi(q,q) &\ge& \epsilon_0 \sum_{j=1}^k \|q\sqrt{\chi_j}\|_{H^1\times L^2}^2 +c_0(p)\left\|q\sqrt{1-\sum_{j=1}^k \chi_j}\right\|_{H^1\times L^2}^2-\frac {\epsilon_0}4 \|q\|_{H^1\times L^2}^2\label{mucho}\\
& -&  \sum_{j=1}^k  \sum_{\lambda=0}^1
|\pi_\lambda^{d_j(s)}(\tilde T^{-1}(q))|^2.\nonumber
\end{eqnarray}
Since \aref{hokm} holds with $\bar \varphi$ replaced by the canonical inner product of $H^1\times L^2$, we use \aref{boundIi} to write
\[
\|q\|_{H^1\times L^2}^2 \le \sum_{j=1}^k \|q\sqrt{\chi_j}\|_{H^1\times L^2}^2 +\left\|q\sqrt{1-\sum_{j=1}^k \chi_j}\right\|_{H^1\times L^2}^2+\frac CA \|q\|_{H^1\times L^2}^2\]
hence for $A$ and $s$ large enough,
\[
\|q\|_{H^1\times L^2}^2 \le 2\sum_{j=1}^k \|q\sqrt{\chi_j}\|_{H^1\times L^2}^2 +2\left\|q\sqrt{1-\sum_{j=1}^k \chi_j}\right\|_{H^1\times L^2}^2
\] 
and 
(ii) follows from \aref{mucho}. This concludes the proof of Claim \ref{clident}.\Box

\subsection{Proof of Lemma \ref{lemdecomp}}
Now we are ready to start the proof of Lemma \ref{lemdecomp}.

\medskip

{\it Proof of Lemma \ref{lemdecomp}}:

\medskip

(i) Since $\psi(y,s)=p|K(y,s)|^{p-1} -\frac{2(p+1)}{(p-1)^2}$ with $K(y,s) = \sum_{j=1}^k e_j \kappa(d_j(s),y)$ by \aref{eqq}, we split $\varphi(r,r')$ into 2 parts as follows:\\
- We first use the definition \aref{defnh0} of the norm in $\H$ to write
\[
\left|\iint \left(\py r_1 \py r_1' (1-y^2)+\frac{2(p+1)}{(p-1)^2}r_1r_1'+r_2r_2'\right)\rho dy\right| \le
C \|r\|_{\H} \|r'\|_{\H}.
\]
- Then, using Claim \ref{lemsobolev-9}, we write
\[
\left|\iint |K(s)|^{p-1}r_1 r_1' \rho dy\right|\le  
C\iint \frac{|r_1| |r_1'|}{1-y^2} \rho dy\le 
C\|r_1\|_{L^2_{\frac \rho{1-y^2}}}\|r_1'\|_{L^2_{\frac \rho{1-y^2}}}
\le C\|r\|_{\H}\|r'\|_{\H}.
\]
Using these two bounds gives the conclusion of (i).

\bigskip

(ii) {\it Proof of \aref{qdecomp}}: It immediately follows from \aref{decomp}, \aref{defa-} and \aref{peur}.

\medskip

\noindent {\it Proof of \aref{48bis}}: The right inequality follows from (i). For the left inequality, we use Claim \ref{cltrans} to translate (ii) of Claim \ref{clident} back to the $y$ variable:\\
for some $\epsilon_2>0$, for $s$ large enough and for all $r\in \H$, 
\begin{equation}\label{rima}
\varphi(r,r) \ge \epsilon_2 \|r\|_{\H}^2 - \frac 1{\epsilon_2} \sum_{i=1}^k \sum_{\lambda=0}^1 |\pi_\lambda^{d_i(s)}(r)|^2.
\end{equation}
Using \aref{rima} with $r(y)= q_-(y,s)$, we write
\begin{equation}\label{imen0}
\varphi(q_-,q_-) \ge \epsilon_2 \|q_-\|_{\H}^2 - \frac 1{\epsilon_2} \sum_{i=1}^k \sum_{\lambda=0}^1 |\pi_\lambda^{d_i}(q_-)|^2.
\end{equation}
Since $\pi_\lambda^{d_i}(F_1^{d_i})= \delta_{\lambda,1}$ by \aref{orth}, we use \aref{qdecomp} to write
\begin{equation}\label{imen1}
\pi_\lambda^{d_i}(q_-) = \pi_\lambda^{d_i}(q) - \sum_{j=1}^k \pi_1^{d_j}(q) \pi_\lambda^{d_i}(F_1^{d_j})
= \sum_{j\neq i} \pi_1^{d_j}(q) \pi_\lambda^{d_i}(F_1^{d_j}).
\end{equation}
Using \aref{defa-}, \aref{defpdi} and \aref{normw}, we see that
\begin{equation}\label{imen2}
|\alpha_1^j|= |\pi_1^{d_j}(q)|
= |\phi(W_\lambda(d_j),q)|\le\|W_\lambda(d_j)\|_{\H}\|q\|_{\H} \le C\|q\|_{\H}.
\end{equation}
Using \aref{defpdi}, integration by parts and the definition \aref{defro} of $\L$, we write
\begin{eqnarray*}
\pi^{d_i}_\lambda(F_\mu(d_j))
&=& \iint \left(W_{\lambda,1}(d_i)F_{\mu,1}(d_j)+
\py W_{\lambda,1}(d_i)\py F_{\mu,1}(d_j)(1-y^2)\right)\rho dy\\
&+& \iint W_{\lambda,2}(d_i)F_{\mu,2}(d_j) \rho dy\\
\\
&=& \iint \left(-\L W_{\lambda,1}(d_i) + W_{\lambda,1}(d_i)\right)F_{\mu,1}(d_j)\rho dy+ \iint W_{\lambda,2}(d_i)F_{\mu,2}(d_j) \rho dy. 
\end{eqnarray*}
Since we have from the definitions  \aref{defkd}, \aref{defWl2-0}, \aref{eqWl1-0} and \aref{deffld} of $\kappa(d,y)$, $W_{\lambda}(d,y)$ and $F_{\mu}(d,y)$, for all $(d,y) \in (-1,1)^2$, 
\begin{equation}\label{neila}
|W_{\lambda,2}(d,y)|+|\L W_{\lambda,1}(d,y) - W_{\lambda,1}(d,y)|\le C \frac {\kappa(d,y)}{1-y^2}
\mbox{ and }|F_{\mu,l}(d,y)|\le C \kappa(d,y), 
\end{equation}
we use (i) of Lemma \ref{lemtech} to write for $s$ large enough,
\begin{equation}\label{mas}
\left|\pi^{d_i}_\lambda(F_\mu(d_i))\right|\le C \iint \kappa(d_i)\kappa(d_j)\frac \rho{1-y^2} dy \le C |\zeta_i-\zeta_j|e^{-\frac 2{p-1}|\zeta_i-\zeta_j|}
\le C\bar J(s)
\end{equation}
by definition \aref{defa-} of $\bar J$. Using \aref{imen1} and \aref{imen2}, we see that for $s$ large enough,
\[
|\pi_\lambda^{d_i}(q_-)| \le C \bar J \|q\|_{\H}.
\]
Using \aref{imen0}, we see that the left inequality in \aref{48bis} follows.

\medskip

\noindent {\it Proof of \aref{equivq}}: The right inequality follows from \aref{imen2}, \aref{contphi} and \aref{48bis}. 
%
%
%
%
For the left inequality in \aref{equivq}, we write from the bilinearity of $\varphi$, \aref{qdecomp}, \aref{contphi} and \aref{79bis}
\begin{eqnarray*}
\varphi(q_-, q_-)&\ge& \varphi(q, q) - C\sum_{i=1}^k |\alpha_1^i|^2
-C\|q\|_{\H}\sum_{i=1}^k |\alpha_1^i|\\
&\ge & \varphi(q, q) - \frac C{\epsilon_2}\sum_{i=1}^k |\alpha_1^i|^2
-\frac{\epsilon_2}2\|q\|_{\H}^2
\end{eqnarray*}
where $\epsilon_1>0$ is introduced in \aref{rima}. Using \aref{rima} with $r=q$, we get the left inequality in \aref{equivq}. This concludes the proof of Lemma \ref{lemdecomp}.\Box.
\section{Projection of equation \aref{eqq} on the different modes}\label{appproj}
We prove Lemma \ref{lemproj} here.
We proceed in 3 parts to prove (i), (ii), and finally (iii).

\bigskip

{\bf Proof of (i): Projection of equation \aref{eqq} on $F_1(d_i(s),\cdot)$ and $F_0(d_i(s),\cdot)$}   

We prove (i) of Lemma \ref{lemproj} here. Fixing some $i=1,...,k$ and projecting equation \aref{eqq} with the projector $\pi^d_\lambda$ \aref{defpdi} (where $\lambda=0$ or $1$), we write (putting on top the main terms)
\begin{eqnarray}
&&\pi^{d_i(s)}_\lambda\left(\partial_s q\right)
= \pi^{d_i(s)}_\lambda\left(L_{d_i(s)}(q)\right)
- e_id_i'(s) \pi^{d_i(s)}_\lambda\vc{\partial_d \kappa(d_i(s),y)}{0}
+\pi^{d_i(s)}_\lambda\vc{0}{R}\nonumber\\ 
&+& \pi^{d_i(s)}_\lambda\vc{0}{f(q_1)}
+ \pi^{d_i(s)}_\lambda\vc{0}{V_i(y,s) q_1}
-\sum_{j\neq i}^ke_jd_j'(s) \pi^{d_i(s)}_\lambda\vc{\partial_d \kappa(d_j(s),y)}{0}. \label{adel1-9}
\end{eqnarray}
Note that we expand the operator $L(q)$ according to \aref{defvi-0}. In the following, we handle each term of \aref{adel1-9} in order to finish the proof of \aref{eqa-star}.

- Using the analysis performed in Claim 5.3 page 104 and Step 1 page 105 in \cite{MZjfa07} for the case of one soliton ($k=1$), we immediately get the following estimates:
\begin{eqnarray}
\left|\pi^{d_i(s)}_\lambda(\partial_s q)-{\alpha^i_\lambda}'(s)\right|
&\le& \frac{C_0}{1-(d_i(s))^2}|d_i'(s)|\|q(s)\|_{\H}\le C_0 |\zeta_i'(s)|\|q(s)\|_{\H},\nonumber\\
\pi^{d_i(s)}_\lambda\left(L_{d_i(s)}(q)\right)&=&\lambda \alpha^i_\lambda(s),\nonumber\\
d_i'(s)\pi^{d_i(s)}_\lambda \vc{\partial_d \kappa(d_i(s),y)}{0}&=&- \frac{2\kappa_0}{(p-1)}\frac{d_i'(s)}{(1-(d_i(s))^2)}\delta_{\lambda,0}= \frac{2\kappa_0}{(p-1)} \zeta_i'(s)\delta_{\lambda,0},\nonumber\\
|f(q_1)|&\le & C\delta_{\{p\ge 2\}}|q_1|^p+C|K|^{p-2}|q_1|^2\label{nonl-9}
\end{eqnarray}
(recall that 
\begin{equation}\label{recall}
d_i(s)=-\tanh \zeta_i(s),\mbox{ hence }\zeta_i'(s)= - \frac{d_i'(s)}{1-d_i(s)^2}\mbox{)}.
\end{equation}

- Since we have from the definitions \aref{eqq}, \aref{defkd} and \aref{defWl2-0} of $R$, $\kappa(d,y)$ and $W_{\lambda,2}(d,y)$ 
\begin{equation}\label{terreur}
|R(y,s)| \le C \sum_{j\neq i} \kappa(d_j,y)^p + \kappa(d_i,y)^{p-1} \kappa(d_j,y)\mbox{ and }|W_{\lambda,2}(d_i,y)|\le C\kappa(d_i,y),
\end{equation}
we use \aref{defpdi}, (i) of Lemma \ref{lemtech} and the definition \aref{defJ} of $J(s)$ to write for $s$ large enough, 
\begin{eqnarray}
\left|\pi^{d_i}_\lambda\vc{0}{R}\right| = \left|\iint W_{\lambda,2}(d_i)R\rho dy\right| &\le& C \sum_{j\neq i} \iint \kappa(d_i) \kappa(d_j)^p \rho dy + \iint \kappa(d_i)^p \kappa(d_j) \rho dy\nonumber\\
&\le & \sum_{j\neq i} e^{-\frac 2{p-1}|\zeta_i-\zeta_j|} \le CJ. \label{coc-9}
\end{eqnarray}

- Using \aref{defpdi}, \aref{terreur}, \aref{nonl-9} and the H\"older inequality, we write
\begin{eqnarray}
\left|\pi_\lambda^{d_i}\vc{0}{f(q_1)}\right|
&\le& C \iint \kappa(d_i)|f(q_1)|\rho dy\nonumber\\
&\le & C \delta_{p\ge 2}\iint \kappa(d_i)|q_1|^p \rho dy
+C\iint \kappa(d_i)|K|^{p-2}|q_1|^2 \rho dy\nonumber\\
&\le &C\delta_{p\ge 2}\|\kappa(d_i)\|_{L^{p+1}_\rho}\|q_1\|_{L^{p+1}_\rho}^p
+CJ_i\|q_1(1-y^2)^{\frac 1{p-1}}\|_{L^\infty}^2\nonumber
\end{eqnarray}
where 
\begin{equation}
J_i =  \iint \kappa(d_i)|K|^{p-2} dy.
\end{equation}
Using (v) of Lemma \ref{lemtech} and Claim \ref{lemsobolev-9}, we see that
\begin{equation}\label{69bis-9}
\left|\pi_\lambda^{d_i}\vc{0}{f(q_1)}\right|
\le C \iint \kappa(d_i)|f(q_1)|\rho dy
\le C \delta_{p\ge 2}  \|q\|_{\H}^p
+C\|q\|_{\H}^2\le C \|q\|_{\H}^2
\end{equation}
where we use \aref{kill} in the last step.

\bigskip

- We claim that for some $\delta_1>0$ and for $s$ large enough, 
\begin{equation}\label{70-9}
\left|\pi^{d_i}_\lambda\vc{0}{V_i q_1}\right|
\le C\|q\|_{\H}^2 + CJ^{1+\delta_1}
\end{equation}
 and
\begin{equation}\label{krim1}
|V_i(y,s)|\le C\un{y_{i-1}}{y_i} \sum_{l\neq i}\kappa(d_i,y)^{p-2}\kappa(d_l,y) + C \sum_{l\neq i} \kappa(d_l,y)^{p-1}\un{y_{l-1}}{y_l}
\end{equation}
where  $y_0=-1$, $y_j= \tanh(\frac{\zeta_j+\zeta_{j+1}}2)$ if $j=1,..,k-1$ and $y_k=1$. In particular, we have,
\[
-1=y_0< -d_1<y_1<-d_2<...<y_j<-d_j<y_{j+1}<...<-d_k<y_k=1
\]
and $\kappa(d_j(s),y_{j+1}(s))=\kappa(d_{j+1}(s),y_{j+1}(s))$ for $j=1,...,k-1$ (use \aref{tkappa} to see this).\\
We first prove \aref{krim1} and then \aref{70-9}.
To prove \aref{krim1},  using \aref{decouple}, we see that:
\begin{itemize}
\item if $y\in (y_{i-1}(s), y_i(s))$, then $|\sum_{l\neq i}e_l \kappa(d_l(s),y)|\le 3 \kappa(d_i(s),y)$, hence from  \aref{defvi-0},
 $|V_i(y,s)|\le C \sum_{l\neq i} \kappa(d_i(s),y)^{p-2}\kappa(d_l(s),y)$;
\item if $y\in  (y_{l-1}(s), y_l(s))$ for some $l\neq i$, then for all $j=1,...,k$, $\kappa(d_j(s),y) \le \kappa(d_l(s),y)$, hence, from \aref{defvi-0}, $|V_i(y,s)|\le C \sum_{j=1}^k \kappa(d_j(s),y)^{p-1}\le C \kappa(d_l(s),y)^{p-1}$.
\end{itemize}
Thus, \aref{krim1} follows. Now, we prove \aref{70-9}.
 Using \aref{defpdi}, \aref{terreur}, Claim \ref{lemsobolev-9} and \aref{krim1}, we write
\begin{eqnarray}
&&\left|\pi^{d_i}_\lambda\vc{0}{V_i q_1}\right|
\le C\iint \kappa(d_i)\left|V_i q_1\right|\rho dy\nonumber\\
&\le& C \|q_1(1-y^2)^{\frac 1{p-1}}\|_{L^\infty}^2
+C\left(\iint \kappa(d_i)\left|V_i\right|(1-y^2)^{\frac 1{p-1}} dy\right)^2\nonumber\\
&\le & C\|q\|_{\H}^2\nonumber\\
&+&C\sum_{l\neq i}\left(\int_{y_{i-1}}^{y_i} \kappa(d_i)^{p-1}\kappa(d_l)(1-y^2)^{\frac 1{p-1}} dy\right)^2
+ \left(\int_{y_{l-1}}^{y_l} \kappa(d_l)^{p-1}\kappa(d_i)(1-y^2)^{\frac 1{p-1}} dy\right)^2.\nonumber
\end{eqnarray}
Using (ii) of Lemma \ref{lemtech},  \aref{70-9} follows.

\bigskip

- Consider $j\neq i$. Since we have from the definitions \aref{defkd} and \aref{deffld} of $\kappa(d,y)$ and $F_0(d,y)$,
\begin{equation}\label{ddkd}
\vc{\partial_d \kappa(d,y)}{0}= -\frac{2\kappa_0}{(p-1)(1-d^2)}F_0(d,y),
\end{equation}
we use \aref{mas} and \aref{recall}
to write
\begin{equation}\label{74-9}
\left|d_j'\pi^{d_i}_\lambda\vc{\partial_d \kappa(d_j)}{0}\right|
\le \frac {C|d_j'|}{1-d_j^2}
\left|\pi^{d_i}_\lambda(F_0(d_j))\right| \le C\bar J|\zeta_j'|
\end{equation}
where $\bar J(s)$ is defined in \aref{defa-}.
Using \aref{adel1-9}, \aref{nonl-9}, \aref{coc-9}, \aref{69bis-9}, \aref{70-9}, \aref{74-9} and \aref{peur}, we write for all $i=1,..,k$ (starting with $\lambda=0$ and then $\lambda=1$), 
\begin{eqnarray}
|\zeta_i'| &\le& C |\zeta_i'| \|q\|_{\H} + CJ + C\|q\|_{\H}^2
+C\bar J\sum_{j\neq i} |\zeta_j'|,\label{fayrouz-9}\\
\left|\ap'-\ap\right|&\le & C |\zeta_i'| \|q\|_{\H} + CJ + C\|q\|_{\H}^2
+C\bar J\sum_{j\neq i} |\zeta_j'|,\label{fayrouz-9bis}
\end{eqnarray}
Since $\|q\|_{\H}+\bar J\to 0$ (see (i) of Lemma \ref{clmod} and \aref{defa-}), summing up \aref{fayrouz-9} in $i$, we get,
\[
\sum_{i=1}^k|\zeta_i'|\le CJ+C\|q\|_{\H}^2.
\]
Plugging this in \aref{fayrouz-9bis}, we get 
\[
\left|\ap'-\ap\right|\le  CJ+C\|q\|_{\H}^2,
\]
which closes the proof of \aref{eqa-star}. This concludes the proof of (i) of Lemma \ref{lemproj}. \Box

\bigskip

{\bf Proof of (ii): Differential inequality satisfied by $A_-(s)$}

We proceed in 2 steps: we first project equation \aref{eqq} with the projector $\pi_-$ defined in \aref{decomp}, and then use that equation to write a differential inequality for $A_- = \varphi(q_-, q_-)$.

\bigskip

{\bf Step 2.1 : Projection of equation \aref{eqq} with $\pi_-$}

In this claim, we project equation \aref{eqq} with the projector $\pi_-$ defined in \aref{decomp}:
\begin{cl}[A partial differential inequality for $q_-$] \label{clrachid} For $s$ large enough, we have
\begin{eqnarray*}
&&\left\|\ps q_--L q_- -\sum_{i=1}^k\pi^{d_i}_1(q)\vc{0}{V_i F_{1,1}(d_i)}-\vc{0}{f(q_1)}-\vc{0}{R}\right\|_{\H}\\
& \le& C J + C\|q\|_{\H}^2
\end{eqnarray*}
where $J(s)$ is defined in \aref{defJ}.
\end{cl}
{\it Proof}: 
Applying the projector $\pi_-$ defined in \aref{decomp} to equation \aref{eqq},  we write
\begin{equation}\label{rachid}
\pi_-\left(\partial_s q\right)
= \pi_-\left(L q\right)
-\sum_{i=1}^ke_id_i' \pi_-\vc{\partial_d \kappa(d_i)}{0}
+ \pi_-\vc{0}{f(q_1)}
+\pi_-\vc{0}{R}.
\end{equation}
In the following, we will estimate each term appearing in this identity.

- Proceeding as for estimate (213) in \cite{MZjfa07} in the case of one soliton, one can straightforwardly control the left-hand term as follows:
\begin{equation}\label{ra1}
\|\pi_-(\ps q) - \ps q_-\|_{\H}\le  CJ\|q\|_{\H}+C\|q\|_{\H}^3.
\end{equation}

- We claim that for $s$ large enough, 
\begin{equation}\label{ha}
\left\|\pi_-(L q) - L q_- -\sum_{i=1}^k\pi^{d_i}_1(q)\vc{0}{V_i F_{1,1}(d_i)}\right\|_{\H}
\le C \|q\|_{\H}^2 + CJ^{1+\delta_1}
\end{equation}
where $\delta_1>0$ is introduced in \aref{70-9}.
Indeed, applying the operator $L$ to \aref{qdecomp} on the one hand, and using \aref{decomp} with $r=L q$ on the other hand, we write
\begin{eqnarray}
L q&=& \sum_{i=1}^k \pi_1^{d_i(s)}(q) L F_1(d_i(s), \cdot) + L q_-\nonumber\\
&=&\sum_{i=1}^k \pi_1^{d_i(s)}(L q) F_1(d_i(s),\cdot)+ 
\sum_{i=1}^k \pi_0^{d_i(s)}(L q) F_0(d_i(s), \cdot)
+\pi_-(Lq).\nonumber
\end{eqnarray}
Therefore, 
\begin{equation}\label{werwar}
\pi_-(L q) - L q_- = \sum_{i=1}^k \pi^{d_i(s)}_1 (q) L F_1 (d_i(s), \cdot)-\pi^{d_i(s)}_1 (L q) F_1(d_i(s),\cdot) -  \sum_{i=1}^k \pi^{d_i(s)}_0 (L q) F_0(d_i(s), \cdot).
\end{equation}
Since we have from \aref{defpdi} and \aref{79bis*}, $\pi^d_\lambda (L_d r) = \phi(W_\lambda(d,\cdot), L_d r) = \phi(L_d^* W_\lambda(d, \cdot), r) = \lambda \pi_\lambda^d(r)$, using this with \aref{defvi-0} and \aref{79bis} gives for $\lambda=0$ or $1$,
\begin{eqnarray}
L F_\lambda(d_i(s),\cdot) 
&=& L_{d_i(s)}F_\lambda(d_i(s),\cdot)  + \vc{0}{V_iF_{\lambda,1}(d_i(s),\cdot)}\nonumber\\
&=& \lambda F_\lambda(d_i(s),\cdot) + \vc{0}{V_iF_{\lambda,1}(d_i(s),\cdot)}\label{w2}\\
\pi^{d_i(s)}_\lambda(L q) 
&=& \pi^{d_i(s)}_\lambda(L_{d_i(s)} q)+ \pi_\lambda^{d_i(s)}\vc{0}{V_i q_1} 
= \lambda \pi^{d_i(s)}_\lambda(q) + \pi_\lambda^{d_i(s)}\vc{0}{V_i q_1}.\label{w1}
\end{eqnarray}
Using \aref{werwar}, \aref{w2} and \aref{w1} together with \aref{70-9} and \aref{79bis}, we get \aref{ha}.

\medskip

- Using the definition \aref{decomp} of the operator $\pi_-$, we see that
\begin{eqnarray}
\pi_-\vc{\partial_d \kappa(d_i, \cdot)}{0}&=& \vc{\partial_d \kappa(d_i, \cdot)}{0}\label{ghazi}\\
& -& \sum_{j=1}^k \pi_1^{d_j}\vc{\partial_d \kappa(d_i, \cdot)}{0}F_1(d_j,\cdot)-  \sum_{j=1}^k \pi_0^{d_j}\vc{\partial_d \kappa(d_i, \cdot)}{0}F_0(d_j,\cdot).\nonumber
\end{eqnarray}
Using \aref{ddkd}, 
it follows from the orthogonality relation \aref{orth} that for $\lambda=0$ or $1$,
\[
\pi_\lambda^{d_i}\vc{\partial_d \kappa(d_i, \cdot)}{0}
= -\frac{2\kappa_0}{(p-1)(1-d^2)}\pi_\lambda^{d_i}(F_0(d_i, \cdot))
= \delta_{\lambda,0}\vc{\partial_d \kappa(d_i, \cdot)}{0}.
\]
Therefore, it follows from \aref{ghazi} that
\[
\pi_-\vc{\partial_d \kappa(d_i, \cdot)}{0}= - \sum_{j\neq i} \pi_1^{d_j}\vc{\partial_d \kappa(d_i, \cdot)}{0}F_1(d_j,\cdot)-  \sum_{j\neq i} \pi_0^{d_j}\vc{\partial_d \kappa(d_i, \cdot)}{0}F_0(d_j,\cdot).
\]
Using \aref{74-9}, \aref{79bis} and \aref{eqa-star}, we see that
\begin{equation}\label{ra2}
\left\|d_i'(s)\pi_-\vc{\partial_d \kappa(d_i(s),\cdot)}{0}\right\|_{\H}
\le C\bar J(s)|\zeta_i'(s)|\le C\bar J(s)\left(\|q(s)\|_{\H}^2 + J(s)\right).
\end{equation}
- From definition \aref{decomp} of the operator $\pi_-$, \aref{79bis}, \aref{69bis-9} and \aref{coc-9}, we have
\begin{eqnarray}
\left\|\pi_-\vc{0}{f(q_1)}- \vc{0}{f(q_1)}\right\|_{\H} &\le& C\sum_{\lambda=1,2;\;i=1}^ k \left|\pi^{d_i(s)}_\lambda\vc{0}{f(q_1)}\right|\le C\|q(s)\|_{\H}^2,\label{ra3}\\
\left\|\pi_-\vc{0}{R}- \vc{0}{R}\right\|_{\H} &\le& C\sum_{\lambda=1,2;\;i=1}^ k \left|\pi^{d_i(s)}_\lambda\vc{0}{R}\right|\le CJ(s).\label{ra4}
\end{eqnarray} 
Using \aref{rachid}, \aref{ra1}, \aref{ha}, \aref{ra2}, \aref{ra3} and \aref{ra4} closes the proof of Claim \ref{clrachid}.\Box

\bigskip

{\bf Step 2.2: A differential inequality on $A_-(s)$} 

By definition \aref{defa-} of $\alpha_-(s)$, it holds that 
\begin{equation}\label{jazz}
\frac 12A_-'(s) = \varphi(\ps q_-, q_-)-\frac {p(p-1)}2 \sum_{i=1}^k e_i d_i' I_i
\end{equation}
 with 
\[
I_i=\iint\pd \kappa(d_i)|K|^{p-2}\left(q_{-,1}\right)^2 \rho dy \mbox{ and }K = \sum_{j=1}^k e_j \kappa(d_j).
\]
Since we have from \aref{ddkd} and \aref{neila}, $|\partial_d \kappa_{d_i}|\le C\frac{\kappa(d_i)}{1-d_i^2}$, 
using \aref{recall}, (v) of Lemma \ref{lemtech}, (i) of Claim \ref{lemsobolev-9} and \aref{91,5},
we see that
\begin{equation}\label{blues}
|d_i'||I_i|\le  \frac {C|d_i'|}{1-d_i^2} \iint \kappa(d_i)|K|^{p-2} dy \|q_{-,1}(1-y^2)^{\frac 1{p-1}}\|_{L^\infty}^2\le   C|\zeta_i'(s)|\|q_-\|_{\H}^2\le   C|\zeta_i'(s)|\|q\|_{\H}^2.
\end{equation}
Using \aref{jazz}, \aref{blues} and 
\aref{eqa-star},
we get
\begin{equation}\label{zz1}
\left|\frac 12A_-'(s)- \varphi(\ps q_-, q_-)\right|
\le C\|q\|_{\H}^2\left(\|q\|_{\H}^2 + J\right).
\end{equation}
Using \aref{91,5},
\aref{contphi} and Claim \ref{clrachid}, we estimate $\varphi(\ps q_-, q_-)$ in the following:
\begin{eqnarray}
&&\left|\varphi(\ps q_-, q_-) - \varphi(L q_-, q_-) -
\iint q_{-,2} f(q_1) \rho dy
-\iint q_{-,2} G \rho(y) dy\right|\nonumber\\
&\le& C\|q_-\|_{\H}\left(J + \|q\|_{\H}^2\right)
\le CJ\sqrt{|A_-|} +C J\bar J\|q\|_{\H} + C\|q\|_{\H}^3\nonumber\\
&\le & CJ\sqrt{|A_-|} + C\|q\|_{\H}^3+C \sum_{m=1}^{k-1}(h(\zeta_{m+1}-\zeta_m))^2\label{zz2}
\end{eqnarray}
where $h$ is defined in \aref{defh} and 
\begin{equation}\label{defG}
G(y,s) = \sum_{i=1}^k \alpha_1^i(s)V_i(y,s) F_{1,1}(d_i(s),y) +R(y,s).
\end{equation}
In the following, we estimate every term of \aref{zz2} in order to finish the proof of \aref{eqam}.

\noindent - Arguing as in page 107 of \cite{MZjfa07}, we write
\begin{equation}\label{zz3}
\varphi(L q_-, q_-) = -\frac 4{p-1} \iint q_{-,2}^2 \frac \rho{1-y^2} dy.
\end{equation}
- Since we have from the definitions \aref{defkd} and \aref{deffld} of $\kappa(d,y)$ and $F_1(d,y)$, 
\begin{equation}\label{boundF12}
F_{1,1}(d,y)= F_{1,2}(d,y)\le C \kappa(d,y),
\end{equation}
using \aref{qdecomp} and \aref{69bis-9}, we write
\begin{equation}
\left|\iint q_{-,2} f(q_1) \rho dy- \iint q_2 f(q_1) \rho dy\right|
\le C\sum_{i=1}^k |\alpha_1^i| \iint \kappa(d_i) |f(q_1)| \rho dy
\le C\|q\|_{\H}^3.\label{moncef}
\end{equation}
If we introduce 
\[
\F(q_1) = \int_0^{q_1} f(\xi) d\xi= \frac{|K+q_1|^{p+1}}{p+1} - \frac{|K|^{p+1}}{p+1}- |K|^{p-1} K q_1- \frac p2 |K|^{p-1} q_1^2,
\]
then it is easy to see that
\begin{equation}\label{boundF}
|\F (q_1)|\le C |q_1|^{p+1} + C \delta_{\{p\ge 2\}} |K|^{p-2}|q_1|^3.
\end{equation}
Introducing $R_- = -\iint \F(q_1) \rho dy$ and using equation \aref{eqq}, we write
\begin{eqnarray}
&&R_-'+ \iint q_2 f(q_1) \rho dy
= R_-'+ \iint \ps q_1 f(q_1) \rho dy + \sum_{i=1}^k d_i' \iint \pd \kappa(d_i) f(q_1)\rho dy\label{gafsia}\\
&=& \sum_{i=1}^k d_i' \iint \left(\pd \kappa(d_i) f(q_1)-\partial_{d_i} \F(q_1)\right) \rho dy
= \frac{p(p-1)}2 
\sum_{i=1}^k d_i' \iint \pd \kappa(d_i) |K|^{p-2}q_1^2\rho dy.\nonumber
\end{eqnarray}

Therefore, using \aref{moncef} and \aref{gafsia},  arguing as for \aref{blues}, using 
(v) of Lemma \ref{lemtech}
and  \aref{eqa-star}, we write
\begin{eqnarray}
\left|\iint q_{2,-} f(q_1) \rho dy + R_-'\right|
&\le& C\|q\|_{\H}^3+C\sum_{i=1}^k \frac{|d_i'|}{1-d_i^2}J_i\|q\|_{\H}^2
\le C\left(\|q\|_{\H}^3+ J\|q\|_{\H}^2\right).\label{zz4}
\end{eqnarray}
Note that from \aref{boundF}, the H\"older inequality and Claim \ref{lemsobolev-9}, we have
\begin{eqnarray}
\left|\iint \F(q_1) \rho dy\right| &\le& 
C\iint |q_1|^{p+1} \rho dy + C \delta_{\{p\ge 2\}} \iint |K|^{p-2}|q_1|^3\rho dy\nonumber\\
&\le &C\|q\|_{\H}^{p+1} + C \delta_{\{p\ge 2\}}  \left(\iint |q_1|^{p+1}\rho dy\right)^{\frac 3{p+1}}\left(\iint |K|^{p+1} \rho dy\right)^{\frac{p-2}{p+1}}\nonumber\\
&\le& C\|q\|_{\H}^{p+1} + C\delta_{\{p\ge 2\}} \|q\|_{\H}^3\le C\|q\|_{\H}^{\bar p+1}\label{zz5}
\end{eqnarray}
where  $\p=\min(p,2)$.

\bigskip

- Using the Cauchy-Schwartz inequality, we write
\begin{equation}\label{zz6}
\left|\iint q_{-,2} G \rho dy\right|
\le \frac 1{p-1}  \iint q_{-,2}^2 \frac \rho{1-y^2} dy + 
C\iint G^2 \rho (1-y^2) dy.
\end{equation}
From the definition \aref{defG} of $G$, we need to handle $R$ and $V_i F_{1,1}$. We start by $R$ first.\\
We claim that 
\begin{equation}\label{krim0}
|R|\le  C \sum_{j=1}^k\kappa(d_j(s),y)^{p-1}\un{y_{j-1}(s)}{y_{j}(s)} \sum_{l\neq j} \kappa(d_l(s),y)
\end{equation}
where  
\begin{equation}\label{defy}
y_0=-1,\;\; y_j= \tanh\left(\frac{\zeta_j+\zeta_{j+1}}2\right)\mbox{ if }j=1,..,k-1\mbox{ and }y_k=1.
\end{equation}
 In particular, we have
\[
-1=y_0< -d_1<y_1<-d_2<...<y_j<-d_j<y_{j+1}<...<-d_k<y_k=1
\]
and $\kappa(d_j(s),y_{j+1}(s))=\kappa(d_{j+1}(s),y_{j+1}(s))$ for $j=1,...,k-1$ (to see this, just use the fact that $\kappa(d,y)(1-y^2)^{\frac 1{p-1}}= \kappa_0\cosh^{-\frac 2{p-1}}(\xi-\zeta_i)$ if $y=\tanh \xi$).\\
To prove \aref{krim0}, we take $y\in (y_{j-1}(s), y_j(s))$ and set $X = (\d\sum_{l\neq j}e_l \kappa(d_l(s),y))/e_j \kappa(d_j(s),y)$. From the fact that $\zeta_{j+1}(s)-\zeta_j(s)\to \infty$, we have $|X|\le 2$ hence
\[
||1+X|^{p-1}(1+X) -1|\le C|X|
\]
and for $y\in (y_{j-1}(s), y_j(s))$ and $s$ large,
\[
||K|^{p-1}K - e_j \kappa(d_j(s),y)^p|\le C \kappa(d_j(s),y)^{p-1} \sum_{l\neq j} \kappa(d_l(s),y).
\]
Since for all $y\in (y_{j-1}(s), y_j(s))$, $\kappa(d_j(s),y)\ge \kappa(d_l(s),y)$ if $l\neq j$, this concludes the proof of \aref{krim0}.\\
Using \aref{krim0}, we see that
\begin{equation}\label{zz7}
\iint R^2\rho (1-y^2) dy \le C\sum_{j=1}^k\sum_{l\neq j}\int_{y_{j-1}}^{y_{j}} \kappa(d_j)^{2(p-1)} \kappa(d_l)^2 \rho (1-y^2) dy\le C \sum_{m=1}^{k-1} h(\zeta_{m+1}-\zeta_m)^2
\end{equation} 
where $h$ is defined in \aref{defh}.\\ 
%
%
%
Now, we handle $V_iF_{1,1}$. 
Using \aref{krim1}, \aref{boundF12} and (i) of Lemma \ref{lemtech}, we see that
\begin{eqnarray}
&&\iint \left(V_i F_{1,1}(d_i)\right)^2 \rho (1-y^2) dy
\le C \sum_{j\neq i}\iint  \kappa(d_i)^2\kappa(d_j)^{2(p-1)}\rho(1-y^2) dy \nonumber\\
&+& C \delta_{\{p \ge 2\}} 
\sum_{j\neq i}\iint \kappa(d_i)^{2(p-1)} \kappa(d_j)^2 \rho (1-y^2) dy\to 0\mbox{ as }s\to \infty.\nonumber
\end{eqnarray}
Hence, using \aref{equivq}, we see that
\begin{equation}\label{zz8}
\left(\alpha^i_1\right)^2\iint \left(V_i F_{1,1}(d_i)\right)^2 \rho (1-y^2) dy 
=o\left(\|q\|_{\H}^2\right).
\end{equation}
%
%
%
Gathering \aref{zz1}, \aref{zz2}, \aref{zz3}, \aref{zz4}, \aref{zz6}, \aref{defG}, \aref{zz7} and \aref{zz8}, we get to the conclusion of \aref{eqam}. Note that the estimate for $R_-(s)$ is given in \aref{zz5}. 

\bigskip

{\bf Proof of (iii): An additional estimate}

We prove estimate \aref{153} here. The proof is the same as in the case of one soliton treated in \cite{MZjfa07}, except for the term involving the interaction term $R(y,s)$ \aref{eqq}. Therefore, arguing exactly as in pages 110 and 112 of \cite{MZjfa07}, we write
\begin{eqnarray*}
\frac d{ds} \int q_1q_2 \rho dy 
\le -\frac 9{10}A_-+ CJ^2+ C\iint q_{-,2}^2 \frac{\rho}{1-y^2} dy +C  \sum_{i=1}^k |\alpha_1^i|^2 + \iint q_1R \rho dy.
\end{eqnarray*}
Since we have from the Cauchy-Schwartz inequality, (i) of Claim \ref{lemsobolev-9}, \aref{qdecomp} and \aref{zz7}
\begin{eqnarray*}
&&\left|\iint q_1R \rho dy\right| \le \left(\iint q_1^2 \frac \rho{1-y^2} dy\right)^{\frac 12}\left(\iint R^2 (1-y^2) \rho dy\right)^{\frac 12}\\
&\le & C\|q\|_{\H}\left(\iint R^2 (1-y^2) \rho dy\right)^{\frac 12}
\le  \frac 1{10}\left(A_- + \sum_{i=1}^k (\alpha_1^i)^2\right) + C \sum_{i=1}^{k-1} h(\zeta_{i+1}- \zeta_i)^2
\end{eqnarray*}
where $h$ is defined in \aref{defh}, this concludes the proof of \aref{153} and the proof of Lemma \ref{lemproj}. \Box
\section{A continuity result in the selfsimilar variable}\label{appclrakia}
We prove Claim \ref{clrakia} here.
%
%
%
%
Consider $\epsilon_0>0$ and from \aref{carprofile01}, fix $\tilde t$ close enough to $T(x_0)$ so that 
\begin{equation}\label{wbg}
\|w_{x_0}(s_0)-w_\infty\|_{L^2_\rho}\le \epsilon_0\mbox{ where }s_0=-\log(T(x_0)-\tilde t).
\end{equation}
Note from \aref{espace0} and the continuity of $x\mapsto T(x)$ that $u(x,\tilde t)$ is well defined for all $x\in [\bar x, x_0+(T(x_0)-\tilde t))$ for some $\bar x<x_0-(T(x_0)-\tilde t)$. Therefore, using the selfsimilar transformation \aref{defw}, we see that
\begin{equation}\label{defym}
w(\cdot, s_0) \in L^2(\bar y, 0)\mbox{ where }\bar y = \frac{\bar x-x_0}{T(x_0)-\tilde t}<-1.
\end{equation}
We aim at proving that for $x'$ close enough to $x_0$, we have
\begin{equation}\label{deftsm}
\left\|\vc{w_{x'}(\tilde s_0(x'))}{\partial_s w_{x'}(\tilde s_0(x'))} - \vc{w_\infty}0\right\|_\H \le 6\epsilon_0\mbox{ where }\tilde s_0(x')= - \log(T(x')-\tilde t).
\end{equation}
For simplicity, we will only prove that 
\begin{equation}\label{fayssal}
\|w_{x'}(\tilde s_0(x'))-w_\infty\|_{L^2_\rho}\le 2 \epsilon_0,
\end{equation}
provided that $x_0-x'$ is small. The estimates involving $\partial_y w_{x'}(\tilde s_0(x'))$ and $\ps w_{x'}(\tilde s_0(x'))$ follow in the same way.\\
Using the selfsimilar transformation \aref{eqw}, we write 
\begin{equation}\label{wawb} 
\forall \tilde y\in (-1,1),\;\;w_{x'}(\tilde y,\tilde s_0(x'))= \theta^{\frac 2{p-1}}w_{x_0}(y,s_0)\mbox{ where }y= \tilde y\theta + \xi,
\end{equation}
\begin{equation}\label{deftx}
\theta =\frac 1{1+ e^{\tilde s_0(x')}(T(x_0)-T(x'))}\to 1\mbox{ and }\xi = (x'-x_0)e^{\tilde s_0(x')}\theta\to 0\mbox{ as }x' \to x_0. 
\end{equation}
Therefore, performing a change of variables, we write for $x_0-x'$ small enough,
\begin{eqnarray*}
&&\|w_{x'}(\tilde s_0(x'))-w_\infty\|_{L^2_\rho}^2=\iint |w_{x'}(\tilde y, \tilde s_0(x'))-w_\infty(\tilde y)|^2 \rho(\tilde y)d\tilde y\\
&= & \int_{-\theta+\xi}^{\theta+\xi}\left|\theta^{\frac 2{p-1}}w_{x_0}(y,s_0)-w_\infty\left(\frac{y-\xi}\theta\right)\right|^2\rho\left(\frac{y-\xi}{\theta}\right)\frac{dy}\theta.
\end{eqnarray*}
Since we have from \aref{deftx}, \aref{deftsm} and the fact that $x\mapsto T(x)$ is $1$-Lipschitz, 
\[
\theta+\xi = \frac{1+(x'-x_0)e^{\tilde s_0(x')}}{1+ e^{\tilde s_0(x')}(T(x_0)-T(x'))}
= \frac{T(x') - \tilde t+x'-x_0}{T(x_0)-\tilde t}\le 1,
\]
it follows that
\[
\|w_{x'}(\tilde s)-w_\infty\|_{L^2_\rho}^2=\int_{\bar y}^1g(\theta, \xi, y)dy
\]
where $\bar y<-1$ is defined in \aref{defym} and
\begin{equation}\label{defg0}
g(\theta, \xi, y) = \frac{1_{\{-\theta+\xi<y<\theta+\xi\}}}\theta\left|\theta^{\frac 2{p-1}}w_{x_0}(y,s_0)-w_\infty\left(\frac{y-\xi}\theta\right)\right|^2\rho\left(\frac{y-\xi}{\theta}\right).
\end{equation}
We claim that in order to conclude, it is enough to prove that for $x_0-x'$ small enough, 
\begin{equation}\label{deftg}
\forall y\in (-\theta+\xi,\theta+\xi),\; g(\theta, \xi, y) \le \tilde g(y)\mbox{ for some }\tilde g\in L^1(\bar y,1).
\end{equation}
Indeed, since we have from \aref{deftx} that 
\[
\forall y\in (\bar y,1),\;\;g(\theta, \xi, y) \to g(1,0,y)\mbox{ as }x'\to x_0,
\]
we use \aref{deftg} to apply the Lebesgue Theorem and obtain that
\[
\|w_{x'}(\tilde s_0(x'))-w_\infty\|_{L^2_\rho}^2=\int_{\bar y}^1g(\theta, \xi, y)dy\to \iint g(1,0,y) dy=\|w_{x_0}(s_0)-w_\infty\|_{L^2_\rho}^2
\]
as $x'\to x_0$. Using \aref{wbg}, we see that for $x_0-x'$ small enough, \aref{fayssal} holds. It remains to prove \aref{deftg} in order to conclude.\\
If $-\theta+\xi\le y \le 0$, then we have $\rho(\frac{y-\xi}\theta)\le 1$. Using \aref{defg0}, \aref{deftx}, \aref{defym} and the definition \aref{defg} of $w_\infty$, we write for $x_0-x'$ small enough: 
\[
g(\theta, \xi, y) \le C(|w_{x_0}(y,s_0)|^2+\|w_\infty\|_{L^\infty(-1,1)}^2) \in L^1(\bar y,0).
\]
If $0\le y\le \theta+\xi$, then we have from \aref{deftx}, $\rho(\frac{y-\xi}\theta)\le C(1-\frac{y-\xi}\theta)^{\frac 2{p-1}}=C(\frac{1-y+\xi}\theta)^{\frac 2{p-1}}\le C(1-y)^{\frac 2{p-1}}\le C\rho(y)$. Therefore, using \aref{defg0} and \aref{wbg}, we write
\[
g(\theta, \xi, y) \le C(|w_{x_0}(y,s_0)|^2+\|w_\infty\|_{L^\infty(-1,1)}^2)\rho(y)\in L^1(0, 1).
\]
Thus, \aref{deftg} holds and so does \aref{fayssal}.\\ 
Since the same technique works for $\left\|\partial_y w_{x'}(\tilde s_0(x'))-\frac{dw_\infty}{dy}\right\|_{L^2_{\rho(1-y^2)}}$ and $\|\partial_s w_{x'}(\tilde s_0(x'))\|_{L^2_\rho}$, estimate \aref{deftsm} follows in the same way. This concludes the proof of Claim \ref{clrakia}. \Box

\section{Computations in the $\xi$ variable}\label{apptech}
In the following, we compute  integrals involving the solitons $\kappa(d,y)$ \aref{defkd}.\\ 
Recalling that $y=-d_i(s)=\tanh \zeta_i(s)$ is the center of the $i$-th soliton $\kappa(d_i(s),y)$, we introduce the following ``separators'' between the solitons: 
\begin{equation}\label{0defy}
y_0=-1,\;\; y_j= \tanh\left(\frac{\zeta_j+\zeta_{j+1}}2\right)\mbox{ if }j=1,..,k-1\mbox{ and }y_k=1
\end{equation}
Note in particular that we have
\[
-1=y_0< -d_1<y_1<-d_2<...<y_j<-d_j<y_{j+1}<...<-d_k<y_k=1
\]
and $\kappa(d_j(s),y_{j+1}(s))=\kappa(d_{j+1}(s),y_{j+1}(s))$ for $j=1,...,k-1$ (to see this, just use the fact that $\kappa(d,y)(1-y^2)^{\frac 1{p-1}}= \kappa_0\cosh^{-\frac 2{p-1}}(\xi-\zeta_i)$ if $y=\tanh \xi$).\\
In the following lemma, we estimate various integrals involving the solitons $\kappa(d_i(s),y)$:
\begin{lem}\label{lemtech}{\bf (A table of integrals involving the solitons)}
We have the following estimates as $s\to \infty$:

(i) If $i\neq j$, $\alpha>0$, $\beta>0$ and $I_1= \d\iint \kappa(d_j)^\alpha\kappa(d_i)^\beta(1-y^2)^{\frac{\alpha+\beta}{p-1}-1} dy$, then:\\
for $\alpha=\beta$, $I_1\sim C_0|\zeta_i-\zeta_j|e^{-\frac {2\beta}{p-1}|\zeta_i-\zeta_j|}$;\\
for $\alpha\neq \beta$, $I_1\sim C_0e^{-\frac 2{p-1}\min(\alpha, \beta)|\zeta_i-\zeta_j|}$  for some $C_0=C_0(\alpha, \beta)>0$.

(ii) If $i\neq j$, $\alpha>0$, $\beta>0$ and $I_2\equiv\d\int_{y_{j-1}}^{y_j}\kappa(d_j)^\alpha\kappa(d_i)^\beta(1-y^2)^{\frac{\alpha+\beta}{p-1}-1} dy$, then:\\
for $\alpha =\beta$, 
$I_2 \le C |\zeta_{j+1}-\zeta_j|e^{-\frac{2\beta}{p-1}|\zeta_{j+1}-\zeta_j|}
+C |\zeta_{j-1}-\zeta_j|e^{-\frac{2\beta}{p-1}|\zeta_{j-1}-\zeta_j|}$;\\
for $\alpha >\beta$, $I_2 \le C e^{-\frac{2\beta}{p-1}|\zeta_{j+1}-\zeta_j|}
+C e^{-\frac{2\beta}{p-1}|\zeta_{j-1}-\zeta_j|}$;\\
for $\beta>\alpha$, $I_2 \le C e^{-\frac{(\alpha+\beta)}{p-1}|\zeta_{j+1}-\zeta_j|}
+C e^{-\frac{(\alpha+\beta)}{p-1}|\zeta_{j-1}-\zeta_j|}$.

(iii) Let $A_{i,j,l} = \d\int_{y_{j-1}}^{y_j}\frac{y+d_i}{1+yd_i}\kappa(d_i)\kappa(d_j)^{p-1}\kappa(d_l)\rho dy$ with $l\neq j$. Then, for some $c_1'''>0$ and $\delta_5(p)>0$, we have:\\
-if $i=j$ and $l=i\pm 1$, then $|A_{i,i,l}- \sgn(l-j)c_1''' e^{-\frac 2{p-1}|\zeta_l-\zeta_i|}|\le  CJ^{1+\delta_5}$,\\ 
 - otherwise, $A_{i,j,l} \le CJ^{1+\delta_5}$, where $J$ is defined in \aref{defJ}.

(iv) If $l\neq j$, then  $B_{i,j,l} \equiv \d\int_{y_{j-1}}^{y_j}\kappa(d_i)\kappa(d_j)^{p-\p}\kappa(d_l)^\p\rho dy\le C J^{1+\delta_6}$ for some $\delta_6(p)>0$ (with $\bar p =\min(p,2)$).

(v) For any $i=1,..,k$, it holds that $J_i\equiv\iint \kappa(d_i)|K|^{p-2} dy \le C$ where $K(y,s)$ is defined in \aref{eqq}.
\end{lem}
{\it Proof}:
(i) With the change of variables  $y=\tanh \xi$, we write
\[
I_1= \kappa_0^{\alpha+\beta}\int_\R \cosh^{-\frac{2\alpha}{p-1}}(\xi-\zeta_j)\cosh^{-\frac{2\beta}{p-1}}(\xi-\zeta_i)d\xi.
\] 
From symmetry, we can assume that $\alpha \ge \beta$ and $\zeta_i> \zeta_j$. Using the change of variables $z=\xi-\zeta_j$, we write
\[
I_1= \kappa_0^{\alpha+\beta}\int_\R \cosh^{-\frac{2\alpha}{p-1}}(z)\cosh^{-\frac{2\beta}{p-1}}(z+\zeta_j-\zeta_i)dz.
\]
When $\alpha >\beta$, we get from Lebesgue's Theorem $I_1 \sim C e^{-\frac {2\beta}{p-1}(\zeta_i-\zeta_j)}$.\\
When $\alpha=\beta$, we write from symmetry and Lebesgue's Theorem
\[
I_1=2\kappa_0^{\alpha+\beta}\int_{-\infty}^{\frac{\zeta_i-\zeta_j}2} \cosh^{-\frac{2\beta}{p-1}}(z)\cosh^{-\frac{2\beta}{p-1}}(z+\zeta_j-\zeta_i)dz\sim C(\zeta_i-\zeta_j)e^{-\frac {2\beta}{p-1}(\zeta_i-\zeta_j)}.
\]

(ii) Since $I_2\le I_1$ and $|\zeta_j-\zeta_i|\ge \min(|\zeta_{j+1}-\zeta_j|, |\zeta_j-\zeta_{j-1}|)$, the result follows from (i) if $\alpha \ge \beta$. When $\alpha <\beta$, we assume that $\zeta_i>\zeta_j$, the other case being parallel. Using the change of variables $y=\tanh \xi$ then $z=\xi-\zeta_j$, we write
\begin{eqnarray*}
I_2&=&\kappa_0^{\alpha+\beta} \int_{-\frac{(\zeta_j-\zeta_{j-1})}2}^{\frac{(\zeta_{j+1}-\zeta_j)}2}\cosh^{-\frac{2\alpha}{p-1}}(z) \cosh^{-\frac {2\beta}{p-1}}(z+\zeta_j-\zeta_i)dz\\
&\sim& e^{-\frac{2\beta}{p-1}(\zeta_i-\zeta_j)}\int_{-\infty}^{\frac{(\zeta_{j+1}-\zeta_j)}2}\cosh^{-\frac{2\alpha}{p-1}}(z) e^{\frac{2\beta}{p-1}z}dz\\
&\sim & C e^{-\frac{2\beta}{p-1}(\zeta_i-\zeta_j)} e^{\frac{2(\beta-\alpha)}{p-1}.\frac{\zeta_{j+1}-\zeta_j}2}\le C e^{-\frac{(\alpha+\beta)}{p-1}(\zeta_{j+1}-\zeta_j)} 
\end{eqnarray*}
since $\zeta_i-\zeta_j\ge \zeta_{j+1} - \zeta_j$, which yields the result. 

\medskip

(iii) If $i=j$, we assume that $\zeta_l >\zeta_i$, since the other case follows by replacing $\xi$ by $-\xi$ (generating a minus sign in the formula). Therefore, it holds that
\begin{equation}\label{linatou}
\zeta_l\ge \zeta_{i+1}.
\end{equation}
Using the change of variables $y=\tanh \xi$, we write
\begin{eqnarray}
A_{i,i,l}&=&\kappa_0^{p+1}\int_{\frac{\zeta_{i-1}+\zeta_i}2}^{\frac{\zeta_i+\zeta_{i+1}}2} \cosh^{-\frac {2p}{p-1}}(\xi-\zeta_i)\tanh (\xi-\zeta_i) \cosh^{-\frac 2{p-1}}(\xi-\zeta_l)d\xi\nonumber\\
&=&\kappa_0^{p+1}\int_{-\frac{(\zeta_i-\zeta_{i-1})}2}^{\frac{(\zeta_{i+1}-\zeta_i)}2}\cosh^{-\frac {2p}{p-1}}(z)\tanh (z) \cosh^{-\frac 2{p-1}}(z+\zeta_i-\zeta_l)dz.\label{asf}
\end{eqnarray}
Since we see from \aref{linatou} that when $z\le \frac{(\zeta_{i+1}-\zeta_i)}2$, it holds that $z+\zeta_i-\zeta_l\le \frac{(\zeta_{i+1}-\zeta_i)}2 +\zeta_i - \zeta_{i+1}=- \frac{(\zeta_{i+1}-\zeta_i)}2\to -\infty$ as $s\to \infty$, we deduce that
\[
\left|\cosh^{-\frac 2{p-1}}(z+\zeta_i-\zeta_l)- 2^{\frac 2{p-1}}e^{\frac {2(z+\zeta_i-\zeta_l)}{p-1}}\right|\le Ce^{\frac {2(z+\zeta_i-\zeta_l)}{p-1}}e^{(\zeta_{i+1}-\zeta_i)}.
\]
Therefore, using \aref{asf} and the definition \aref{defJ} of $J$, we see that 
\[
|A_{i,i,l}-c_1"e^{-\frac 2{p-1}(\zeta_l-\zeta_i)}|\le Ce^{-\frac 2{p-1}(\zeta_l-\zeta_i)}e^{(\zeta_{i+1}-\zeta_i)}\le C J^{1+\frac{p-1}2}
\]
where
\begin{eqnarray*}
c_1" &=& 2^{\frac 2{p-1}}\kappa_0^{p+1} \int_\R\cosh^{-\frac {2p}{p-1}}(z)\tanh (z) e^{\frac {2z}{p-1}}dz\\
&=&2^{\frac 2{p-1}}\kappa_0^{p+1} \int_0^{\infty}\cosh^{-\frac {2p}{p-1}}(z)\tanh (z)(e^{\frac {2z}{p-1}}-e^{-\frac {2z}{p-1}}) dz>0,
\end{eqnarray*}
which gives the result when $l=i+1$.\\
If $l\ge i+2$, then $e^{-\frac 2{p-1}(\zeta_l-\zeta_i)}\le e^{-\frac 2{p-1}(\zeta_l-\zeta_{i+1)}}
e^{-\frac 2{p-1}(\zeta_{i+1}-\zeta_i)}\le J^2$ and the result follows as well.

\medskip

Now, if $j\neq i$, then we have from the Cauchy-Schwartz inequality, 
\[
A_{i,j,l}\le \left(\int_{y_{j-1}}^{y_j}\kappa(d_j)^{p-1}\kappa(d_i)^2\rho dy\right)^{1/2}
\left(\int_{y_{j-1}}^{y_j}\kappa(d_j)^{p-1}\kappa(d_l)^2\rho dy\right)^{1/2},
\]
and the conclusion follows from (ii) and the definition \aref{defJ} of $J(s)$.

(iv) If $i=j$, then the result follows from (ii).
If $i\neq j$, using the H\"older inequality with $P=\p+1$ and $Q=\frac{\p+1}{\p}$, we write
\[
B_{i,j,l}\le \left(\int_{y_{j-1}}^{y_j}\kappa(d_i)^{\p+1}\kappa(d_j)^{p-\p}\rho dy\right)^{\frac 1{\p+1}}
\left(\int_{y_{j-1}}^{y_j}\kappa(d_l)^{\p+1}\kappa(d_j)^{p-\p}\rho dy\right)^{\frac {\p}{\p+1}},
\]
and the result follows from (ii) and the definition \aref{defJ} of $J(s)$.
\bigskip

(v) Using the change of variables $y=\tanh \xi$, we write 
\begin{equation}\label{defkb}
J_i = \kappa_0^{p-1}\int_{\R} \cosh^{-\frac 2{p-1}}(\xi-\zeta_i)|\bar K(\xi,s)|^{p-2} d\xi
\mbox{ where } \bar K(\xi,s) = \sum_{j=1}^k e_j \cosh^{-\frac 2{p-1}}(\xi-\zeta_j).
\end{equation}
If $p \ge 2$, then $|\bar K(\xi,s)| \le C$ and $|J_i(s)| \le C$.\\
If $p <2$ and the $e_j$ are the same, then $|\bar K(\xi,s)| \ge \cosh^{-\frac 2{p-1}}(\xi-\zeta_i)$ and $|J_i(s)|\le \int_{\R} \cosh^{-2}(\xi-\zeta_i)d \xi \le C$.\\ 
It remains to treat the delicate case where $p<2$ with the $e_j$ not all the same. Taking advantage of the decoupling in the sum of the solitons (see \aref{decouple}), we write
\begin{equation}\label{b0}
J_i = \kappa_0^{p-1}\sum_{j=1}^k \int_{\theta_{j-1}+A}^{\theta_j+A} \cosh^{-\frac 2{p-1}}(\xi-\zeta_i)|\bar K(\xi,s)|^{p-2} d\xi
\end{equation}
where $\theta_0 = - \infty$, $\theta_j = \frac{\zeta_j+\zeta_{j+1}}2$  if $j=1,..,k-1$, $\theta_k=\infty$ and $A=A(p)$ is fixed such that
\begin{equation}\label{mehdi}
e^{\frac {2A}{p-1}} \ge 2 e^{-\frac {2A}{p-1}}.
\end{equation} 
This partition isolates each soliton in the definition of $\bar K(\xi,s)$. It is shifted by $A$ since $\bar K(\xi,s)$ may be zero for some $z_j(s) \sim \theta_j(s)$ if $e_j e_{j+1}=-1$, giving rise to a singularity in $|\bar K(\xi,s)|^{p-2}$, integrable though delicate to control.\\
Consider some $j=1,...,k-1$.\\
If $e_j= e_{j+1}$, then we have from \aref{decouple} and \aref{mehdi} for all $\xi \in (\theta_{j-1}+A, \theta_j+A)$, $|\bar K(\xi,s)|\ge C(A)\cosh^{-\frac 2{p-1}}(\xi- \zeta_j)$ and $\cosh^{-\frac 2{p-1}}(\xi-\zeta_i) \le C(A) \cosh^{-\frac 2{p-1}}(\xi-\zeta_j)$, hence
\begin{equation}\label{b1}
 \int_{\theta_{j-1}+A}^{\theta_j+A} \cosh^{-\frac 2{p-1}}(\xi-\zeta_i)|\bar K(\xi,s)|^{p-2} d\xi \le C(A)\int_{\theta_{j-1}+A}^{\theta_j+A} \cosh^{-2} (\xi-\zeta_j) d\xi \le C(A).
\end{equation}
If $e_j = -e_{j+1}$, then  $\bar K(z_j(s), s)=0$ with $z_j(s) \sim \theta_j(s)$, which makes $|\bar K(\xi,s)|^{p-2}$ singular at $\xi=z_j(s)$. 
For this we split the integral over the interval $(\theta_{j-1}+A, \theta_j+A)$ into two parts, below and above $\theta_j-A$:\\
- the part on the interval $(\theta_{j-1}+A, \theta_j-A)$ is bounded by the same argument as in the case $e_j = e_{j+1}$;\\
- the part on the interval $(\theta_j-A, \theta_j+A)$. Since we have from the definition \aref{defkb} of $\bar K(\xi,s)$
\[
\partial_\xi \bar K(\xi,s) = -\frac 2{p-1}\sum_{l=1}^k e_l\sinh(\xi-\zeta_l) \cosh^{-\frac 2{p-1}-1}(\xi- \zeta_l),
\]
it follows that for all $\xi \in (\theta_j-A, \theta_j+A)$, $|\partial_\xi \bar K(\xi,s)|\ge C(A) e^{-\frac 2{p-1}(\theta_j-\zeta_j)}= C(A) e^{-\frac{\zeta_{j+1}-\zeta_j}{p-1}}$ for some $C(A)>0$, hence
\[
|\bar K(\xi,s)|^{p-2} = |\bar K(\xi,s)-\bar K(z_j(s),s)|^{p-2} \le C(A)|\xi-z_j(s)|^{p-2} e^{-\frac{p-2}{p-1}(\zeta_{j+1}-\zeta_j)}.
\]
 Therefore, since for all $\xi \in (\theta_j-A, \theta_j+A)$, $\cosh^{-\frac 2{p-1}}(\xi-\zeta_i) \le C(A)\cosh^{-\frac 2{p-1}}(\xi-\zeta_j) \le C(A) \cosh^{-\frac 2{p-1}}(\theta_j-\zeta_j)\le C(A) e^{-\frac{\zeta_{j+1} - \zeta_j}{p-1}}$, it follows that
\begin{eqnarray}
\int_{\theta_j-A}^{\theta_j+A} \cosh^{-\frac 2{p-1}}(\xi-\zeta_i)|\bar K(\xi,s)|^{p-2} d\xi &\le& C(A) e^{-(\zeta_{j+1}- \zeta_j)}\int_{\theta_j-A}^{\theta_j+A}|\xi-z_j(s)|^{p-2} d\xi\nonumber\\
& \le&  C(A) e^{-(\zeta_{j+1}- \zeta_j)}\label{b2}
\end{eqnarray}
because $z_j(s) \sim \theta_j(s)$ as $s\to \infty$. Therefore, (v) follows from \aref{b0}, \aref{b1} and \aref{b2}. \Box


\begin{thebibliography}{10}

\bibitem{Apndeta95}
S.~Alinhac.
\newblock {\em Blowup for nonlinear hyperbolic equations}, volume~17 of {\em
  Progress in Nonlinear Differential Equations and their Applications}.
\newblock Birkh\"auser Boston Inc., Boston, MA, 1995.

\bibitem{Afle02}
S.~Alinhac.
\newblock A minicourse on global existence and blowup of classical solutions to
  multidimensional quasilinear wave equations.
\newblock In {\em Journ\'ees ``\'Equations aux D\'eriv\'ees Partielles''
  (Forges-les-Eaux, 2002)}, pages Exp. No. I, 33. Univ. Nantes, Nantes, 2002.

\bibitem{AMimrn01}
C.~Antonini and F.~Merle.
\newblock Optimal bounds on positive blow-up solutions for a semilinear wave
  equation.
\newblock {\em Internat. Math. Res. Notices}, (21):1141--1167, 2001.

\bibitem{CFarma85}
L.~A. Caffarelli and A.~Friedman.
\newblock Differentiability of the blow-up curve for one-dimensional nonlinear
  wave equations.
\newblock {\em Arch. Rational Mech. Anal.}, 91(1):83--98, 1985.

\bibitem{CFtams86}
L.~A. Caffarelli and A.~Friedman.
\newblock The blow-up boundary for nonlinear wave equations.
\newblock {\em Trans. Amer. Math. Soc.}, 297(1):223--241, 1986.

\bibitem{GSVjfa92}
J.~Ginibre, A.~Soffer, and G.~Velo.
\newblock The global {C}auchy problem for the critical nonlinear wave equation.
\newblock {\em J. Funct. Anal.}, 110(1):96--130, 1992.

\bibitem{KL1cpde93}
S.~Kichenassamy and W.~Littman.
\newblock Blow-up surfaces for nonlinear wave equations. {I}.
\newblock {\em Comm. Partial Differential Equations}, 18(3-4):431--452, 1993.

\bibitem{KL2cpde93}
S.~Kichenassamy and W.~Littman.
\newblock Blow-up surfaces for nonlinear wave equations. {I}{I}.
\newblock {\em Comm. Partial Differential Equations}, 18(11):1869--1899, 1993.

\bibitem{Ltams74}
H.~A. Levine.
\newblock Instability and nonexistence of global solutions to nonlinear wave
  equations of the form {$Pu\sb{tt}=-Au+{\cal F}(u)$}.
\newblock {\em Trans. Amer. Math. Soc.}, 192:1--21, 1974.

\bibitem{MZajm03}
F.~Merle and H.~Zaag.
\newblock Determination of the blow-up rate for the semilinear wave equation.
\newblock {\em Amer. J. Math.}, 125:1147--1164, 2003.

\bibitem{MZimrn05}
F.~Merle and H.~Zaag.
\newblock Blow-up rate near the blow-up surface for semilinear wave equations.
\newblock {\em Internat. Math. Res. Notices}, (19):1127--1156, 2005.

\bibitem{MZjfa07}
F.~Merle and H.~Zaag.
\newblock Existence and universality of the blow-up profile for the semilinear
  wave equation in one space dimension.
\newblock {\em J. Funct. Anal.}, 253(1):43--121, 2007.

\bibitem{MZcmp08}
F.~Merle and H.~Zaag.
\newblock Openness of the set of non characteristic points and regularity of
  the blow-up curve for the $1$ d semilinear wave equation.
\newblock {\em Comm. Math. Phys.}, 282:55--86, 2008.

\bibitem{MZisol09}
F.~Merle and H.~Zaag.
\newblock Isolatedness of characteristic points for a semilinear wave equation
  in one space dimension.
\newblock 2009.
\newblock in preparation.

\end{thebibliography}
\def\cprime{$'$}

\noindent{\bf Address}:\\
Universit\'e de Cergy Pontoise, D\'epartement de math\'ematiques, 
2 avenue Adolphe Chauvin, BP 222, 95302 Cergy Pontoise cedex, France.\\
\vspace{-7mm}
\begin{verbatim}
e-mail: merle@math.u-cergy.fr
\end{verbatim}
Universit\'e Paris 13, Institut Galil\'ee, 
Laboratoire Analyse, G\'eom\'etrie et Applications, CNRS UMR 7539,
99 avenue J.B. Cl\'ement, 93430 Villetaneuse, France.\\
\vspace{-7mm}
\begin{verbatim}
e-mail: Hatem.Zaag@univ-paris13.fr
\end{verbatim}

\end{document}